\documentclass[twoside]{amsart}
\usepackage{amsmath,amsfonts,amsthm,mathrsfs}
\usepackage[unicode]{hyperref}

\usepackage[alphabetic,initials]{amsrefs}

\usepackage{amssymb}
\usepackage{enumitem}

\usepackage[dpi=1270,PostScript=dvips,heads=LaTeX]{diagrams}
\diagramstyle{small,labelstyle=\scriptstyle}
\newcommand{\rToDots}{\rTo[l>=3.3em]~{\raisebox{0.1pt}{$\,\dotsb\,$}}}
\newarrow{Equal}   =====
\newarrow{IntoBold}   {boldhook}---> 

\usepackage{chngcntr}

\makeatletter
\let\@@seccntformat\@seccntformat
\renewcommand*{\@seccntformat}[1]{%
  \expandafter\ifx\csname @seccntformat@#1\endcsname\relax
    \expandafter\@@seccntformat
  \else
    \expandafter
      \csname @seccntformat@#1\expandafter\endcsname
  \fi
    {#1}%
}
\newcommand*{\@seccntformat@subsection}[1]{%
  \textbf{\csname the#1\endcsname.}
}
\makeatother

\makeatletter
\let\@paragraph\paragraph
\renewcommand*{\paragraph}[1]{%
	\vspace{0.3\baselineskip}%
	\@paragraph{\textit{#1}}%
}
\makeatother

\counterwithin{equation}{subsection}
\counterwithout{subsection}{section}
\counterwithin{figure}{subsection}

\newcommand{\subsecref}[1]{\S\ref{#1}}

\newtheorem{theorem}[equation]{Theorem}
\newtheorem*{theorem*}{Theorem}
\newtheorem{lemma}[equation]{Lemma}
\newtheorem{corollary}[equation]{Corollary}
\newtheorem{proposition}[equation]{Proposition}
\newtheorem{conjecture}[equation]{Conjecture}
\newtheorem{condition}[equation]{Condition}
\theoremstyle{definition}
\newtheorem{definition}[equation]{Definition}
\theoremstyle{remark}

\newtheorem*{note}{Note}
\newtheorem*{notation}{Notation}
\theoremstyle{plain}

\makeatletter
\let\old@caption\caption
\renewcommand*{\caption}[1]{%
	\setcounter{figure}{\value{equation}}%
	\stepcounter{equation}%
	\old@caption{#1}\relax%
}
\makeatother

\newcommand{\MHM}{\operatorname{MHM}}
\newcommand{\MHS}{\operatorname{MHS}}

\newcommand{\Dmod}{\mathcal{D}}
\newcommand{\Mmod}{\mathcal{M}}
\newcommand{\Nmod}{\mathcal{N}}

\newcommand{\decal}[1]{\lbrack #1 \rbrack}

\newcommand{\ltriangle}[4][]%
{\begin{diagram}[#1]%
	{#2} &\rTo& {#3} &\rTo& {#4} &\rTo& {#2 \decal{1}}%
\end{diagram}}

\newcommand{\shH}{\shf{H}}

\newcommand{\norm}[1]{\lVert#1\rVert}
\newcommand{\bignorm}[1]{\bigl\lVert#1\bigr\rVert}
\newcommand{\abs}[1]{\lvert #1 \rvert}
\newcommand{\bigabs}[1]{\bigl\lvert #1 \bigr\rvert}

\newcommand{\eps}{\varepsilon}

\newcommand{\tensor}{\otimes}

\newcommand{\End}[1]{\operatorname{End}(#1)}
\newcommand{\shHom}{\mathscr{H}\hspace{-2.7pt}\mathit{om}}
\newcommand{\shExt}{\mathscr{E}\hspace{-1.5pt}\mathit{xt}}
\newcommand{\Hom}{\operatorname{Hom}}
\newcommand{\Ext}{\operatorname{Ext}}

\newcommand{\NN}{\mathbb{N}}
\newcommand{\ZZ}{\mathbb{Z}}
\newcommand{\QQ}{\mathbb{Q}}
\newcommand{\RR}{\mathbb{R}}
\newcommand{\CC}{\mathbb{C}}
\newcommand{\HH}{\mathbb{H}}

\newcommand{\pder}[2]{\frac{\partial #1}{\partial #2}}

\newcommand{\vfeld}[1]{\frac{\partial}{\partial #1}}
\newcommand{\vfelds}[1]{\partial / \partial #1}

\newcommand{\menge}[2]{\bigl\{ \thinspace #1 \thinspace\thinspace \big\vert%
\thinspace\thinspace #2 \thinspace \bigr\}}
\newcommand{\Menge}[2]{\Bigl\{ \thinspace #1 \thinspace\thinspace \Big\vert%
\thinspace\thinspace #2 \thinspace \Bigr\}}

\DeclareMathOperator{\im}{im}

\DeclareMathOperator{\ad}{ad}

\DeclareMathOperator{\Res}{Res}
\DeclareMathOperator{\Spec}{Spec}

\DeclareMathOperator{\id}{id}

\renewcommand{\Im}{\operatorname{Im}}
\renewcommand{\Re}{\operatorname{Re}}

\DeclareMathOperator{\rat}{rat}

\DeclareMathOperator{\Sym}{Sym}
\DeclareMathOperator{\Gr}{Gr}
\DeclareMathOperator{\DR}{DR}

\newcommand{\define}[1]{\emph{#1}}

\newcommand{\lie}[2]{\lbrack #1, #2 \rbrack}

\newcommand{\bigglie}[2]{\biggl[ #1, #2 \biggr]}
\newcommand{\biglie}[2]{\bigl[ #1, #2 \bigr]}

\newcommand{\shf}[1]{\mathscr{#1}}
\newcommand{\OX}{\shf{O}_X}
\newcommand{\OmX}[1]{\Omega_X^{#1}}

\newcommand{\restr}[1]{\big\vert_{#1}}

\newcommand{\argbl}{-}

\def\overbar#1#2#3{{%
	\setbox0=\hbox{\displaystyle{#1}}%
	\dimen0=\wd0
	\advance\dimen0 by -#2 
	\vbox {\nointerlineskip \moveright #3 \vbox{\hrule height 0.3pt width \dimen0}%
		\nointerlineskip \vskip 1.5pt \box0}%
}}

\newcommand{\dst}{\Delta^{\ast}}
\newcommand{\dstn}[1]{(\dst)^{#1}}

\newcommand{\into}{\hookrightarrow}
\newcommand{\onto}{\twoheadrightarrow}

\usepackage{ifthen}
\DeclareMathOperator{\SL}{SL}
\newcommand{\sltwo}{\mathfrak{sl}_2(\CC)}

\newcommand{\Dcheck}{\check{D}}
\newcommand{\Npl}{N^{+}}
\newcommand{\glie}{\mathfrak{g}}

\newcommand{\qlie}{\mathfrak{q}}

\newcommand{\Wn}[1]{W^{(#1)}}

\newcommand{\Cn}[1]{C^{(#1)}}

\newcommand{\HR}{H_{\RR}}
\newcommand{\HZ}{H_{\ZZ}}
\newcommand{\HQ}{H_{\QQ}}
\newcommand{\HC}{H_{\CC}}
\newcommand{\Phit}{\tilde{\Phi}}

\newcommand{\shHO}{\shf{H}_{\shf{O}}}
\newcommand{\shHOe}{\shHO^{\mathrm{e}}}
\newcommand{\shHZ}{\shf{H}_{\ZZ}}
\newcommand{\shHQ}{\shf{H}_{\QQ}}
\newcommand{\shHC}{\shf{H}_{\CC}}

\newcommand{\Hd}{H^{\vee}}
\newcommand{\Vd}{V^{\vee}}
\newcommand{\VdC}{V_{\CC}^{\vee}}
\newcommand{\HdC}{H_{\CC}^{\vee}}

\newcommand{\shV}{\shf{V}}
\newcommand{\shVd}{\shf{V}^{\vee}}
\newcommand{\shVZ}{\shV_{\ZZ}}
\newcommand{\shVC}{\shV_{\CC}}
\newcommand{\shVQ}{\shV_{\QQ}}
\newcommand{\shVO}{\shV_{\shf{O}}}
\newcommand{\shVOd}{\shV_{\shf{O}}^{\vee}}
\newcommand{\class}[1]{\lbrack #1 \rbrack}
\newcommand{\VC}{V_{\CC}}
\newcommand{\VR}{V_{\RR}}
\newcommand{\VZ}{V_{\ZZ}}
\newcommand{\VQ}{V_{\QQ}}
\newcommand{\VRone}{V_{\RR,1}}
\newcommand{\VCone}{V_{\CC,1}}
\newcommand{\VZone}{V_{\ZZ,1}}

\newcommand{\vZ}{v_{\ZZ}}
\newcommand{\Tnu}{T_{\nu}}
\newcommand{\Gammanu}{\Gamma_{\nu}}
\newcommand{\Nnu}{N_{\nu}}
\newcommand{\Nmodnu}{\Nmod_{\nu}}

\newcommand{\Jb}{\bar{J}}
\newcommand{\Jbtor}{\Jb_{\mathrm{tor}}}
\newcommand{\Xb}{\bar{X}}
\newcommand{\Yb}{\bar{Y}}
\newcommand{\OXb}{\shf{O}_{\Xb}}
\newcommand{\TZ}{T_{\ZZ}}
\newcommand{\TZx}{T_{\ZZ,x}}

\newcommand{\ZZX}{\ZZ_X}
\newcommand{\ZZXb}{\ZZ_{\Xb}}
\newcommand{\nub}{\bar{\nu}}
\newcommand{\JBPS}{\bar{J}^{\mathit{BPS}}}
\newcommand{\PC}{P_{\CC}}

\newcommand{\OW}{\shf{O}_W}
\newcommand{\Pb}{\bar{P}}
\newcommand{\OPb}{\shf{O}_{\Pb}}
\newcommand{\OmW}[1]{\Omega_W^{#1}}
\newcommand{\famX}{\mathscr{X}}
\newcommand{\pil}{\pi_{\ast}}
\newcommand{\Hvan}[1]{H_{\mathrm{van}}^{#1}}
\newcommand{\Jvan}{J_{\mathrm{van}}}
\newcommand{\Kvan}{K_{\mathrm{van}}}

\newcommand{\jl}{j_{\ast}}
\newcommand{\jlreg}{j_{\ast}^{\mathrm{reg}}}

\newcommand{\jlsl}{j_{!\ast}}
\newcommand{\ju}{j^{\ast}}
\newcommand{\fu}{f^{\ast}}
\newcommand{\fl}{f_{\ast}}
\newcommand{\fus}{f^{!}}
\newcommand{\iu}{i^{\ast}}
\newcommand{\ius}{i^{!}}

\newcommand{\pZ}{p_{\ZZ}}
\newcommand{\pZu}{p_{\ZZ}^{\ast}}
\newcommand{\qu}{q^{\ast}}
\newcommand{\qus}{q^{!}}

\newcommand{\iku}[1]{i_{#1}^{\ast}}
\newcommand{\ikus}[1]{i_{#1}^{!}}

\newcommand{\DmodXb}{\Dmod_{\Xb}}

\newcommand{\ODelta}{\shO_{\Delta}}

\newcommand{\QHXb}{\QQ_{\Xb}^H}
\newcommand{\QHDelta}{\QQ_{\Delta}^H}
\newcommand{\DD}{\mathbb{D}}
\DeclareMathOperator{\can}{can}
\DeclareMathOperator{\Var}{Var}
\newcommand{\psione}[1]{\psi_{#1,1}}
\newcommand{\phione}[1]{\phi_{#1,1}}

\newcommand{\shF}{\shf{F}}
\newcommand{\shG}{\shf{G}}
\newcommand{\shE}{\shf{E}}
\newcommand{\shFd}{\shF^{\vee}}
\newcommand{\shEd}{\shE^{\vee}}
\newcommand{\Ed}{E^{\ast}}
\DeclareMathOperator{\Map}{Map}
\newcommand{\OY}{\shf{O}_Y}

\newcommand{\alphast}{\alpha^{\ast}}
\newcommand{\betast}{\beta^{\ast}}
\newcommand{\shO}{\shf{O}}

\newcommand{\Blog}{\mathcal{B}(\log)}
\newcommand{\DSigma}{D_{\Sigma}}
\newcommand{\DSigmap}{D_{\Sigma'}'}
\newcommand{\Phib}{\bar{\Phi}}

\newcounter{thmA}
\newtheorem{theoremA}[thmA]{Theorem}

\begin{document}

%========================================================
\title[Complex analytic N{\'e}ron models]{Complex analytic N{\'e}ron models for
arbitrary families of intermediate Jacobians}
\author[C.~Schnell]{Christian Schnell}
\address{Department of Mathematics, Statistics \& Computer Science \\
University of Illinois at Chicago \\
851 South Morgan Street \\
Chicago, IL 60607}
\email{cschnell@math.uic.edu}
\begin{abstract}
Given a family of intermediate Jacobians (for a polarized variation of Hodge
structure of weight $-1$) on a Zariski-open subset of a complex manifold, we
construct an analytic space that naturally extends the family. 
Its two main properties are: (a) the horizontal and holomorphic sections are
precisely the admissible normal functions without singularities; (b) the graph of any
admissible normal function has an analytic closure inside our space. As a
consequence, we obtain a new proof for the zero locus conjecture of M.~Green and
P.~Griffiths.  The construction uses filtered $\Dmod$-modules and M.~Saito's theory
of mixed Hodge modules; it is functorial, and does not require normal crossing or
unipotent monodromy assumptions. 
%We also show that our space maps continuously to the identity component of the
%N\'eron model defined by P.~Brosnan, G.~Pearlstein, and M.~Saito.
\end{abstract}
\dedicatory{Dedicated to Herb Clemens on the occasion of his 70th birthday}
\maketitle
%========================================================

\section{Overview}

\subsection{Introduction}

Not too long ago, during a lecture at the Institute for Pure and Applied Mathematics,
P.~Griffiths mentioned the problem of constructing N\'eron models for arbitrary
families of intermediate Jacobians. In other words, given a family of intermediate
Jacobians over a Zariski-open subset $X$ of a complex manifold $\Xb$, one should
construct a space that extends the family to all of $\Xb$. This has to be done in
such a way that normal functions extend to sections of the N\'eron model.

It is known that two additional conditions need to be imposed to make this into a
reasonable question. Firstly, the family of intermediate Jacobians should come from a
\emph{polarizable} variation of Hodge structure, which we may normalize to be of weight
$-1$. Secondly, one should only consider \emph{admissible} normal functions. The
N\'eron model is then expected to have the following structure: (1) Over each point of $\Xb$,
its fiber should be a countable union of complex Lie groups. (2) The components over a
point $x \in \Xb - X$ where the variation degenerates should be indexed by a
countable group, whose elements are the possible values for the \define{singularity}
at $x$ of admissible normal functions---an invariant introduced by M.~Green and
P.~Griffiths \cite{GG1} that measures whether the cohomology class of a normal
function is trivial in a neighborhood of $x$. (3) The horizontal sections of the identity
component of the N\'eron model should be the admissible normal functions
without singularities.

The existence of N\'eron models with good properties has useful consequences,
for instance, a proof of the following conjecture by M.~Green and P.~Griffiths:
\begin{conjecture} \label{conj:GG}
Let $\nu$ be an admissible normal function on an algebraic variety $X$. Then the
zero locus $Z(\nu)$ is an algebraic subvariety of $X$.
\end{conjecture}
By Chow's Theorem, it suffices to show that the closure of $Z(\nu)$ inside a
projective compactification $\Xb$ remains analytic; this is almost automatic if we
assume that $\nu$ can be extended to a section of a N\'eron model over $\Xb$ with
good properties. M.~Saito has established Conjecture~\ref{conj:GG} for
$\dim X = 1$ by this method \cite{Saito-GGK}; an entirely different approach has been
pursued by P.~Brosnan and G.~Pearlstein \cites{BP1,BP2}, who have announced a full
proof this summer \cite{BP3}.

In this paper, we largely solve P.~Griffiths' problem, by constructing an analytic
space that has all the properties expected for the identity component of the N\'eron
model---in particular, its horizontal and holomorphic sections are precisely the
admissible normal functions without singularities.  We also show that the graph of
any admissible normal function has an analytic closure inside our space; one
consequence is a new proof for Conjecture~\ref{conj:GG}. Lastly, we describe the
construction of an analytic N\'eron model for admissible normal functions with
torsion singularities. Based on some examples, we argue that this is the most general
setting in which the N\'eron model exists as an analytic space or even as a Hausdorff
space.

The construction that is proposed here is very natural and suitably functorial; it is
motivated by unpublished work of H.~Clemens on the family of hypersurface sections of
a smooth projective variety (briefly reviewed in \subsecref{subsec:background}
below). An important point is that no assumptions on the singularities of $D = \Xb -
X$, or on the local monodromy of the variation of Hodge structure are needed.
(Whereas the traditional approach would be to make $D$ into a divisor with normal
crossings by using resolution of singularities, and then to pass to a finite cover to
get unipotent monodromy.) We accomplish this generality with the help of M.~Saito's
theory of mixed Hodge modules \cite{Saito-MHM}.

Two other solutions to the problem have been given recently. One is by P.~Brosnan,
G.~Pearlstein, and M.~Saito \cite{BPS}, whose N\'eron model is a topological space to
which admissible normal functions extend as continuous sections. They also show that
the base manifold $\Xb$ can be stratified in such a way that, over each stratum,
their space is a family of complex Lie groups, and the extended normal function a
holomorphic section. Unfortunately, it is not clear from the construction whether the
resulting space is Hausdorff; and when the local monodromy of $\shH$ is not
unipotent, the fibers of their N\'eron model can be too small, even in one-parameter
degenerations of abelian varieties. We address both issues in \subsecref{subsec:BPS}
below, by showing that there is always a continuous and surjective map from
$\Jb(\shH)$ to the identity component of their N\'eron model.

A second solution is contained in a preprint by K.~Kato, C.~Nakayama, and S.~Usui
\cite{KNU}, who use classifying spaces of pure and mixed nilpotent orbits to define a
N\'eron model in the category $\Blog$. At present, their construction is only
available for $\dim X = 1$; but it is expected to work in general, at least when $D$
is a normal crossing divisor and $\shH$ has unipotent monodromy. It seems likely that
there will be a connection between the identity component of their N\'eron model and
the subset of $\Jb(\shH)$ defined by the horizontality condition in
\subsecref{subsec:horizontality}. This question is currently under investigation by
T.~Hayama.

\subsection{Background}
\label{subsec:background}

The idea for constructing the analytic space $\Jb(\shH)$ goes back to unpublished
work of H.~Clemens, for the case of hypersurface sections of an even-dimensional
variety. To motivate what follows, we briefly describe this development.

Let $W$ be a smooth projective variety of dimension $2m$, and consider the family
of its hypersurface sections of large degree, parametrized by the projective space $\Pb = \lvert
\OW(d) \rvert$ with $d \gg 0$.  Denote by $D \subseteq \Pb$ the dual variety; then $P
= \Pb - D$ parametrizes smooth hypersurfaces. Let $\pi \colon \famX \to P$ be the
universal family, and let $R^{2m-1} \pil \ZZ_{\famX}(m)$ be the variation of Hodge
structure on the cohomology of the fibers, normalized to be of weight $-1$. 

Consider one of the smooth hypersurface sections $X$. A basic fact, due to P.~Griffiths
in the case of projective space, and to M.~Green \cite{Green} in general, is that the
vanishing cohomology of $X$ is generated by residues of meromorphic forms; moreover,
the Hodge filtration is essentially the filtration by pole order. In particular, for
$d \gg 0$, the residue map $\Res \colon H^0 \bigl( W, \OmW{2m}(mX) \bigr) \to F^m
\Hvan{2m-1}(X, \CC)$ is surjective. H.~Clemens observed that, consequently, the
intermediate Jacobian
\[
	\Jvan(X) = \frac{\bigl( F^m \Hvan{2m-1}(X, \CC) \bigr)^{\vee}}%
		{\Hvan{2m-1}(X, \ZZ)}
\]
is a subspace of the bigger object
\[
	\Kvan(X) = \frac{\bigl( H^0 \bigl( W, \OmW{2m}(mX) \bigr)^{\vee}}%
		{\Hvan{2m-1}(X, \ZZ)}.
\]
The original motivation for introducing $\Kvan(X)$ was to extend the Abel-Jacobi map
to ``topological cycles,'' and to obtain a form of Jacobi inversion for such cycles.
But it is also clear that $H^0 \bigl( W, \OmW{2m}(mX) \bigr)$ is isomorphic to the
space of sections of the line bundle $\OmW{2m} \tensor \OW(m)$; therefore the
numerator in the definition of $\Kvan(X)$ is essentially independent of $X$, and makes sense
even when $X$ becomes singular. This suggests that residues might be useful in
extending the family of intermediate Jacobians from $P$ to $\Pb$.

Let $\shH \subseteq R^{2m-1} \pil \ZZ_{\famX}(m)$ be the variation of Hodge structure
on the vanishing cohomology, and $\bigl( \shHO, \nabla \bigr)$ the corresponding flat
vector bundle. It can be shown that the residue calculus extends to
the family of all hypersurface sections, including the singular ones, in the
following way: Let $j \colon P \into \Pb$ be the inclusion, and define subsheaves
$F_p \Mmod$ of $\jl \shHO$ by the condition that a section in $H^0 \bigl( U \cap P,
\shHO \bigr)$ belongs to $H^0(U, F_p \Mmod)$ iff it is the residue of a meromorphic
$2m$-form on $U \times W$ with a pole of order at most $m+p$ along the incidence
variety. Let $\Mmod$ be the union of the $F_p \Mmod$; then $\Mmod$ is a holonomic
$\Dmod$-module on $\Pb$, extending the flat vector bundle, and $F_{\bullet} \Mmod$ is
a good filtration. It was proved in \cite{Schnell} that $(\Mmod, F)$ underlies a
polarized Hodge module on $\Pb$, namely the intermediate extension $M = \jlsl \shH
\decal{\dim P}$ of the variation of Hodge structure (at least when $d \gg 0$).  This
is how filtered $\Dmod$-modules and M.~Saito's theory introduce themselves into the problem.

The nice geometry of the family of hypersurfaces, especially the fact that $\Pb$ is a
projective space, is the primary motivation for trying to construct the N\'eron model
without resolving singularities and without passing to a finite cover.  Moreover,
each sheaf $F_p \Mmod$ is the quotient of $H^0 \bigl( W, \OmW{2m}(m+p) \bigr) \tensor
\OPb(m+p)$; therefore $T(F_0 \Mmod)$ is a submanifold of the anti-ample vector bundle
with sheaf of sections $H^0 \bigl( W, \OmW{2m}(m) \bigr)^{\vee} \tensor \OPb(-m)$.
This important fact gives $\Jb(\shH) \to \Pb$ many good properties that will be
explained in a separate article (currently in preparation); it may also place
restrictions on sections of $\Jb(\shH)$, i.e., on normal functions without
singularities.  This is of interest because M.~Green and P.~Griffiths
have related the existence of singularities of normal functions to the
Hodge conjecture \cites{GG1,GG2}.

\subsection{Summary of the main results}
\label{subsec:summary}

We now describe the construction of the analytic space $\Jb(\shH)$, and summarize the
main results of the paper. Throughout, we let $\shH$ be a polarized variation of
Hodge structure of weight $-1$, defined on a Zariski-open subset $X$ of a complex
manifold $\Xb$. We denote the corresponding family of intermediate Jacobians by
$J(\shH) \to X$.

To begin with, consider a single polarized Hodge structure $H$ of weight $-1$. To
emphasize the analogy with what comes later, we view the Hodge filtration as an
increasing filtration $F_{\bullet} \HC$ by setting $F_p \HC = F^{-p} \HC$. We also
let $\HZ$ be the integral lattice, and $Q \colon \HZ \tensor \HZ \to \ZZ$ the
polarization. Since $H$ has weight $-1$, the polarization induces an isomorphism $\HC
/ F_0 \HC \simeq (F_0 \HC)^{\vee}$; this justifies defining the \define{intermediate
Jacobian} as
\[
	J(H) = (F_0 \HC)^{\vee} / \HZ,
\]
where the map $\HZ \to (F_0 \HC)^{\vee}$ is given by $h \mapsto Q(h, \argbl)$. The
advantage of this point of view is that an extension of mixed Hodge structures
\begin{diagram}
0 &\rTo& H &\rTo& V &\rTo& \ZZ(0) &\rTo& 0
\end{diagram}
of ``normal function type'' determines a point in $J(H)$ with only one choice: after
dualizing the extension, one has $F_{-1} \VdC \simeq F_0 \HC$ since $H$ is polarized;
now any element $\vZ \in \VZ$ lifting $1 \in \ZZ$ defines a linear functional on
$F_0 \HC$, and hence a point in $J(H)$.

Similarly, the sheaf of sections of the family $J(\shH) \to X$ is given by $(F_0
\shHO)^{\vee} / \shHZ$, where $\shHZ$ is the local system underlying the variation,
and $F_{\bullet} \shHO$ the Hodge filtration on the associated vector bundle.  To
extend this in a natural way to $\Xb$, we view $\shH$ as a polarized Hodge module on
$X$; according to M.~Saito's theory, it can be extended in a canonical manner to a
polarized Hodge module $M = \jlsl \shH \decal{\dim X}$ on $\Xb$. The holonomic $\Dmod$-module
$\Mmod$ underlying $M$ is known to be the \define{minimal extension} of the flat
vector bundle $(\shHO, \nabla)$; in particular, its de Rham complex $\DR(\Mmod)$ is
isomorphic to the intersection complex $\operatorname{\mathit{IC}}_{\!\Xb}(\shHC)$. The
$\Dmod$-module comes with a good filtration $F = F_{\bullet} \Mmod$ by
$\OXb$-coherent subsheaves; $F_p \Mmod$ is difficult to describe in general,
but may be viewed as a natural extension of the Hodge bundle $F_p \shHO$.

Guided by the above, we define the space $\Jb(\shH)$ in such a way that its sheaf of
holomorphic sections is $(F_0 \Mmod)^{\vee} / \jl \shHZ$. Namely, we let $T(F_0
\Mmod)$ be the analytic spectrum of the symmetric algebra of $F_0 \Mmod$ (see
\subsecref{subsec:spaces}), and $\TZ$ the \'etal\'e space of the sheaf $\jl \shHZ$.
Using the polarization, we show that there is a natural holomorphic map $\eps \colon
\TZ \to T(F_0 \Mmod)$ (see \subsecref{subsec:construction}). The main technical
result of the paper is that the image of $\eps$ is a closed analytic subset of
$T(F_0 \Mmod)$. 

\begin{theoremA} \label{thm:A}
The map $\eps \colon \TZ \to T(F_0 \Mmod)$ is a proper holomorphic embedding.
Consequently, the fiberwise quotient space $T(F_0 \Mmod) / \TZ$ is an
analytic space, and in particular Hausdorff.
\end{theoremA}

The second statement follows from the first by simple topological arguments (see
\subsecref{subsec:quotients}). We now define $\Jb(\shH) = T(F_0 \Mmod) / \TZ$; this
is an analytic space over $\Xb$ that naturally extends the family of intermediate
Jacobians. Further evidence that it is a good candidate for the identity component of
the N\'eron model is given by the following list of properties:
\begin{enumerate}
\item Every normal function on $X$ that is admissible (relative to $\Xb$) and without
singularities extends to a holomorphic section of $\Jb(\shH)$
(Proposition~\ref{prop:nub}). In fact, the process that gives the extension is
analogous to the one for a single Hodge structure, explained above.
\item There is a notion of horizontality for sections of $\Jb(\shH)$, and
the holomorphic and horizontal sections are precisely the admissible normal functions
without singularities (Proposition~\ref{prop:sections}). 
\item The construction is functorial, in the following sense: Given a holomorphic map
$f \colon \Yb \to \Xb$ such that $Y = f^{-1}(X)$ is dense in $\Yb$, let $\fu \shH$
denote the pullback of the variation to $Y$. Then there is a canonical holomorphic map
\[
	\Yb \times_{\Xb} \Jb(\shH) \to \Jb(\fu \shH),
\]
compatible with normal functions (Proposition~\ref{prop:Jb-functoriality}).
\item There is a continuous and surjective map from $\Jb(\shH)$ to the identity component
of the N\'eron model defined in \cite{BPS}, compatible with normal functions, whose
effect is to contract parts of certain fibers (Lemma~\ref{lem:blowdown}).
\end{enumerate}

A few words about the proof of Theorem~\ref{thm:A}. We use results from M.~Saito's
theory, in particular nearby and vanishing cycle functors and their description in
terms of the $V$-filtration of M.~Kashiwara and B.~Malgrange, to reduce the general
problem to the case where $\Xb = \Delta^n$, $X = \dstn{n}$, and the local system
$\shHZ$ has unipotent monodromy (see \subsecref{subsec:reduction}). In that case,
there is an explicit description of the sheaf $F_0 \Mmod$ in terms of P.~Deligne's
canonical extension of $(\shHO, \nabla)$: for every $k \geq 0$, $F_0 \Mmod$ contains
all $k$-th derivatives of sections in $F_{-k} \shHOe = F^k \shHOe$. In particular, we
have a holomorphic map $T(F_0 \Mmod) \to T(F_0 \shHOe)$. In general, the image of $\TZ$
in $T(F_0 \shHOe)$ is badly behaved, which comes from the fact that
sections of $F_0 \shHOe$ are not sufficient to separate sections of $\shHZ$ ``in the
limit.'' The following result shows that $F_0 \Mmod$ has enough additional sections to
overcome this problem.

\begin{theoremA} \label{thm:B}
Let $\HC$ denote the space of sections of $\shHC$ on the universal covering space
$\HH^n$, let $Q$ be the polarization, and let $N_1, \dotsc, N_n$ be the logarithms of
the monodromy operators. Also let $\sigma_1, \dotsc, \sigma_r$ be a collection
of holomorphic sections that generate $F_0 \Mmod$ on $\Delta^n$. Then there are constants $C > 0$
and $\alpha > 0$, such that for every $z \in \HH^n$ and every real vector $h \in \HR$, 
\begin{equation} \label{eq:B}
	\max_{k \geq 0} \norm{(y_1 N_1 + \dotsb + y_n N_n)^k h} \leq 
		C \cdot \max_{1 \leq j \leq r} \bigabs{Q \bigl( h, \sigma_j(z) \bigr)},
\end{equation}
provided that $y_j = \Im z_j \geq \alpha$ and $0 \leq \Re z_j \leq 1$ for all $j=1,
\dotsc, n$.
\end{theoremA}

The estimate \eqref{eq:B}, which is proved in \subsecref{subsec:proof} below, quickly
leads to the proof of Theorem~\ref{thm:A} in the normal crossing case. We obtain it
essentially by linear algebra methods, using only familiar consequences of the
$\SL_2$-Orbit Theorem \cite{CKS}.

Perhaps surprisingly, the space $\Jb(\shH)$ is also useful for the study of normal
functions with nontrivial singularities. Of course, such normal functions cannot
be extended to holomorphic sections; nevertheless, the following is true.

\begin{theoremA} \label{thm:C}
Let $\nu \colon X \to J(\shH)$ be a normal function, admissible relative to
$\Xb$. Then the topological closure of its graph inside $\Jb(\shH)$ is a closed
analytic subset.
\end{theoremA} 

This result clearly implies that the closure of the zero locus of $\nu$ is an
analytic subset of $\Xb$, and thus gives a different proof for
Conjecture~\ref{conj:GG}. To prove Theorem~\ref{thm:C} in the normal crossing case
(see \subsecref{subsec:graphs-NC}), we require one consequence of the $\SL_2$-Orbit
Theorem of \cite{KNU-SL}, namely the boundedness of the canonical splitting in mixed
nilpotent orbits. Beyond that, only simple arguments from linear algebra are needed.
The reader who is mainly interested in the proof of Conjecture~\ref{conj:GG} can
focus on Part~\ref{sec:normal-crossing} of the paper, where mixed Hodge modules play
no role.

For admissible normal functions with \emph{torsion} singularities, it turns out
(Proposition~\ref{prop:torsion}) that there is always a maximal extension whose graph
is closed inside of $\Jb(\shH)$. As a consequence, it is possible to construct a
N\'eron model for this class of normal functions by a gluing construction (see
\subsecref{subsec:torsion}).

\begin{theoremA} \label{thm:D}
There is an analytic space $\Jbtor(\shH) \to \Xb$ whose holomorphic and horizontal
sections are precisely the admissible normal functions with torsion singularities. It
contains $\Jb(\shH)$ as the identity component, and has similar functoriality
properties.
\end{theoremA}

Unfortunately, it appears that admissible normal functions with torsion singularities
are the biggest class for which there exists a N\'eron model that is an analytic
space (or a Hausdorff space). The reason is the following: Over a point in $\Xb$
where an admissible normal function has a non-torsion singularity, the closure of its
graph can have a fiber of positive dimension. This happens even in very simple
examples, such as two-parameter families of elliptic curves.  As we argue in
\subsecref{subsec:impossible} below, it is therefore unlikely that there can be a
N\'eron model that (a) graphs all admissible normal functions, (b) has a reasonable
identity component, and (c) is Hausdorff as a topological space. Nevertheless, the
result of Theorem~\ref{thm:C} in itself is probably sufficient to study singularities
of normal functions in the way proposed in \cite{GG2}, without having such a more general
N\'eron model.

Three examples have been included in Part~\ref{sec:examples}, to illustrate
different aspects of the construction. On the other hand, given the length of the
paper, we have not included any background on mixed Hodge modules,
degenerations of variations of Hodge structure, or admissible normal functions. Here,
the reader should consult the following sources: (1) for mixed Hodge modules, the survey
paper \cite{Saito-survey}; (2) for degenerations of variations of Hodge
structure, the paper \cite{CKS}; (3) for a discussion of admissibility, the papers
\cite{Kashiwara-study} and \cite{Saito-ANF}.

\subsection{History of N\'eron models}

N\'eron models originated in a construction of A.~N\'eron for abelian varieties
\cite{Neron}. Let $A$ be an abelian variety, defined over the field of functions $K$
of a Dedekind domain $D$. Then the \define{N\'eron model} for $A$ is a smooth and
commutative group scheme $\mathcal{A}$ over $R$, such that $\mathcal{A}(S) =
\mathcal{A}(S \times_R K)$ for any smooth morphism $S \to R$; more details can be
found in the book \cite{BLR}. The definition means that $\mathcal{A}$ is the natural
extension of $A$ from the open subset $\Spec K$ to all of $\Spec R$.

In the complex-analytic setting, a family of abelian varieties is a special case of a
polarized variation of Hodge structure. After P.~Griffiths popularized the use of normal
functions in Hodge theory, N\'eron models for one-parameter degenerations of more
general variations were constructed by several people. S.~Zucker \cite{Zucker}
introduced a generalized intermediate Jacobian for hypersurfaces with one ordinary
double point, and used it to define the identity component of a N\'eron model in
Lefschetz pencils. H.~Clemens \cite{Clemens} extended this to the construction of a
N\'eron model for one-parameter degenerations with certain restrictions on the local
monodromy. In his paper on admissible normal functions, M.~Saito \cite{Saito-ANF}
generalized both constructions to arbitrary one-parameter degenerations, and also
constructed a compactification of the ``Zucker extension'' (which, however, is
usually not Hausdorff). 

The recent interest in N\'eron models stems from the work by M.~Green, P.~Griffiths, and
M.~Kerr \cite{GGK}, who observed that a subspace of the Zucker extension is sufficient
to graph admissible normal functions without singularities. Briefly summarized, their
construction works as follows: Let $\Xb$ be a smooth curve, and $\shH$ a polarized
variation of Hodge structure with unipotent monodromy, defined on a Zariski-open
subset $X$. At each of the points $x \in \Xb - X$, a choice of local coordinate
determines an asymptotic mixed Hodge structure; the monodromy-invariant part $H =
\ker(T - \id)$ is independent of that choice. The identity component of the N\'eron model
in \cite{GGK} has the generalized intermediate Jacobian $J(H) = \HC / (F^0 \HC +
\HZ)$ as its fiber over $x$. M.~Green, P.~Griffiths, and M.~Kerr also defined the full
N\'eron model that graphs arbitrary admissible normal functions, and computed its
finite group of components at each point of $\Xb - X$. Their construction produces
a so-called ``slit'' analytic space; M.~Saito \cite{Saito-GGK} has shown that the
resulting topological space is Hausdorff.

As mentioned above, a construction of a N\'eron model for $\Xb$ of arbitrary
dimension has been proposed by P.~Brosnan, G.~Pearlstein, and M.~Saito \cite{BPS}.
They observe that, at each point $x \in \Xb$, the stalk $H_x$ of the sheaf $R^1 \jl
\shHZ$ carries a mixed Hodge structure of weight $-1$, and therefore defines a
generalized intermediate Jacobian $J(H_x) = \Ext_{\MHS}^1 \bigl( \ZZ(0), H_x \bigr)$.
The identity component of their N\'eron model is the disjoint union of the complex
Lie groups $J(H_x)$, topologized in a rather tricky way by reduction to the
normal-crossing case. The full N\'eron model is then obtained by gluing.
As pointed out in \cite{BPS}, the construction does not seem to work
very well in the case of non-unipotent local monodromy.

The most recent work, also alluded to above, is by K.~Kato, C.~Nakayama, and S.~Usui
\cite{KNU}. The variation of Hodge structure $\shH$ determines a period map $\Phi \colon X \to \Gamma
\backslash D$, and according to the general theory in \cite{KU-book}, it can be 
extended to a map $\Phib \colon \Xb \to \Gamma \backslash \DSigma$, where $\DSigma$
is a space of nilpotent orbits. For $\dim X = 1$, they show that there is a good
choice of compatible fan $\Sigma'$, such that an admissible normal function defines a
map from $\Xb$ into a space $\DSigmap$ of nilpotent orbits of normal function type.
The N\'eron model can then be constructed as the fiber product
\begin{diagram}[width=3em,midshaft]
	\Jb_{\Sigma'}(\shH) &\rTo& \Gamma' \backslash \DSigmap \\
	\dTo && \dTo \\
	\Xb &\rTo^{\Phib}& \Gamma \backslash \DSigma
\end{diagram}
in the category $\Blog$. It is hoped that a similar construction will work as long
as $\Xb - X$ is a normal crossing divisor and $\shHZ$ has unipotent local monodromy.
The authors point out that the construction is not entirely canonical, since
it depends on the choice of fan $\Sigma'$.

For families of complex abelian varieties, there is a complete construction of a
N\'eron model by A.~Young \cite{Young}, at least when $\Xb - X$ is a divisor with
normal crossings and the local monodromy of $\shHZ$ is unipotent. His construction
uses toric geometry; the identity component of his model agrees with an
older construction by Y.~Namikawa \cite{Namikawa} for degenerations of abelian
varieties (and, therefore, with the model that is proposed in this paper), and is in
particular a complex manifold. When all components are considered together, the space
is however not Hausdorff.

\subsection{Conventions}
\label{subsec:conventions}

In dealing with filtrations, we index increasing filtrations (such as weight
filtrations, or Hodge filtrations on $\Dmod$-modules) by lower indices, and decreasing
filtrations (such as Hodge filtrations on vector spaces, or $V$-filtrations on left
$\Dmod$-modules) by upper indices. We may pass from one to the other by the
convention that $F^{\bullet} = F_{-\bullet}$. To be consistent, shifts in the
filtration thus have different effects in the two cases:
\[
	F \decal{1}^{\bullet} = F^{\bullet+1}, \quad \text{while} \quad
	F \decal{1}_{\bullet} = F_{\bullet-1}.
\]
This convention agrees with the notation used in M.~Saito's papers. 

When dealing with mixed Hodge modules and mixed Hodge structures (or variations of
mixed Hodge structure) at the same time, we usually consider the Hodge filtrations on
the latter as increasing filtrations.

In this paper, we work with \emph{left} $\Dmod$-modules, and $\Dmod$-module always
means left $\Dmod$-module (in contrast to \cite{Saito-MHM}, where right
$\Dmod$-modules are used).

When $M$ is a mixed Hodge module, the effect of a Tate twist $M(k)$ on the underlying
filtered $\Dmod$-module $(\Mmod, F)$ is as follows: 
\[
	(\Mmod, F)(k) = \bigl( \Mmod, F \decal{k} \bigr) = 
		\bigl( \Mmod, F_{\bullet-k} \bigr).
\]

For a regular holonomic $\Dmod$-module $\Mmod$ that is defined on the complement of a
divisor $D \subseteq X$, we let $\jlreg \Mmod$ be the direct image in the category of
regular holonomic $\Dmod$-modules; its sections have poles of finite order along $D$.

The dual of a complex vector space $V$ will be denoted by $\Vd = \Hom_{\CC}(V, \CC)$.
Similar notation is used for mixed Hodge structures and for coherent sheaves.

The individual sections of the paper are numbered consecutively, and are referred to
with a paragraph symbol (such as \subsecref{subsec:conventions}).

\subsection{Acknowledgments}

The paper owes a lot to Greg Pearlstein, and to the techniques that he has developed
for dealing with the zero locus problem. In particular, I learned about the normal
form for period maps and the $\SL_2$-Orbit Theorem of \cite{KNU} from him; and an
earlier proof of the main theorem in the normal crossing case used several other
ideas from his work. I am very grateful for his help.

I thank Morihiko Saito for many valuable comments about the construction, mostly
in the case of one variable, and for help with questions about mixed Hodge modules.
He also corrected my earlier (and mistaken) opinions about the comparison with the
N\'eron model of \cite{BPS}.

Patrick Brosnan suggested Proposition~\ref{prop:torsion} and asked about the
construction of a N\'eron model for normal functions with torsion singularities. I
thank him for these useful contributions to the paper.

Most of all, I thank my former thesis adviser, Herb Clemens, to whom
the basic idea behind the construction is due. I believe that, when applied to the
family of hypersurface sections of a smooth projective variety, the N\'eron model
proposed here is the space that he envisioned with ``$J \subseteq K$''.

%%%%%%%%%%%%%%%%%%%%%%%%%%%
\section{The construction of the analytic space}

\subsection{Intermediate Jacobians}
\label{subsec:Jacobians}

Let $H = \bigl( \HC, F_{\bullet} \HC, \HZ, Q \bigr)$ be a polarized Hodge structure
of weight $-1$ (we set $F_p \HC = F^{-p} \HC$ for consistency with later sections). 
There are two ways of defining the associated intermediate Jacobian; the
usual definition,
\[
	J_1(H) = \frac{\HC}{F_0 \HC + \HZ},
\]
does not use the fact that $H$ is polarized. But since $Q(F_0 \HC, F_0 \HC) = 0$, the
polarization induces an isomorphism $\HC / F_0 \HC \simeq (F_0 \HC)^{\vee}$, and
therefore 
\[
	J_2(H) = \frac{(F_0 \HC)^{\vee}}{\HZ} \simeq J_1(H),
\]
where the map $\HZ \into (F_0 \HC)^{\vee}$ is given by $h \mapsto Q(h, \argbl)$.
One theme of this paper is that the \emph{second} definition is the correct
one. With this in mind, we briefly review the correspondence between extensions of
mixed Hodge structure of the form
\begin{equation} \label{eq:ext-MHS}
\begin{diagram}
	0 &\rTo& H &\rTo& V &\rTo& \ZZ(0) &\rTo& 0
\end{diagram}
\end{equation}
and points of $J_2(H)$.

Given an extension as in \eqref{eq:ext-MHS}, the underlying sequence of $\ZZ$-modules
\begin{diagram}
	0 &\rTo& \HZ &\rTo& \VZ &\rTo& \ZZ &\rTo& 0
\end{diagram}
splits non-canonically, and so we can find $\vZ \in \VZ$ mapping to $1 \in
\ZZ$; it is unique up to elements of $\HZ$. In the usual way of proceeding, one 
chooses a second lifting $v_F \in F_0 \VC$, using the surjectivity of $F_0 \VC \to
\CC$, and observes that $v_F - \vZ$ gives a well-defined point in $J_1(H)$. On the
other hand, only one choice is necessary to obtain a point in $J_2(H)$. First off, the fact that
$H$ is polarized implies that the dual Hodge structure $\Hd$ is isomorphic to
$H(-1)$. Dualizing \eqref{eq:ext-MHS}, we obtain a second exact sequence
\begin{diagram}
	0 &\rTo& \ZZ(0) &\rTo& \Vd &\rTo& H(-1) &\rTo& 0,
\end{diagram}
and the strictness of morphisms of Hodge structure gives $F_{-1} \VdC
\simeq F_0 \HC$. Now $\vZ$ defines a linear operator on $\VdC$, and hence on $F_{-1}
\VdC$; taking the ambiguity in choosing $\vZ$ into account, we therefore get a
well-defined point in the quotient
\[
	\frac{(F_0 \HC)^{\vee}}{\HZ} = J_2(H).
\]

\begin{lemma} \label{lem:J1-J2}
Under the isomorphism between $J_1(H)$ and $J_2(H)$ induced by $Q$, the two
constructions give rise to the same point.
\end{lemma}

\begin{proof}
If $A$ and $B$ are two mixed Hodge structures, then $F_p \Hom(A_{\CC},B_{\CC})$
consists of all maps $\phi \colon A_{\CC} \to B_{\CC}$ with $\phi(F_k A_{\CC})
\subseteq F_{k+p} B_{\CC}$. Therefore
\[
	F_{-1} \HdC = \menge{\psi \colon \HC \to \CC}{\psi(F_0 \HC) = 0},
\]
and the isomorphism with $F_0 \HC$ is given by taking $h \in F_0 \HC$ to the
functional $\psi_h = Q(h, \argbl)$. Similarly,
\[
	F_{-1} \VdC = \menge{\phi \colon \VC \to \CC}{\phi(F_0 \VC) = 0},
\]
which maps to $F_{-1} \HdC$ by restriction. 

Now $\vZ$ operates on $F_{-1} \VdC$ by taking $\phi \colon \VC \to \CC$ to $\phi(\vZ)$.
Given $h \in F_0 \HC$, let $\phi_h \in F_{-1} \VdC$ be the unique extension of $\psi_h
\in F_{-1} \HdC$. We compute that
\[
	\phi_h(\vZ) = \phi_h(\vZ - v_F) = \psi_h(\vZ - v_F) = Q(h, \vZ - v_F) 
		= Q(v_F - \vZ, h).
\]
But under the isomorphism $J_1(H) \simeq J_2(H)$, the point in $J_1(H)$
defined by $v_F - \vZ$ is exactly represented by the class of the linear map
\[
	F_0 \HC \to \CC, \quad h \mapsto Q(v_F - \vZ, h),
\]
and so the two constructions define the same point, as asserted.
\end{proof}

From now on, we shall use the second definition exclusively.
\begin{definition}
Let $H$ be a polarized integral Hodge structure of weight $-1$ with polarization $Q$. The
\define{intermediate Jacobian} of $H$ is the complex torus 
\[
	J(H) = (F_0 \HC)^{\vee} / \HZ,
\]
where the map $\HZ \into (F_0 \HC)^{\vee}$ is given by $h \mapsto Q(h, \argbl)$.
\end{definition}

This seems a good point to introduce a small generalization of the intermediate
Jacobian, which appears in the construction of \cite{BPS}. 

\begin{definition}
Let $H$ be an integral mixed Hodge structure of weight $\leq -1$. The
\define{generalized intermediate Jacobian} of $H$ is the complex Lie group
\[
	J(H) = (F_0 \PC)^{\vee} / \HZ,
\]
where $P = \Hd(1)$ is an integral mixed Hodge structure of weight $\geq -1$.
\end{definition}

Note that extensions as in \eqref{eq:ext-MHS} are still classified by the complex Lie
group
\[
	\Ext_{\MHS}^1 \bigl( \ZZ(0), H \bigr) \simeq \frac{\HC}{F_0 \HC + \HZ}.
\]
The reason for the above definition is that $\HC / F_0 \HC \simeq (F_0 \PC)^{\vee}$,
because $F_0 \PC = \menge{\phi \colon \HC \to \CC}{\phi(F_0 \HC) = 0}$.
Just as in the pure case, the generalized intermediate Jacobian parametrizes
extensions of mixed Hodge structure: an extension as in
\eqref{eq:ext-MHS} determines a point in $J(H)$ by using $\vZ \in \VZ$ as a linear
functional on $F_0 \PC$.

\subsection{Outline of the construction}
\label{subsec:construction-outline}

Let $\Xb$ be a complex manifold of dimension $n$, and let $X = \Xb - D$ be the
complement of a closed analytic subset. Let $\shH = \bigl( \shHO, \nabla, F_{\bullet}
\shHO, \shHZ, Q \bigr)$ be a polarized variation of Hodge structure of weight $-1$ on
$X$. To introduce some notation, we recall that this means the following: $\shHO$ is
a holomorphic vector bundle with a flat connection $\nabla$, and $\shHZ$ is a local
system of free $\ZZ$-modules such that $\ker \nabla \simeq \CC \tensor_{\ZZ} \shHZ$.
The $F_p \shHO$ are holomorphic subbundles of $\shHO$ that satisfy Griffiths'
transversality condition $\nabla \bigl( F_p \shHO \bigr) \subseteq \OmX{1} \tensor
F_{p+1} \shHO$.  Finally, $Q \colon \shHZ \tensor \shHZ \to \ZZ_X$ is alternating,
nondegenerate, flat for the connection $\nabla$, and satisfies $Q(F_p \shHO, F_q
\shHO) = 0$ if $p+q \leq 0$.

\begin{note}
Here and in what follows, we often consider the flat vector bundle $(\shHO, \nabla)$
as a special case of a left $\Dmod$-module; it is then more natural to write the Hodge
filtration as an increasing filtration, by setting $F_p \shHO = F^{-p} \shHO$. 
\end{note}

Each Hodge structure in the variation has its associated intermediate Jacobian
(defined as in \subsecref{subsec:Jacobians}); they fit together into a holomorphic
fiber bundle that we denote by $J(\shH) \to X$. By definition, its sheaf of
holomorphic sections is given by $(F_0 \shHO)^{\vee} / \shHZ$. To extend $J(\shH)$ to
a space over $\Xb$, we let $M$ be the polarized Hodge module on $\Xb$, obtained from
the variation $\shH$ by intermediate extension via the inclusion map $j \colon X
\into \Xb$ \cite{Saito-MHM}*{Theorem~3.21}. Then $M$ is a polarized Hodge module of
weight $n-1$ with strict support equal to all of $\Xb$. Its underlying perverse sheaf
$\rat M$ is simply the intersection complex of the local system $\QQ \tensor \shHZ$.

Let $(\Mmod, F)$ be the filtered left $\Dmod$-module underlying $M$. This means that
$\Mmod$ is a filtered holonomic $\DmodXb$-module, and $F = F_{\bullet} \Mmod$ is an
increasing filtration of $\Mmod$ by $\OXb$-coherent subsheaves that is good in the
sense of \cite{Borel}. The condition on the strict support implies that $\Mmod$ is
the minimal extension of the flat vector bundle $(\shHO, \nabla)$ from $X$ to $\Xb$.
The coherent sheaves $F_p \Mmod$ are natural extensions of the Hodge bundles;
in particular, $\ju(F_0 \Mmod) \simeq F_0 \shHO$.

As mentioned in \subsecref{subsec:summary}, it is sensible to define $\Jb(\shH)$ as
the quotient $T(F_0 \Mmod) / \TZ$, where $\TZ$ is the \'etal\'e space of the sheaf
$\jl \shHZ$, and $T(F_0 \Mmod)$ is as in \subsecref{subsec:spaces} below. Note that
both are analytic spaces, whose sheaves of sections are $\jl \shHZ$ and $(F_0
\Mmod)^{\vee}$, respectively. To carry this through, we first construct a holomorphic
map $\eps \colon \TZ \to T(F_0 \Mmod)$ that generalizes the embedding of the local
system $\shHZ$ into the vector bundle $T(F_0 \shHO)$. We then prove that the $\eps$
is a closed embedding, and that the fiberwise quotient $T(F_0 \Mmod) / \TZ$ is an
analytic space (in particular, Hausdorff), provided that the following condition is
satisfied:
\begin{condition} \label{eq:main-condition}
The map $\eps \colon \TZ \to T(F_0 \Mmod)$ is injective, and $\eps(\TZ)$ is a closed
analytic subset of $T(F_0 \Mmod)$.
\end{condition}
An important result of this paper is that Condition~\ref{eq:main-condition}
is true without assumptions on the complement $\Xb - X$ or on the local
monodromy of $\shHZ$.

\subsection{The analytic space associated to a coherent sheaf}
\label{subsec:spaces}

Let $X$ be an analytic space, and $\shF$ a coherent analytic sheaf on $X$. In this
section, we describe how to associate to $\shF$ an analytic space $T(\shF) \to X$,
relatively Stein, whose sheaf of sections is $\shFd = \shHom(\shF, \OX)$. We also
describe its most basic properties.

The construction of $T(\shF)$ is very simple: let $\Sym_{\OX}(\shF)$ be
the symmetric algebra in $\shF$, and define
\[
	T(\shF) = \Spec_X \bigl( \Sym_{\OX}(\shF) \bigr).
\]
When $\shF = \OX(E)$ is the sheaf of sections of a holomorphic vector bundle $E \to X$, we
recover the dual vector bundle since $T(\shF) = E^{\ast}$.
This leads to the following more concrete description of $T(\shF)$.
Let $j \colon U \into X$ be any open subset of $X$ that is Stein. Then $\ju \shF$ can
be written as a quotient of locally free sheaves on $U$,
\begin{diagram}
	\shE_1 &\rTo^{\varphi}& \shE_0 &\rTo& \ju \shF &\rTo& 0.
\end{diagram}
Let $\Ed_0 \to U$ be the holomorphic vector bundle whose sheaf of sections is
$\shEd_0$; similarly define $\Ed_1$. Then $\varphi$ induces a map of vector
bundles $\Ed_0 \to \Ed_1$, and $T(\ju \shF) \subseteq \Ed_0$ is the preimage of the
zero section. The reason is that $\Sym_{\OX}(\ju \shF)$ is the quotient of
$\Sym_{\OX}(\shE_0)$ by the ideal generated by $\varphi(\shE_1)$.

From the local description, it follows that $T(\shF) \to X$ is relatively Stein,
meaning that the preimage of every Stein open subset is again Stein; moreover,
every fiber is a linear space of some dimension, and over any analytic subset
of $X$ where the fiber dimension is constant, $T(\shF)$ is a holomorphic vector
bundle. The space has the following universal property.

\begin{lemma}
For any holomorphic map $f \colon Y \to X$ from an analytic space $Y$, 
\[
	\Map_X \bigl( Y, T(\shF) \bigr) \simeq 
		\Hom_{\OY} \bigl( \fu \shF, \OY \bigr).
\]
\end{lemma}

\begin{proof}
Holomorphic maps $Y \to T(\shF)$ over $X$ are in one-to-one correspondence with
morphisms of $\OX$-algebras $\Sym_{\OX}(\shF) \to \fl \OY$, hence with maps of
$\OX$-modules $\shF \to \fl \OY$, and therefore also with maps of $\OY$-modules
$\fu \shF \to \OY$.
\end{proof}

In particular, the sheaf of holomorphic sections of $T(\shF) \to X$ is exactly
$\shFd$. The next lemma shows that the construction of $T(\shF)$ behaves well under
pullback by arbitrary holomorphic maps. It follows that the fiber over
a point $x \in X$ is the dual of the vector space $\shF \tensor_{\OX} \shO_{X,x} /
\mathfrak{m}_x$.

\begin{lemma} \label{lem:T-functoriality}
For any holomorphic map $f \colon Y \to X$, we have
\[
	Y \times_X T(\shF) \simeq T(\fu \shF).
\]
\end{lemma}

\begin{proof}
This is true because $\fu \Sym_{\OX}(\shF) \simeq \Sym_{\OY}(\fu \shF)$, by the universal property
of the symmetric algebra.
\end{proof}

\begin{lemma} \label{lem:T-surjective}
Let $\shF \onto \shG$ be a surjective map of coherent sheaves. Then the induced map
$T(\shG) \to T(\shF)$ is a closed embedding.
\end{lemma}

\begin{proof}
The statement is local on $X$, and so we may assume without loss of generality that
$X$ is a Stein manifold. By writing $\shF$ as the quotient of a locally free sheaf
$\shE_0$, we can find compatible presentations
\begin{diagram}
\shE_1 &\rTo^{\varphi}& \shE_0 &\rTo& \shF &\rTo& 0 \\
\dTo & & \dEqual & & \dOnto \\
\shE_2 &\rTo^{\psi}& \shE_0 &\rTo& \shG &\rTo& 0.
\end{diagram}
Obviously, we now have $T(\shG) \subseteq T(\shF) \subseteq \Ed_0$, proving the
assertion.
\end{proof}

\begin{note}
Another analytic space with sheaf of sections $\shFd$ would be $T(\shF^{\vee\vee})$,
obtained by replacing $\shF$ by its double dual. Since the sheaf $\shF^{\vee\vee}$ is
reflexive, this may seem a more natural choice at first glance. But the problem is
that taking the dual does not commute with pullbacks by non-flat maps;
this second choice of space is therefore not sufficiently functorial for our purposes.
\end{note}

\subsection{Quotients of certain complex manifolds}
\label{subsec:quotients}

In this section, we show how conditions analogous to Condition~\ref{eq:main-condition}
allow one to prove that certain quotient spaces of holomorphic vector bundles are
again complex manifolds. The general situation is the following.
Let $p \colon E \to X$ be a holomorphic vector bundle on a complex manifold $X$.
Let $\shG$ be a sheaf of finitely generated abelian groups on $X$, and suppose that
we have a map of sheaves $\shG \to \OX(E)$. It defines a map of complex manifolds
$\eps \colon G \to E$, where $G$ is the \'etal\'e space of the sheaf $\shG$. We shall
require the following two conditions:
\begin{enumerate}[label=(\roman{*}), ref=(\roman{*})]
\item The image $\eps(G) \subseteq E$ is a closed analytic subset of $E$.
\label{en:quotients-i}
\item The map $\eps$ is injective.
\label{en:quotients-ii}
\end{enumerate}
For a point $x \in X$, we let $E_x = p^{-1}(x)$ and $G_x = \eps^{-1}(E_x)$ be the fibers. 
The second condition is equivalent to the injectivity of the individual maps
$G_x \to E_x$; note that $\eps(G_x)$ is then automatically a discrete subset of $E$, being both
closed analytic and countable.

\begin{lemma} \label{lem:quotients-1}
For any point $g \in G$, there is an open neighborhood of $\eps(g) \in E$ whose
intersection with $\eps(G)$ is the image of a local section of $G$.
\end{lemma}

\begin{proof}
As an analytic subset, $\eps(G)$ has a decomposition into (countably many) irreducible
components, and there is a small open neighborhood of $e = \eps(g)$ that meets only
finitely many of them. Shrinking that neighborhood, if necessary, we can find an open
set $U$ containing $e$, such that $\eps(G) \cap U$ has finitely many irreducible
components, each passing through the point $e$. Since $\eps$ is injective by
\ref{en:quotients-ii}, there can be only one such component $Z$; noting that $\eps(G_x)$
is discrete in $E$, we may further shrink $U$ and assume that $Z
\cap E_x = \{e\}$. For dimension reasons, we then have $\dim Z = \dim X$. Now $G$ is
the \'etal\'e space of the sheaf $\shG$, and so we can find a local section of $G$, 
defined in a suitable neighborhood $V$ of the point $x = p(e) \in X$, with the property
that $\gamma(x) = g$. It follows that $Z = \eps \bigl( \gamma(V) \bigr)$, as claimed.
\end{proof}

\begin{lemma}
The map $\eps \colon G \to E$ is a closed embedding.
\end{lemma}

\begin{proof}
First of all, $\eps$ is a proper map. To see this, let $g_n \in G$ be any sequence of
points in $G$ such that $\eps(g_n)$ converges to a point $e \in E$. By
\ref{en:quotients-i}, the limit is of the form $e = \eps(g)$ for some $g \in G$. By the
preceding lemma, there is an open neighborhood $U$ containing $e$, and a local section $\gamma
\colon V \to G$, such that $U \cap \eps(G) = \eps \bigl( \gamma(V) \bigr)$ and $g =
\gamma(x)$. We conclude that $g_n = \gamma(x_n)$ for some choice of $x_n \in V$. But
now $x_n = p \bigl( \eps(g_n) \bigr) \to x$, and therefore $g_n \to g$; this establishes
the properness of $\eps$. Lemma~\ref{lem:quotients-1}  also shows that $\eps \colon G
\to \eps(G)$ is a local isomorphism. Since $\eps$ is in addition proper and
injective, it has to be a closed embedding.
\end{proof}

The lemma justifies identifying $G$ with its image in $E$; from now on, we regard $G$
as a closed submanifold of $E$. We are then interested in taking the fiberwise
quotient of $E$ by $G$. Let $\sim$ be the equivalence relation on $E$ defined by
\[
	e \sim e' \quad \text{if and only if $p(e) = p(e')$ and $e - e' \in G$.}
\]
Let $q \colon E \to E / \sim$ be the map to the quotient, endowed with the quotient
topology.

\begin{lemma} \label{lem:quotients-3}
The map $q$ is open.
\end{lemma}

\begin{proof}
Let $U \subseteq E$ be any open set; we need to verify that $q^{-1} \bigl( q(U)
\bigr)$ is again open. It suffices to show that for any sequence of points $e_n$
that converges to some $e \in q^{-1} \bigl( q(U) \bigr)$, all but finitely many of the
$e_n$ also belong to $q^{-1} \bigl( q(U) \bigr)$. Since $q(e) \in q(U)$, there
exists $e' \in U$ with $e \sim e'$, hence $e' - e \in G$. Let $\gamma \colon V
\to G$ be a local section such that $e' = e + \gamma \bigl( p(e) \bigr)$. If we
put $e_n' = e_n + \gamma \bigl( p(e_n) \bigr)$, then $e_n' \to e'$, and so $e_n' \in
U$ for large $n$. But then $e_n \sim e_n'$ also belongs to $q^{-1} \bigl( q(U)
\bigr)$.
\end{proof}

\begin{lemma} \label{lem:quotients-4}
The quotient $E / \sim$ is Hausdorff.
\end{lemma}

\begin{proof}
Since $q$ is open, the quotient $E/\sim$ is Hausdorfff if and only if the equivalence
relation $\sim$ is closed in $E \times E$. Suppose that we have a sequence of points
$(e_n, e_n')$ with $e_n \sim e_n'$, such that $(e_n, e_n') \to (e,e') \in E \times
E$. Since $p$ is continuous, we deduce that $p(e) = p(e')$. But then $e_n'
- e_n \in G$ converges to $e'-e$, and because $G$ is closed, it follows that $e' - e
\in G$, and so $e' \sim e$. This proves that $\sim$ is indeed a closed subset of $E
\times E$.
\end{proof}

\begin{proposition} \label{prop:quotients-5}
If the two conditions in \ref{en:quotients-i} and \ref{en:quotients-ii} are
satisfied, then the quotient space $E / \sim$ is a complex manifold, and the map $q$
is holomorphic.
\end{proposition}

\begin{proof}
From Lemma~\ref{lem:quotients-1} and the fact that $q$ is open, it follows
that any sufficiently small open set on $E$ is mapped homeomorphically onto its image
in $E / \sim$, and thus can serve as a local chart on the quotient. Being
Hausdorff, $E / \sim$ is then a complex manifold, and the quotient map $q$ is
holomorphic by construction.
\end{proof}

\subsection{The construction of the quotient}
\label{subsec:construction}

In this section, we shall prove that the quotient $T(F_0 \Mmod) / \TZ$ is an analytic
space, provided that Condition~\ref{eq:main-condition} is satisfied.

We first explain how to embed the \'etal\'e space of the sheaf $\jl \shHZ$ into the
analytic space $T(F_0 \Mmod)$. On $X$, where we have a variation of Hodge structure
of weight $-1$, it is clear how to do this. To extend the embedding to all of $\Xb$,
we need to know that sections of $\jl \shHZ$ can act, via the polarization $Q$, on
arbitrary sections of the $\Dmod$-module $\Mmod$.

\begin{lemma} \label{lem:shHC-action}
Let $U \subseteq \Xb$ open, and let $h \in \Gamma \bigl( U, \jl \shHC
\bigr)$ and $\sigma \in \Gamma(U, \Mmod)$ be any two sections. Then the holomorphic function
$Q(h, \sigma)$ on $U \cap X$ extends holomorphically to all of $U$.
\end{lemma}

\begin{proof}
Restricting to $U$, we may assume that $U = \Xb$. Let $D = \Xb - X$, which we may
assume to be a divisor (the statement being trivial otherwise). 
Fix a section $h \in \Gamma \bigl( \Xb, \jl \shHC \bigr)$; note that $h$ is in
particular flat. Now consider the map of holonomic $\Dmod$-modules
\[
	Q(h, \argbl) \colon j^{-1} \Mmod \to \OX;
\]
by the adjointness between the functors $j^{-1}$ and $\jlreg$, it induces a map
\[
	Q(h, \argbl) \colon \Mmod \to \jlreg \OX,
\]
where sections of $\jlreg \OX = \OXb(\ast D)$ have at worst poles along $D$.
Let $\Mmod_h$ be the preimage of $\OXb$ under this map; it is again holonomic, and
clearly has the same restriction to $X$ as $\Mmod$ itself. Because $\Mmod$ is a
minimal extension, we have to have $\Mmod_h = \Mmod$, thus proving the assertion.
\end{proof}

\begin{note}
A more elementary proof goes as follows: Since $\Xb$ is a complex manifold, it
suffices to show that $Q(h, \argbl)$ extends over a general point of $D$. After
restricting to a curve that meets $D$ at a smooth point, and is non-characteristic
for $\Mmod$, we can therefore reduce to the case of a variation of Hodge structure on
$\dst$, where the statement is easily proved by looking at the canonical
extension.
\end{note}

Let $\TZ \to \Xb$ be the \'etal\'e space of the
sheaf $\jl \shHZ$; as a set, $\TZ$ is the union of all the stalks of the sheaf,
topologized to make every section continuous. For every point in $\TZ$, there is a
unique local section of $\jl \shHZ$ that passes through that point. By using such
local sections as charts, $\TZ$ acquires the structure of a complex manifold, making
the projection map and every section of the sheaf holomorphic. Note that the map
$\pZ \colon \TZ \to \Xb$ is locally an isomorphism, and therefore flat.

Lemma~\ref{lem:shHC-action} gives us a map of sheaves $\jl \shHZ \to (F_0
\Mmod)^{\vee}$, and therefore a holomorphic section of $\pZu (F_0 \Mmod)^{\vee} \simeq
(\pZu F_0 \Mmod)^{\vee}$ on $\TZ$. By the universal property of $T(F_0 \Mmod)$ in
Lemma~\ref{lem:T-functoriality}, this means that we have a holomorphic map 
\begin{equation} \label{eq:definition-eps}
	\eps \colon \TZ \to T(F_0 \Mmod)
\end{equation}
from the complex manifold $\TZ$ to the analytic space $T(F_0 \Mmod)$. 

\begin{proposition}
Assume that Condition~\ref{eq:main-condition} is satisfied. Then $\eps \colon \TZ \to
T(F_0 \Mmod)$ is a closed embedding.
\end{proposition}

\begin{proof}
The question is clearly local on $\Xb$; thus we may assume that $\Xb$ is a Stein
manifold. As explained in \subsecref{subsec:spaces}, we present $F_0 \Mmod$ as a quotient
of locally free sheaves, 
\begin{equation} \label{eq:presentation}
\begin{diagram}
	\shE_1 &\rTo& \shE_0 &\rTo& F_0 \Mmod &\rTo& 0,
\end{diagram}
\end{equation}
and let $\varphi \colon \Ed_0 \to \Ed_1$ be the corresponding map of vector bundles;
then $T(F_0 \Mmod) = \varphi^{-1}(0)$ is a closed analytic subset of $\Ed_0$. Because of
Condition~\ref{eq:main-condition}, the map from $\TZ$ to $\Ed_0$ satisfies the two conditions
in \subsecref{subsec:quotients}; we can now apply Lemma~\ref{lem:quotients-1}
to conclude that $\TZ \to \Ed_0$, and therefore also $\eps$ itself, is a closed embedding. 
\end{proof}

From now on, we identify $\TZ$ with its image in $T(F_0 \Mmod)$.  Next, we deduce
from the general results in \subsecref{subsec:quotients} that the quotient $T(F_0 \Mmod) /
\TZ$ is an analytic space.

\begin{proposition}
Assume that Condition~\ref{eq:main-condition} is satisfied. Then the fiberwise
quotient $T(F_0 \Mmod) / \TZ$ is an analytic space over $\Xb$.
\end{proposition}

\begin{proof}
This is again a local problem, and so we continue to assume that $\Xb$ is a Stein
manifold, and that $F_0 \Mmod$ has a presentation as in \eqref{eq:presentation}. Let
$p \colon \Ed_0 \to \Xb$ be the projection, and let $\sim$ be the equivalence
relation on $\Ed_0$ given by
\[
	e \sim e' \quad \text{if and only if $p(e) = p(e')$ and $e'-e \in \TZ$}.
\]
The quotient space $Y = \Ed_0 / \TZ = \Ed_0 / \sim$ is a complex manifold by
Proposition~\ref{prop:quotients-5}, and the quotient map $q \colon \Ed_0 \to Y$ is
holomorphic. (Note that the quotient is in particular Hausdorff, as proved in
Lemma~\ref{lem:quotients-4}.)

The map $\varphi \colon \Ed_0 \to \Ed_1$ takes the submanifold $\TZ$ into the zero
section of $\Ed_1$. This implies that we have a factorization $\varphi = \psi \circ
q$, with $\psi \colon Y \to \Ed_1$ holomorphic. Remembering that $T(F_0 \Mmod) =
\varphi^{-1}(0)$, we see that the quotient $T(F_0 \Mmod) / \TZ$ is naturally
identified with the closed subset $\psi^{-1}(0)$ of $Y$, and is thus an analytic
space as well.
\end{proof}

\subsection{The V-filtration and pullbacks of Hodge modules}
\label{subsec:V-filtration}

In this section, we briefly review the $V$-filtration, and then study the behavior of
the Hodge filtration under pullbacks of mixed Hodge modules. This will be used in
\subsecref{subsec:functoriality} below to prove the functoriality of our
construction.

Let $X$ be a complex manifold, and $Z \subseteq X$ a submanifold of codimension one.
We first look at the local setting where $Z$ is the zero locus of a holomorphic
function $t$; set $\partial = \vfelds{t}$. Let $I_Z = t \cdot \OX$ be the corresponding
ideal. Then
\[
	V^0 \Dmod_X = \menge{D \in \Dmod_X}{D \cdot I_Z \subseteq I_Z}.
\]
Now let $\Mmod$ be a left $\Dmod$-module on $X$. There is at most one decreasing
filtration $V = V^{\bullet} \Mmod$, indexed by $\QQ$, satisfying the following
conditions:
\newcounter{saveenum}
\begin{enumerate}[label=(\roman{*}), ref=(\roman{*})]
\item Each $V^{\alpha} M$ is a coherent $V^0 \Dmod_X$-module.
\item The filtration is exhaustive, meaning that $\Mmod = \bigcup_{\alpha}
V^{\alpha} \Mmod$, and left-con\-ti\-nu\-ous, meaning that $V^{\alpha} \Mmod =
\bigcap_{\beta < \alpha} V^{\beta} \Mmod$.
\item The filtration is discrete, meaning that any bounded interval contains only finitely
many $\alpha \in \QQ$ such that $\Gr_V^{\alpha} \Mmod = V^{\alpha} \Mmod / V^{>
\alpha} \Mmod$ is nonzero.
\item One has $t \cdot V^{\alpha} \Mmod \subseteq V^{\alpha+1} \Mmod$ and
$\partial \cdot V^{\alpha} \Mmod \subseteq V^{\alpha-1} \Mmod$.
\item For $\alpha \gg 0$, the filtration satisfies $V^{\alpha} \Mmod = t \cdot V^{\alpha-1} \Mmod$.
\item The operator $t \partial - \alpha + 1$ is nilpotent on $\Gr_V^{\alpha} \Mmod$.
\setcounter{saveenum}{\value{enumi}}
\end{enumerate}
When $\Mmod$ is regular and holonomic, M.~Kashiwara \cite{Kashiwara-V} and
B.~Malgrange \cite{Malgrange} have shown that such a filtration exists;
it is called the \define{$V$-filtration} of $\Mmod$, relative to the closed submanifold $Z$.
It is easy to see from the conditions that $t \colon V^{\alpha-1} \Mmod \to
V^{\alpha} \Mmod$ is an isomorphism for $\alpha > 1$, and that $\partial \colon
\Gr_V^{\alpha+1} \Mmod \to \Gr_V^{\alpha} \Mmod$ is an isomorphism for $\alpha \neq 0$.

Now consider the case when $(\Mmod, F)$ is a filtered $\Dmod$-module. In that case,
the $V$-filtration is said to be \define{compatible with $F$}, and $(\Mmod, F)$ is
called \define{quasi-unipotent and regular along $Z$} if, in addition to the above:
\begin{enumerate}[label=(\roman{*}), ref=(\roman{*})]
\setcounter{enumi}{\value{saveenum}}
\item For every $p$, one has $F_p V^{\alpha} \Mmod = t \cdot F_p V^{\alpha-1} \Mmod$,
provided that $\alpha > 1$.
\item For every $p$, one has $F_p \Gr_V^{\alpha} \Mmod = \partial \cdot
F_{p-1} \Gr_V^{\alpha+1} \Mmod$, provided that $\alpha < 0$.
\end{enumerate}
When $(\Mmod, F)$ is the filtered $\Dmod$-module underlying a polarized mixed Hodge
module, then the $V$-filtration exists and is compatible with $F$; moreover, each
$\Gr_V^{\alpha} \Mmod$, with the induced filtration, again underlies a mixed Hodge
module on $Z$. In fact, this is built into M.~Saito's definition
\cite{Saito-MHM}*{2.17} of the category of mixed Hodge modules.

The $V$-filtration is essential for the construction of nearby cycles, vanishing
cycles, and the various pullback operations on mixed Hodge modules \cite{Saito-MHM}. Suppose that $M$
is a mixed Hodge module on $X$, with underlying filtered
$\Dmod$-module $(\Mmod, F)$. To begin with, let $i \colon Z \into X$ be the inclusion
of a submanifold that is defined by a single holomorphic equation $t$. In this
situation, one can associate to $M$ two mixed Hodge modules on $Z$:
\begin{enumerate}[label=(\alph{*}), ref=(\alph{*})]
\item The (unipotent) nearby cycles $\psione{t} M$. Their underlying filtered
$\Dmod$-module is $\bigl( \Gr_V^1 \Mmod, F \bigr)$, where the Hodge filtration is induced
by that on $\Mmod$.
\item The vanishing cycles $\phione{t} M$. Their underlying filtered
$\Dmod$-module is given by $\bigl( \Gr_V^0 \Mmod, F \decal{-1} \bigr)$.
\end{enumerate}
The two standard maps $\can \colon \psione{t} M \to \phione{t} M$ and $\Var \colon
\phione{t} M \to \psione{t} M(-1)$ are morphisms of mixed Hodge modules; on the level
of $\Dmod$-modules, $\can$ is multiplication by $\partial$, and $\Var$ multiplication
by $t$. The axioms imply that $t \partial$ is nilpotent on $\Gr_V^1 \Mmod$; it
corresponds to $(2 \pi i)^{-1} N$, where $N$ is the logarithm of the monodromy around
$Z$ on the nearby cycles $\psione{t} M$.

The pullback $\iu M$ is an object in the derived category $D^b \MHM(Z)$; by
\cite{Saito-MHM}*{Corollary~2.24}, it is represented by the complex (in
degrees $-1$ and $0$)
\begin{diagram}[width=3em,midshaft]
\iu M = \Big\lbrack \psione{t} M &\rTo^{\can}& \phione{t} M \Big\rbrack \decal{1}.
\end{diagram}
Each cohomology module $H^k \iu M$ is again a mixed Hodge module on $M$,
nonzero only for $k=-1,0$. Note that pulling back does not increase weights: if $M$
has weight $\leq w$, then $H^k \iu M$ has weight $\leq w+k$
\cite{Saito-MHM}*{Proposition~2.26}. Analogously, $\ius M$ is
represented by the complex (in degrees $0$ and $1$)
\begin{diagram}[width=3em,midshaft]
\ius M = \Big\lbrack \phione{t} M &\rTo^{\Var}& \psione{t} M(-1) \Big\rbrack,
\end{diagram}
and $H^k \ius M$ has weight $\geq w+k$ if $M$ has weight $\geq w$.

% It is easy to see from the two complexes that we have a morphism of mixed Hodge
% modules $H^{-1} \iu M \to H^1 \ius M(1)$. When $M$ is pure of weight $w$, then the
% first module has weight $\leq w-1$ and the second weight $\geq w-1$, and the map
% induces an isomorphism $\Gr_{w-1}^W H^{-1} \iu M \simeq \Gr_{w-1}^W H^1 \ius M(1)$.

We now describe how the operation $\ius$ interacts with the Hodge filtration on the
underlying $\Dmod$-modules. 

\begin{lemma} \label{lem:map-ius-1}
Let $i \colon Z \into X$ be the inclusion of a submanifold, defined by a single
holomorphic equation $t$. Let $M^i = H^1 \ius M(1)$, and write $(\Mmod^i, F)$ for its
underlying filtered $\Dmod$-module on $Z$. 
\begin{enumerate}[label=(\roman{*}), ref=(\roman{*})]
\item There are canonical injective maps of coherent sheaves $F_p \Mmod^i \into \iu F_p \Mmod$.
\label{en:map-ius-1i}
\item If $M$ is smooth, then $M^i$ is the pullback of the corresponding
variation of mixed Hodge structure, and the map in \ref{en:map-ius-1i} is an
isomorphism.
\label{en:map-ius-1ii}
\end{enumerate}
\end{lemma}

\begin{proof}
Since $H^1 \ius M(1)$ is the cokernel of $\Var(1) \colon \phione{t} M(1) \to
\psione{t} M$, its underlying $\Dmod$-module $\Mmod^i$ is the cokernel of 
the map $\Gr_V^0 \Mmod \to \Gr_V^1 \Mmod$ given by multiplication by $t$.
Now $\Var$ is a morphism of mixed Hodge modules, and hence strict for $F$; this implies
that $F_p \Mmod^i$  is the cokernel of $F_p \Gr_V^0 \Mmod \to F_p \Gr_V^1 \Mmod$.
Equivalently,
\[
	F_p \Mmod^i = \frac{F_p V^1 \Mmod}{F_p V^{>1} \Mmod + t \cdot F_p V^0 \Mmod}
		= \frac{F_p V^1 \Mmod}{t \cdot F_p V^0 \Mmod},
\] 
where we have used the compatibility of $V$ with the Hodge filtration to conclude
that $F_p V^{>1} \Mmod = t \cdot F_p V^{>0} \Mmod$. We clearly have a map
\[
	\frac{F_p V^1 \Mmod}{t \cdot F_p V^0 \Mmod} \to 
		\frac{F_p \Mmod}{t \cdot F_p \Mmod},
\]
and since the quotient on the right is equal to $\iu F_p \Mmod$, we obtain half of
the assertion in \ref{en:map-ius-1i}. To show that the map is injective, it suffices
to prove that the intersection $V^1 \Mmod \cap (t \cdot F_p \Mmod)$ is contained
in $t \cdot F_p V^0 \Mmod$. So let $m \in V^1 \Mmod$, and suppose that $m = t m'$ for
some $m' \in F_p \Mmod$. Since the $V$-filtration is exhaustive, we can let $\alpha
\leq 0$ be the largest rational number with $m' \in V^{\alpha} \Mmod$. Now the
multiplication map
\[
	t \colon \Gr_V^{\alpha} \Mmod \to \Gr_V^{\alpha+1} \Mmod;
\]
is an isomorphism for $\alpha \neq 0$; since $t m' = m \in V^1 \Mmod$, we conclude
that $\alpha = 0$. Therefore $m \in t \cdot F_p V^0 \Mmod$, as desired.

As for \ref{en:map-ius-1ii}, note that when $M$ is smooth, $\Mmod$ is a flat vector
bundle. In that case, the $V$-filtration is essentially the $I_Z$-adic filtration, since $V^{\alpha}
\Mmod = I_Z^{\lceil \alpha \rceil -1} \Mmod$ (which equals $\Mmod$ if $\alpha \leq
0$), and so $\Gr_V^1 \Mmod = \iu \Mmod$, while $\Gr_V^0 \Mmod = 0$. It is then
immediate from the construction above that the map is an isomorphism.
\end{proof}

More generally, suppose that $i \colon Z \to X$ is the inclusion of a submanifold of
codimension $d$, defined by holomorphic equations $t_1, \dotsc, t_d$. In that case,
the functors $\iu$ and $\ius$ are obtained by iterating the construction above
\cite{Saito-MHM}*{p.~263}; thus $\ius M$ is the single complex
associated to the $d$-fold complex of mixed Hodge modules
\begin{equation} \label{eq:ius-def}
	\bigl( \phione{t_1} \xrightarrow{\Var} \psione{t_1}(-1) \bigr) \circ \dotsb 
		\circ \bigl( \phione{t_d} \xrightarrow{\Var} \psione{t_d}(-1) \bigr)(M).
\end{equation}
As before, we set $M^i = H^d \ius M(d)$, and denote the underlying filtered
$\Dmod$-module by $(\Mmod^i, F)$. Then $M^i$ is a quotient of the iterated nearby
cycles $\psione{t_1} \dotsb \psione{t_d} M$, and the statement of the previous lemma
continues to hold.

\begin{lemma} \label{lem:map-ius-d}
Let $i \colon Z \into X$ be the inclusion of a submanifold of codimension $d$,
defined by $d$ holomorphic equations $t_1, \dotsc, t_d$. Let $M^i = H^d \ius M(d)$,
and write $(\Mmod^i, F)$ for its underlying filtered $\Dmod$-module on $Z$. 
\begin{enumerate}[label=(\roman{*}), ref=(\roman{*})]
\item There are canonical maps of coherent sheaves $\beta^i \colon F_p \Mmod^i \to \iu F_p \Mmod$.
\label{en:map-ius-di}
\item If $M$ is smooth, then $M^i$ is the pullback of the corresponding
variation of mixed Hodge structure, and the map $\beta^i$ is an isomorphism.
\label{en:map-ius-dii}
\end{enumerate}
\end{lemma}

\begin{proof}
Arguing by induction on the codimension, we may suppose that we have already
constructed the map
\begin{equation} \label{eq:map-ius-d}
	\beta^{i_1} \colon F_p \Mmod^{i_1} \to \iku{1} F_p \Mmod,
\end{equation}
where $Z_1$ is the submanifold defined by $t_2 = \dotsb = t_d = 0$, and $i_1 \colon
Z_1 \into X$ is the inclusion. Clearly, $\Mmod^{i_1}$ underlies the mixed Hodge
module $M^{i_1} = H^{d-1} \ikus{1} M(d-1)$. Then $Z$ is of codimension one in $Z_1$,
and we let $i_0 \colon Z \into Z_1$ be the inclusion map. Since $\ikus{0} \circ
\ikus{1} \simeq \ius$, we have a spectral sequence
\[
	E_2^{p,q} = H^p \ikus{0} H^q \ikus{1} M \Longrightarrow H^{p+q} \ius M.
\]
From the complex in \eqref{eq:ius-def}, it is clear that $H^p \ikus{0} = 0$ unless
$p=0,1$, and $H^q \ikus{1} M = 0$ unless $0 \leq q \leq d-1$. Therefore the spectral
sequence degenerates, and 
\[
	H^1 \ikus{0} M^{i_1}(1) \simeq H^1 \ikus{0} H^{d-1} \ikus{1} M(d) \simeq 
		H^d \ius M(d) = M^i.
\]
Lemma~\ref{lem:map-ius-1}, applied to $M^{i_1}$, thus produces $F_p \Mmod^i \to
\iku{0} F_p \Mmod^{i_1}$. Compose this with the map $\iku{0} F_p \Mmod^{i_1} \to
\iku{0} \iku{1} F_p \Mmod \simeq \iu F_p \Mmod$ derived from \eqref{eq:map-ius-d}
to get \ref{en:map-ius-di}. The assertion in
\ref{en:map-ius-dii} follows directly from Lemma~\ref{lem:map-ius-1}.
\end{proof}

The maps $\beta^i$ in Lemma~\ref{lem:map-ius-d} are independent of the choice of
equations $t_1, \dotsc, t_d$ for $Z$; this follows from M.~Saito's proof
\cite{Saito-MHM}*{p.~259} that the functor $\ius$ is well-defined---in fact, the single complex
associated to \eqref{eq:ius-def} is a well-defined object in the derived category
of mixed Hodge modules on $Z$.

\begin{lemma} \label{lem:map-ius-indep}
The map $\beta^i \colon F_p \Mmod^i \to \iu F_p \Mmod$ depends only on the inclusion
of the submanifold $i \colon Z \into X$, but not on the choice of generators
for the ideal $I_Z$.
\end{lemma}

We now consider the functor $\fus$ for a general map $f \colon Y \to X$ between complex
manifolds. 

\begin{proposition} \label{prop:map-fus}
Let $f \colon Y \to X$ be a holomorphic map of complex manifolds, and let $M$ be a
mixed Hodge module on $X$ with underlying filtered $\Dmod$-module $(\Mmod, F)$.
Define $N^f = H^{d_X-d_Y} \fus M(d_X - d_Y)$, and denote its underlying
filtered $\Dmod$-module by $(\Nmod^f, F)$. Then we have canonical maps
\[
	\beta^f \colon F_p \Nmod^f \to \fu F_p \Mmod,
\]
which are isomorphisms if $M$ is smooth.
\end{proposition}

\begin{proof}
Let $d = d_X - d_Y$. To compute $\fus$, one factors $f$ as
\begin{diagram}
	Y &\rTo^i& W &\rTo^q& X,
\end{diagram}
with $q$ smooth of relative dimension $k$, and $i$ a closed embedding of codimension
$k+d$. Then $H^d \fus M(d) = H^{k+d} \ius H^{-k} \qus M(d) = H^{k+d} \ius N^q(k+d)$,
where we have set $N^q = H^{-k} \qus M(-k) = H^k \qu M \simeq \QQ_Y^H \decal{k} \boxtimes M$.
Evidently, the filtered $\Dmod$-module underlying $N^q$ is $(\Nmod^q, F) = (\qu
\Mmod, \qu F)$, with the pullbacks taken in the category of quasi-coherent sheaves;
in particular, $F_p \Nmod^q = \qu F_p \Mmod$. On the other hand, we have
$N^f = H^{k+d} \ius N^q(k+d)$, and so Lemma~\ref{lem:map-ius-d} shows that there is a
canonical map $\beta^i \colon F_p \Nmod^f \to \iu F_p \Nmod^q$. Compose this and the isomorphism
$\iu F_p \Nmod^q \simeq \fu F_p \Mmod$ to obtain the desired map. For smooth $M$, the
map is an isomorphism because of Lemma~\ref{lem:map-ius-d}. It remains to prove that
the map we have constructed is independent of the factorization $f = q \circ i$; this
is the content of the following lemma.
\end{proof}

\begin{lemma}
Let $f = q_1 \circ i_1 = q_2 \circ i_2$ be two factorizations of $f \colon Y \to X$
into a closed embedding $i_j \colon Y \into W_j$ and a smooth morphism $q_j \colon
W_j \to X$. Then the two resulting maps $F_p \Nmod^f \to \fu F_p \Mmod$ are equal.
\end{lemma}

\begin{proof}
Let $W = W_1 \times_Y W_2$ be the fiber product; both projections $p_j \colon W \to
W_j$ are smooth. Because of the commutative diagram
\begin{diagram}[width=2.5em]
&& W_1 \\
 & \ruTo^{i_1} & \uTo_{p_1} & \rdTo^{q_1} \\
Y &\rTo^i& W &\rTo^q& X \\
 & \rdTo^{i_2} & \dTo_{p_2} & \ruTo^{q_2} \\
&& W_2,
\end{diagram}
it suffices to show that the factorizations $q_j \circ i_j$ both give the same map as
$q \circ i$. The construction in Proposition~\ref{prop:map-fus} is clearly
insensitive to factorizing $q = q_2 \circ p_2$ since all three maps are smooth; this
reduces the problem to considering the maps $i_j = p_j \circ i$.

We may thus assume, without loss of generality, that $f \colon Y \to X$ is a closed
embedding. Let $d$ be the codimension of $Y$ in $X$, and $r$ the relative dimension
of the map $q \colon W \to X$. Then what we need to prove is the commutativity of
\begin{equation} \label{eq:map-fus-indep}
\begin{diagram}
F_p \Nmod^f &\rTo^{\beta^f}& \fu F_p \Mmod \\
\dTo^{\beta^i} && \uTo^{\simeq} \\
\iu F_p \Nmod^q &\rEqual& \iu \qu F_p \Mmod,
\end{diagram}
\end{equation}
with $\beta^i$ and $\beta^f$ as in Lemma~\ref{lem:map-ius-d}. This is a local
question; we may therefore assume that $X = Y \times \Delta^d$ and $W = Y \times
\Delta^d \times \Delta^r$. After factorizing the closed embeddings (which is permissible by
Lemma~\ref{lem:map-ius-indep}), we eventually reduce the whole problem to the special
case $Y = X$ and $W = Y \times \Delta$, with $i \colon Y \into Y \times \Delta$ the
inclusion and $q \colon Y \times \Delta \to Y$ the projection. Evidently, $\beta^f$
is now the identity map. Let $t$ be the coordinate function on $\Delta$. One easily
checks that the $V$-filtration on $\Nmod^q = \qu \Mmod$ is given by $V^{\alpha}
\Nmod^q = t^{\lceil \alpha \rceil -1} \qu \Mmod$ (by which we mean $\qu \Mmod$ if
$\alpha \leq 0$); then $\Gr_V^1 \Nmod^q = \qu \Mmod / t \cdot \qu \Mmod = \iu \qu
\Mmod \simeq \Mmod$, and so the diagram in \eqref{eq:map-fus-indep} does
commute as asserted.
\end{proof}

\subsection{Functoriality}
\label{subsec:functoriality}

In this section, we prove that our construction of the space $\Jb(\shH)$ is
functorial, in a sense made precise below.

Let $f \colon \Yb \to \Xb$ be a map of complex manifolds, such that $Y = f^{-1}(X)$
is dense in $\Yb$ (we also write $f \colon Y \to X$ for the induced map). As above, let $\shH$ be
a polarized variation of Hodge structure of weight $-1$ on $X$, let $M$ be the
polarized Hodge module on $\Xb$ obtained by intermediate extension, and $(\Mmod, F)$
its underlying filtered $\Dmod$-module. We denote the pullback of the variation of
Hodge structure by $\shH' = \fu \shH$, its intermediate extension to $\Yb$ by $M'$, and the
underlying filtered $\Dmod$-module by $(\Mmod', F)$.

\begin{lemma} \label{lem:F0-functoriality}
We have a canonical map of coherent sheaves
\[
	F_0 \Mmod' \to \fu F_0 \Mmod,
\]
whose restriction to $Y$ is the obvious isomorphism of Hodge bundles.
\end{lemma}

\begin{proof}
Let $n = \dim X$ and $m = \dim Y$, and note that $M$ has weight $n-1$. Since the
functor $\fus$ does not decrease weights, the mixed Hodge module 
\[
	N^f = H^{n-m} \fus M(n-m)
\]
has weight $\geq m-1$. The pure Hodge module $W_{m-1} N^f$ is therefore a
submodule of $N^f$. The restriction of $W_{m-1} N^f$ to $Y$ is canonically isomorphic to
the variation of Hodge structure $\shH'$; in the decomposition by strict support,
the component with strict support $\Yb$ has to be isomorphic to $M'$. Since the
decomposition is canonical, we get a uniquely defined map $M' \into W_{m-1} N^f \into
N^f$. Passing to the Hodge filtrations on the underlying $\Dmod$-modules, we thus
have a canonical map of coherent sheaves
$F_0 \Mmod' \into F_0 \Nmod^f$. We compose this with the map $\beta^f \colon F_0 \Nmod^f
\to \fu F_0 \Mmod$ in Proposition~\ref{prop:map-fus} to get the first half of the
assertion; the second follows directly from Proposition~\ref{prop:map-fus}.
\end{proof}

\begin{proposition} \label{prop:Jb-functoriality}
Let $f \colon \Yb \to \Xb$ be a morphism of complex manifolds, such that $Y =
f^{-1}(X)$ is dense in $\Yb$. If we let $\shH' = \fu \shH$ be the pullback of
the variation of Hodge structure $\shH$ from $X$ to $Y$, we have a
canonical holomorphic map
\[
	\Yb \times_{\Xb} \Jb(\shH) \to \Jb(\shH')
\]
over $\Yb$, whose restriction to $Y$ is the evident isomorphism between the
two families of intermediate Jacobians.
\end{proposition}

\begin{proof}
First consider the spaces $T(F_0 \Mmod)$ and $T(F_0 \Mmod')$ that appear in
the construction of $\Jb(\shH)$ and $\Jb(\shH')$. By
Lemma~\ref{lem:T-functoriality}, $\Yb \times_{\Xb} T(F_0 \Mmod) \simeq T(\fu F_0
\Mmod)$.  On the other hand, Lemma~\ref{lem:F0-functoriality} provides us with a
canonical map $T(\fu F_0 \Mmod) \to T(F_0 \Mmod')$. Composing the two, we obtain a
canonical holomorphic map
\[
	\Yb \times_{\Xb} T(F_0 \Mmod) \to T(F_0 \Mmod')
\]
over $\Yb$; over $Y$, the left-hand side restricts to the pullback of the vector
bundle associated with $(F^0 \shHO)^{\vee}$, the right-hand side to the vector
bundle associated with $(F^0 \shHO')^{\vee}$, and the map to the obvious
isomorphism between them. Since $\Yb \times_{\Xb} \TZ$ is easily seen to map into
$\TZ'$, we get the assertion for the quotient spaces as well.
\end{proof}

\subsection{Restriction to points}
\label{subsec:points}

In this section, we describe how $T(F_0 \Mmod)$ behaves upon restriction to points, by
relating its fibers to Hodge-theoretic information. Let $i \colon \{x\} \to \Xb$ be
the inclusion of a point.  Define the mixed Hodge structure $H = H^{-n} \iu M$, which
has weight $\leq -1$. 

\begin{lemma} \label{lem:H-integral}
The mixed Hodge structure $H = H^{-n} \iu M$ is naturally defined over $\ZZ$, with $\HZ$ isomorphic
to the stalk of the sheaf $\jl \shHZ$ at the point $x$. Consequently, $\HZ$ 
embeds into the stalk of $\shHZ$ at any nearby point $x_0 \in X$, and the
quotient $\shH_{\ZZ,x_0} / \HZ$ is torsion-free.
\end{lemma}

\begin{proof}
There is a natural map from the stalk of the sheaf $\jl \shHC$ to $\HC$,
\[
	\lim_{U \ni x} H^0 \bigl( U \cap X, \shHC \bigr) \to \HC,
\]
given as follows: Let $t_1, \dotsc, t_n$ be local holomorphic coordinates centered at
$x$, and $\partial_j = \vfelds{t_j}$; then a local section of $\jl \shHC$ is a
section $s \in H^0(U, \Mmod)$ that satisfies $\partial_1 s = \dotsb = \partial_n s = 0$. On
the other hand, if $V_j$ denotes the $V$-filtration relative to $t_j=0$, then $\HC$ is
by construction a subspace of $\Gr_{V_1}^1 \dotsb \Gr_{V_n}^1 \Mmod$, consisting of
those elements $h$ for which each $\partial_j h$ is zero in $\Gr_{V_1}^1 \dotsb
\Gr_{V_j}^0 \dotsb \Gr_{V_n}^1 \Mmod$. It is easy to see from the axioms in
\subsecref{subsec:V-filtration} that $\partial_n s = 0$ implies $s \in V_n^1
\Mmod$; for similar reasons, $s$ defines a point in $\Gr_{V_1}^1 \dotsb \Gr_{V_n}^1 \Mmod$, and hence
a point in $\HC$. The resulting map is an isomorphism; this means that the mixed
Hodge structure $H$ is defined over $\ZZ$, with integral lattice $\HZ$ isomorphic to
the stalk of $\jl \shHZ$. 

Now let $U \subseteq \Xb$ be a small open ball around $x$, and $x_0 \in U \cap X$.
The stalk of $\jl \shHZ$ at $x$ is naturally identified with the subgroup of
$\shH_{\ZZ,x_0}$ consisting of classes invariant under the action by $\pi_1(U \cap X,
x_0)$. Since $\shH_{\ZZ,x_0}$ is torsion-free, it is then easy to deduce the second
assertion.
\end{proof}

The following result will be used in two places: to reduce the proof of
Condition~\ref{eq:main-condition} from the general case to the normal crossing case
(in \subsecref{subsec:reduction}); and to relate the space $\Jb(\shH)$ to the N\'eron
model of \cite{BPS} (in \subsecref{subsec:BPS}).

\begin{lemma} \label{lem:restriction-points}
Let $i \colon \{x\} \to \Xb$ be the inclusion of a point, and set $H = H^{-n} \iu M$
and $P = H^n \ius M(n) \simeq \Hd(1)$. Then the canonical map $F_0 \PC \to \iu
F_0 \Mmod$ induces a surjective map of complex Lie groups
\[
	\Jb(\shH)_x \onto J(H),
\]
where $J(H)$ denotes the generalized intermediate Jacobian for $H$ (see
\subsecref{subsec:Jacobians}).
\end{lemma}

\begin{proof}
The mixed Hodge structure $P = H^n \ius M(n)$ has weight $\geq -1$. By
duality, $\DD \bigl( \iu M \bigr) \simeq \ius \DD(M) \simeq \ius M(-1)$, and
therefore $P \simeq \Hd(1)$; this shows that $P$ is also defined over $\ZZ$. The
pairing between $\HC$ and $\PC$ is induced by the polarization $Q$, as follows: By
Lemma~\ref{lem:shHC-action}, a section $h \in H^0 \bigl( U \cap X, \shHC \bigr)$
defines a map of $\Dmod$-modules $Q(h, \argbl) \colon \Mmod \restr{U} \to \shO_U$.
Since a map of $\Dmod$-modules automatically respects the $V$-filtration, it induces
a map
\[
	Q(h, \argbl) \colon \Gr_{V_1}^1 \dotsb \Gr_{V_n}^1 \Mmod \restr{U} \to 
		\Gr_{V_1}^1 \dotsb \Gr_{V_n}^1 \shO_U
\]
By construction, $\PC$ is a quotient of the mixed Hodge structure on the left (by
$t_1, \dotsc, t_n$), and therefore we obtain a linear map $Q(h, \argbl) \colon \PC
\to \CC$.

Now Proposition~\ref{prop:map-fus} gives us a map of vector spaces
\begin{equation} \label{eq:map-PC}
	F_0 \PC \to \iu F_0 \Mmod.
\end{equation}
From the description above, it is clear that \eqref{eq:map-PC} is compatible with the
action by $\HC$. By Lemma~\ref{lem:T-functoriality}, the fiber of $T(F_0 \Mmod)$ at
the point $x$ is exactly $(\iu F_0 \Mmod)^{\vee}$, and so \eqref{eq:map-PC} induces a
linear map
\[
	T(F_0 \Mmod)_x \simeq (\iu F_0 \Mmod)^{\vee} \to (F_0 \PC)^{\vee}.
\]
We observe that this map is surjective: indeed, the map $\HC \to (F_0 \PC)^{\vee}$ is
trivially surjective, and by the discussion above, it factors through $T(F_0 \Mmod)_x$.
If we now take the quotient by $\TZx \simeq \HZ$, we arrive at a surjective map
\[
	\frac{T(F_0 \Mmod)_x}{\TZx} \onto \frac{(F_0 \PC)^{\vee}}{\HZ}
\]
of complex Lie groups, as asserted.
\end{proof}

\subsection{Restriction to curves}
\label{subsec:curves}

In this section, we investigate how $T(F_0 \Mmod)$ behaves upon restriction to
curves, and use the result to show that the subset $\TZ$ is closed under limits along
analytic arcs. Throughout, we let $f \colon \Delta \to \Xb$ be a holomorphic map
such that $f(\dst) \subseteq X$, and set $x = f(0)$; the understanding is that $x
\in \Xb - X$.

We define $\shH'$ to be the pullback of the variation of Hodge structure $\shH$ to
$\dst$. If we let $M'$ be its intermediate extension, then $M'$ is a polarized Hodge
module of weight $0$ on $\Delta$; as usual, we shall denote its underlying filtered
$\Dmod$-module by $(\Mmod', F)$. We also let $N^f = H^{n-1} \fus M(n-1)$, which is a
mixed Hodge module of weight $\geq 0$. As in Lemma~\ref{lem:F0-functoriality},
decomposition by strict support means that we have canonical maps $M' \into W_0 N^f
\into N^f$, and consequently, a map of coherent sheaves $F_0 \Mmod' \to F_0 \Nmod^f$.
Proposition~\ref{prop:map-fus} gives us $F_0 \Nmod^f \to \fu F_0 \Mmod$, and so we
have two holomorphic maps
\begin{equation} \label{eq:Delta-holo}
	\Delta \times_{\Xb} T(F_0 \Mmod) \to T(F_0 \Nmod^f) \to T(F_0 \Mmod')
\end{equation}
of analytic spaces over $\Delta$.

We now study the fibers of the three spaces over $0 \in \Delta$. To begin with, let
$i \colon \{x\} \to \Xb$, and define the mixed Hodge structures $H = H^{-n} \iu M$
and $P = H^n \ius M(n) \simeq \Hd(1)$ as in \subsecref{subsec:points}. Recall that
$H$ is of weight $\leq -1$ and defined over $\ZZ$, with integral lattice $\HZ$
isomorphic to the stalk of $\jl \shHZ$ at the point $x$; on the other hand, $P$ is of
weight $\geq -1$. As in Lemma~\ref{lem:restriction-points}, we have a canonical map
\[
	T(F_0 \Mmod)_x \to (F_0 \PC)^{\vee} \simeq \HC / F_0 \HC.
\]
Similarly, let $i_0 \colon \{0\} \to \Delta$, and define $H' = H^{-1} \iku{0} M'$
(of weight $\leq -1$) and $P' = H^1 \ikus{0} M'(1)$ (of weight $\geq -1$); we also have a
second map
\[
	T(F_0 \Mmod')_0 \to \HC' / F_0 \HC'.
\]

To get information about the mixed Hodge module $N^f$, we note that $\ikus{0} \circ
\fus \simeq \ius$. This means that there is a spectral sequence
\[
	E_2^{p,q} = H^p \ikus{0} H^q \fus M \Longrightarrow H^{p+q} \ius M.
\]
Because $\Delta$ is one-dimensional, $H^p \ikus{0} = 0$ unless $p=0,1$; therefore the
spectral sequence degenerates, and we find that $H^1 \ikus{0} H^{n-1} \fus M \simeq
H^n \ius M$, using that $H^n \fus M = 0$. Consequently,
\[
	H^1 \ikus{0} N^f(1) \simeq H^1 \ikus{0} H^{n-1} \fus M(n) \simeq H^n \ius M(n) = P,
\]
and as before, this leads to a linear map
\[
	T(F_0 \Nmod^f)_0 \to \HC / F_0 \HC.
\]

Since the various maps we produce are compatible, we now end up with the following
commutative diagram that relates the fibers of the analytic spaces in
\eqref{eq:Delta-holo} to the mixed Hodge structures $H$ and $H'$:
\begin{equation} \label{eq:Delta-fibers}
\begin{diagram}
T(F_0 \Mmod)_x &\rTo& T(F_0 \Nmod^f)_0 &\rTo& T(F_0 \Mmod')_0 \\
\dTo && \dTo && \dTo \\
\HC / F_0 \HC &\rEqual& \HC / F_0 \HC &\rTo& \HC' / F_0 \HC'
\end{diagram}
\end{equation}

We can use the discussion above to show that $\eps(\TZ) \subseteq T(F_0 \Mmod)$ is closed
under limits along analytic curves, in the following sense.

\begin{lemma} \label{lem:restriction-curves}
Let $g \colon \Delta \to T(F_0 \Mmod)$ be a holomorphic map with the property that
$g(\dst) \subseteq \eps(\TZ) \cap p^{-1}(X)$, where $p \colon T(F_0 \Mmod) \to \Xb$
is the projection. Then we actually have $g(\Delta) \subseteq \eps(\TZ)$.
\end{lemma}

\begin{proof}
Set $f = p \circ g$, and let $\shH' = \fu \shH$ be the pullback of the variation to
$\dst$; we also use the other notation introduced above. Since $g(\dst) \subseteq
\eps(\TZ)$, it corresponds to an integral section 
\[
	h' \in \HZ' \simeq H^0 \bigl( \Delta, \shHZ' \bigr).
\]
Over $\dst$, the two spaces $\Delta \times_{\Xb} T(F_0 \Mmod)$ and $T(F_0 \Mmod')$
are isomorphic to the dual of a Hodge bundle, and the section of this bundle defined
by $g$ is nothing but $Q'(h', \argbl)$. By Lemma~\ref{lem:shHC-action}, it extends
to a holomorphic section of $T(F_0 \Mmod')$ over $\Delta$.

Using \eqref{eq:Delta-fibers}, the value $g(0) \in T(F_0 \Mmod)_x$ determines a
point in the quotient $\HC / F_0 \HC$. Since the point in $\HC' / F_0 \HC'$ coming
from the section $Q'(h', \argbl)$ is evidently the image of $h'$, the
commutativity of the diagram implies that
\[
	h' \in F_0 \HC' + \im(\HC \to \HC').
\]
Now $H \into H'$ is a morphism of mixed Hodge structures; let $H''$ be the quotient,
still a mixed Hodge structure of weight $\leq -1$. The image of $h'$ in $\HC''$ is
both rational (since $h'$ is) and in $F_0 \HC''$, and therefore equal to zero;
consequently, $h' \in \HQ$. Now Lemma~\ref{lem:H-integral} implies that we
automatically have $h' \in \HZ$: indeed, $\HZ$ and $\HZ'$ are both subgroups of the
stalk of $\shHZ$ at some nearby point $f(t_0)$, and $\shH_{\ZZ,f(t_0)} / \HZ$ has
no torsion.  

But then $Q(h', \argbl)$ defines a holomorphic section of
$\eps(\TZ) \subseteq T(F_0 \Mmod)$ in a neighborhood of the point $x$. Over the image
of $\dst$, it is an extension of $g$; since both are holomorphic, this means that
$g(\Delta) \subseteq \eps(\TZ)$, as claimed.
\end{proof}

\subsection{Reduction to the normal crossing case}
\label{subsec:reduction}

This section is devoted to reducing the proof of Condition~\ref{eq:main-condition} to
the case of a divisor with normal crossings and unipotent local monodromy.
Evidently, the problem is local on $\Xb$, and so we may assume that $\Xb = \Delta^n$
is a polydisk, and that $\Xb - X$ is a divisor (possibly singular and with several
components). Let $\shH$ be a polarized variation of Hodge structure of weight $-1$ on
$X$. Recall that $\TZ$ is the \'etal\'e space of the sheaf $\jl \shHZ$, and that we
had constructed a holomorphic map $\eps \colon \TZ \to T(F_0 \Mmod)$ in
\eqref{eq:definition-eps}, using the polarization.

We begin by showing that $\eps$ is injective. For this, it is clearly sufficient to
prove that the map on fibers, $\TZx \to T(F_0 \Mmod)_x$, is injective. This follows
rather easily from the results of \subsecref{subsec:points}, as follows.

\begin{lemma}
For $x \in X$, let $\TZx = \pZ^{-1}(x)$ and $T(F_0 \Mmod)_x = p^{-1}(x)$ denote the
fibers of $\TZ$ and $T(F_0 \Mmod)$, respectively. Then $\eps$ is injective, and
embeds $\TZx$ into $T(F_0 \Mmod)_x$ as a discrete subset.
\end{lemma}

\begin{proof}
Let $i \colon \{x\} \to X$ be the inclusion of the point, and let $H = H^{-n} \iu M$, which
is an integral mixed Hodge structure of weight $\leq -1$ with $\HZ \simeq \TZx$. Also
define the mixed Hodge structure $P = H^n \ius M(n)$ (of weight $\geq -1$); it
satisfies $P \simeq \Hd(1)$. According to the discussion in
\subsecref{subsec:points}, have a surjective linear map
\[
	T(F_0 \Mmod)_x \onto (F_0 \PC)^{\vee} \simeq \HC / F_0 \HC.
\]
But since $H$ has weight $\leq -1$, the set of integral points $\HZ$ maps injectively
and hence discretely into $\HC / F_0 \HC$. Consequently, the map $\eps$ also embeds $\TZx$
into $T(F_0 \Mmod)_x$ as a discrete subset, proving the assertion.
\end{proof}

For the remainder of this section, we assume that Condition~\ref{eq:main-condition}
is satisfied whenever $\Xb - X$ is a divisor with normal crossings and $\shH$ has
unipotent local monodromy. We show that it then holds in general.

\begin{lemma}
The closure of $\eps(\TZ)$ in $T(F_0 \Mmod)$ is an analytic subset.
\end{lemma}

\begin{proof}
Since the underlying local system $\shHZ$ is defined over $\ZZ$, the local
monodromy is at least quasi-unipotent by a theorem due to A.~Borel \cite{Schmid}*{Lemma~4.5}. 
Taking a finite branched cover, unbranched over $X$, and resolving singularities, we
construct a proper holomorphic
map $f \colon \Yb \to \Xb$ from a complex manifold $\Yb$ of dimension $n$, with the
following properties: $Y = f^{-1}(X)$ is dense in $\Yb$; the restriction of $f$ to
$Y$ is finite and \'etale; the complement $\Yb - Y$ is a divisor with normal
crossings; and the pullback of $\shHZ$ to $Y$ has unipotent monodromy.

Let $\shH' = \fu \shH$, let $M'$ be its intermediate extension to a polarized Hodge
module on $\Yb$, and $(\Mmod', F)$ the underlying filtered $\Dmod$-module. By
Lemma~\ref{lem:F0-functoriality}, we have a commutative diagram of holomorphic maps
\begin{diagram}[midshaft,width=4em]
	\TZ &\lTo& \Yb \times_{\Xb} \TZ &\rTo& \TZ' \\
	\dTo^{\eps} && \dTo^{\id \times \eps} && \dTo^{\eps'} \\
	T(F_0 \Mmod) &\lTo& \Yb \times_{\Xb} T(F_0 \Mmod) &\rTo^{\Phi}& T(F_0 \Mmod').
\end{diagram}
Assuming Condition~\ref{eq:main-condition} for the variation $\shH'$, we know that
$\eps'(\TZ')$ is a closed analytic subset of $T(F_0 \Mmod')$. Then $\Phi^{-1} \bigl(
\eps'(\TZ') \bigr)$ is a closed analytic subset of $\Yb \times_{\Xb} T(F_0 \Mmod)$.
The projection to $T(F_0 \Mmod)$ is proper, since $f$ is a proper map, and so the
image of $\Phi^{-1} \bigl( \eps'(\TZ') \bigr)$ in $T(F_0 \Mmod)$ is again a closed
analytic subset by Grauert's theorem. The part of it that lies over $X$ is equal to
$\eps(\TZ) \cap p^{-1}(X)$, and so the closure of $\eps(\TZ)$ must be an analytic subset
(and, in fact, one of the components of the image). 
\end{proof}

To conclude the reduction to the normal crossing case, we use the results about
restriction to curves in \subsecref{subsec:curves} to show that taking the closure
does not actually add any points to $\eps(\TZ)$.

\begin{lemma}
$\eps(\TZ)$ is a closed analytic subset of $T(F_0 \Mmod)$.
\end{lemma}

\begin{proof}
The restriction of $\eps(\TZ)$ to $X$ clearly has the same closure as $\eps(\TZ)$
itself. Since the closure is analytic, any of its points belongs to the image
of a holomorphic map $g \colon \Delta \to T(F_0 \Mmod)$, such that $g(\dst)$ is
contained in $\eps(\TZ) \cap p^{-1}(X)$. Lemma~\ref{lem:restriction-curves} shows
that $g(\Delta) \subseteq \eps(\TZ)$, and this proves that $\eps(\TZ)$ is
itself closed.
\end{proof}

\begin{corollary}
If Condition~\ref{eq:main-condition} is true whenever $\Xb - X$ is a divisor with
normal crossings and $\shH$ has unipotent local monodromy, then it is true in
general.
\end{corollary}

%%%%%%%%%%%%%%%%%%%%%%%%%%%
\section{The normal crossing case with unipotent monodromy}
\label{sec:normal-crossing}

\subsection{Introduction}

This part of the paper contains the proof of Condition~\ref{eq:main-condition} in the
case where $\Xb = \Delta^n$, $X = \dstn{n}$, and $\shH$ is a polarized variation of
Hodge structure of weight $-1$ on $X$ with unipotent monodromy. In this situation,
there is an explicit description of the filtered $\Dmod$-module $(\Mmod, F)$ in terms
of Deligne's canonical extension $\shHOe$; in particular,
\[
	F_0 \Mmod = \sum_{k \geq 0} F_k \Dmod \cdot F_{-k} \shHOe
\]
consists of all sections in $F_0 \shHOe$, all first-order derivatives of sections in
$F_{-1} \shHOe$, and so on. This means that we have a natural map 
\[
	T(F_0 \Mmod) \to T(F_0 \shHOe).
\]
It is known that the image of $\TZ$ in the vector bundle $T(F_0
\shHOe)$ is not well-behaved (the quotient is the so-called ``Zucker extension'').
But we shall see that $T(F_0 \Mmod)$, which only maps to a very restricted subset of
$T(F^0 \shHOe)$, solves this problem.

If we pull the variation back to the universal covering space $\HH^n$, it can be
viewed as family of Hodge structures $\Phit(z)$ on a fixed vector space $\HC$,
polarized by a fixed alternating form $Q$. Let $\sigma_1, \dotsc, \sigma_m$ be a
collection of sections that generate $F_0 \Mmod$ over $\Delta^n$. At each point $z
\in \HH^n$, we define
\[
	B(z,h) = \sup \menge{\abs{Q(h, \sigma_j(z))}}{j=1, \dotsc, m},
\]
noting that it gives a norm on $\HR$ because the Hodge structures in question have
weight $-1$. We can compare $B(z,h)$ to a fixed norm on $\HR$ by an inequality of the
form
\[
	\norm{h} \leq C(z) \cdot B(z,h),
\]
where $C(z)$ only depends on $z$. The main idea is to show that, even as the
imaginary parts of $z_1, \dotsc, z_n$ tend to infinity, the constant $C(z)$ remains bounded.

It is illustrative to compare this with the situation for the canonical extension.  Of
course, we could similarly define a quantity $B_0(z,h)$, using only sections of $F_0
\shHOe$, and compare it to $\norm{h}$ by a constant $C_0(z)$. It is then not hard to see
that $B_0(z,h)$ is equivalent to the Hodge norm in the Hodge structure $\Phit(z)$.
The norm estimates of \cite{CKS} and \cite{Kashiwara-HS} show that $C_0(z)$ need not
be bounded: it generally grows like a certain polynomial in the imaginary parts of
$z_1, \dotsc, z_n$. This points to a very interesting analogy, observed by H.~Clemens.
Namely, in the definition of $B(z,h)$, we are controlling not just the 
various holomorphic functions $Q(h, \sigma)$, for $\sigma \in F_0 \shHOe$, but also
some of their derivatives---indeed, the additional sections in $F_0 \Mmod$ arise
precisely as derivatives of sections of $F_0 \shHOe$. If we now think of the Hodge norm
as an $L^2$-norm, and of $B(z,h)$ as a kind of Sobolev norm, then the fact that the
uniform norm $\norm{h}$ is bounded by a fixed multiple of $B(z,h)$ resembles the
well-known Sobolev inequality.

\subsection{The normal form of a period map}
\label{subsec:normal-form}

In this section, we set up some basic notation, and describe how to represent the
period map $\Phit(z)$ in terms of the limit mixed Hodge structure coming from the
$\SL_2$-Orbit Theorem. All the results cited here can be found in
\cite{CKS}*{Section~4}.

We consider a variation of polarized Hodge structure of weight $-1$ on
$\dstn{n}$. Let $s = (s_1, \dotsc, s_n)$ be the standard holomorphic coordinates on
$\Delta^n$. Throughout, we shall make the assumption that the monodromy of the variation
around each divisor $s_j = 0$ is unipotent. As usual, let $\HH^n \to \dstn{n}$ be the
universal covering space, with $s_j = e^{2 \pi i z_j}$. Let $N_j$ be the logarithm of
the monodromy transformation around the divisor $s_j = 0$. 

The pullback of the variation to $\HH^n$ can be viewed as a varying Hodge filtration
$\Phit \colon \HH^n \to D$ on a fixed vector space $\HC$, where $D$ is a suitable
period domain.  Since the variation is integral and polarized, there is a fixed
underlying integral lattice $\HZ$, and a bilinear form $Q \colon \HZ \tensor \HZ \to
\ZZ$ that is alternating and nondegenerate.  As usual, we denote by $G_{\RR} =
\operatorname{Aut}(\HR, Q)$ the real Lie group determined by the pairing, and by
$\mathfrak{g}_{\RR}$ its Lie algebra.  By the Nilpotent Orbit Theorem
\cite{Schmid}*{Theorem~4.12}, we have
\[
	e^{-\sum z_j N_j} \Phit(z) = \Psi(s),
\]
with $\Psi \colon \Delta^n \to \Dcheck$ holomorphic.  Let $\Wn{n} = W(N_1, \dotsc,
W_n)$ be the monodromy weight filtration for the cone $\Cn{n} = C(N_1, \dotsc, N_n)$,
and set $W = \Wn{n} \decal{-1}$. Then $\bigl( W, \Psi(0) \bigr)$ is a
mixed Hodge structure, polarized by $Q$ and any element of $\Cn{n}$, in the sense of
\cite{CKS}*{Definition~2.26}. Let $\delta \in
L_{\RR}^{-1,-1} \bigl( W, \Psi(0) \bigr)$ be the unique real element for which
$\bigl( W, e^{-i \delta} \Psi(0) \bigr)$ is $\RR$-split
\cite{CKS}*{Proposition~2.20}, and define $F = e^{-i \delta} \Psi(0) \in \Dcheck$.
Note that $\delta$ commutes with every $N_j$.  Let 
\[
	I^{p,q} = I^{p,q}(W, F) = W_{p+q} \cap F^p \cap \overline{F^q}
\]
be Deligne's canonical decomposition of the $\RR$-split mixed Hodge structure $(W, F)$.

The Lie algebra $\glie$ inherits a decomposition
\[
	\glie = \bigoplus_{p,q} \glie^{p,q},
\]
with $\glie^{p,q}$ consisting of those $X$ that satisfy $X I^{a,b} \subseteq I^{a+p,b+q}$.
Then we have
\[
	\glie = \glie^F \oplus \qlie = \glie^F \oplus \bigoplus_{p < 0}
		\glie^{p,q},
\]
and $\qlie$ is a nilpotent Lie subalgebra of $\glie$.  This decomposition makes it
possible to write $e^{-i \delta} \Psi(s) = e^{\Gamma(s)} F$ for a unique holomorphic
map $\Gamma \colon \Delta^n \to \qlie$ with $\Gamma(0) = 0$. We can therefore put the
period map into the standard form
\begin{equation} \label{eq:period-std}
	\Phit(z) = e^{i \delta} e^{\sum z_j N_j} e^{\Gamma(s)} F = e^{X(z)} F,
\end{equation}
with $X(z) \in \qlie$ and hence nilpotent.

The horizontality of the period map implies the following relationship between
$\Gamma(s)$ and the nilpotent operators $N_j$; it can be found in
\cite{CK}*{Proposition~2.6}.

\begin{lemma} \label{lem:Gamma}
Let $\Phit(z) = e^{i \delta} e^{\sum z_j N_j} e^{\Gamma(s)} F$ be the normal form of a period map.
\begin{enumerate}
\item We have
\[
	d \Bigl( e^{\sum z_j N_j} e^{\Gamma(s)} \Bigr) =
		e^{\sum z_j N_j} e^{\Gamma(s)} 
			\biggl( d \, \Gamma_{-1}(s) + \sum_{j=1}^n N_j dz_j \biggr).
\]
\item	For every $j = 1, \dotsc, n$, the commutator
\[
	\biglie{N_j}{e^{\Gamma(s)}} = N_j e^{\Gamma(s)} - e^{\Gamma(s)} N_j
\]
vanishes along $s_j = 0$.
\end{enumerate}
\end{lemma}

\begin{proof}
Since $\Gamma(s)$, all the $N_j$, and $\delta$ belong to the nilpotent Lie algebra $\qlie$,
we can write
\[
	e^{i \delta} e^{\sum z_j N_j} e^{\Gamma(s)} = e^{X(z)}
\]
for a unique holomorphic $X \colon \HH^n \to \qlie$. From the definition of $\qlie$,
we have $X(z) = X_{-1}(z) + X_{-2}(z) + \dotsb$, with $X_{p}(z) \in \bigoplus_q
\glie^{p,q}$. Note that 
\[
	X_{-1}(z) = \Gamma_{-1}(s) + \sum_{j=1}^n z_j N_j + i \delta_{-1}.
\]
Horizontality of the period map $e^{X(z)} F$ is equivalent to the condition that
\[
	e^{-X(z)} \cdot d \bigl( e^{X(z)} \bigr) = d X_{-1}(z),
\]
which gives the first assertion (because $\delta$ is constant). Writing the condition
out explicitly, we get
\[
	e^{-\Gamma(s)} \sum_{j=1}^n N_j dz_j \cdot e^{\Gamma(s)} + 
		e^{-\Gamma(s)} \cdot d \bigl( e^{\Gamma(s)} \bigr) = 
	d \, \Gamma_{-1}(s) + \sum_{j=1}^n N_j dz_j.
\]	
Now $d s_j = 2 \pi i s_j \cdot dz_j$; thus if we evaluate the identity on the tangent
vector field $\vfelds{z_j}$, we get
\[
	N_j e^{\Gamma(s)} - e^{\Gamma(s)} N_j = 
		2 \pi i s_j \cdot \biggl( e^{\Gamma(s)} \vfeld{s_j} \Gamma_{-1}(s) 
			- \vfeld{s_j} e^{\Gamma(s)} \biggr).
\]
We then obtain the second assertion by setting $s_j = 0$.
\end{proof}

The fact that the commutator $\lie{N_j}{e^{\Gamma(s)}}$ vanishes along the divisor
$s_j = 0$ has the following highly useful consequence.

\begin{lemma} \label{lem:adN}
Let $y_j = \Im z_j$, and suppose that $y_1 \geq \dotsb \geq y_n \geq 1$ and $0 \leq
\Re z_j \leq 1$. Define the nilpotent operator $N = y_1 N_1 + \dotsb + y_n N_n$. Then
there is a constant $C > 0$ and an integer $m$, both independent of $z$, such that 
\[
	\bignorm{(\ad N)^k e^{\Gamma(s)}} \leq C \cdot \sum_{j=1}^n y_j^m e^{-2 \pi y_j}
\]
for any $k \geq 1$.
\end{lemma}

\begin{proof}
Since $\Gamma(s)$ is holomorphic in $s = (s_1, \dotsc, s_n)$, and $\Gamma(0) = 0$, we
can write
\begin{align*}
	e^{\Gamma(s)} - \id &= 
		\Bigl( e^{\Gamma(s_1, \dotsc, s_n)} - e^{\Gamma(0, s_2, \dotsc, s_n)} \Bigr)
		+ \dotsb + \Bigl( e^{\Gamma(0, \dotsc, 0, s_n)} - \id \Bigr) \\
		&= s_1 B_1(s) + \dotsc + s_n B_n(s),
\end{align*}
where each $B_j(s)$ is an operator that depends holomorphically on $s$. Moreover,
$B_j(s)$ commutes with $N_1, \dotsc, N_{j-1}$, and $\norm{B_j(s)}$ is uniformly
bounded, independent of $s$. We then compute that (for $k \geq 1$)
\[
	(\ad N)^k e^{\Gamma(s)} = (\ad N)^k \sum_{j=1}^n s_j B_j(s) 
		= \sum_{j=1}^n s_j \bigl( \ad (y_j N_j + \dotsb + y_n N_n) \bigr)^k B_j(s).
\]
Since each $N_j$ is nilpotent, $y_j \geq \dotsb \geq y_n$, and $\abs{s_j} = e^{-2 \pi y_j}$,
the assertion follows by taking norms.
\end{proof}

\subsection{Sections of the minimal extension}
\label{subsec:minimal-extension}

Now let $M$ be the intermediate extension of the variation of Hodge structure to a
polarized Hodge module on $\Delta^n$. In this section, we review M.~Saito's description
of the underlying filtered left $\Dmod$-module $(\Mmod, F)$.  Let $\shHO$ be the holomorphic vector
bundle on $\dstn{n}$ underlying the variation, and $\nabla$ the induced flat
connection on $\shHO$. Since the local monodromies are unipotent, $\shHO$ can be
canonically extended to a vector bundle $\shHOe$ on $\Delta^n$, such that the
connection has logarithmic poles along $s_1 \dotsb s_n = 0$ with nilpotent residues
\cite{Deligne}*{Proposition~5.2}.  More explicitly, for each $v \in \HC$, the map
\[
	\HH^n \to \HC, \qquad z = (z_1, \dotsc, z_n) \mapsto e^{\sum z_j N_j} v
\]
descends to a holomorphic section of $\shHO$ on $\dstn{n}$, and $\shHOe$ is the
locally free subsheaf of $\jl \shHO$ generated by all such sections. Using the standard
form of the period map in \eqref{eq:period-std}, the maps
\begin{equation} \label{eq:sections-canext}
	\HH^n \to \HC, \qquad z \mapsto e^{X(z)} v = e^{i \delta} e^{\sum z_j N_j} e^{\Gamma(s)} v
\end{equation}
also induce a collection of sections that generate $\shHOe$.

The Nilpotent Orbit Theorem implies that the Hodge bundles $F^p \shHO$ extend
uniquely to holomorphic subbundles $F^p \shHOe$ of the canonical extension. Each
$F^p \shHOe$ is generated by the sections in \eqref{eq:sections-canext} for $v \in
F^p$.

Now $\Mmod$, the minimal extension of $(\shHO, \nabla)$ to a holonomic $\Dmod$-module
on $\Delta^n$, is simply the $\Dmod$-submodule of $\jl \shHO$ generated by $\shHOe$.
Moreover, the Hodge filtration on $\Mmod$ is given by 
\[
	F_p \Mmod = \sum_{k \geq 0} F_k \Dmod \cdot F^{k-p} \shHOe.
\]
It satisfies $F_k \Dmod \cdot F_p \Mmod \subseteq F_{k+p} \Mmod$, and each $F_p
\Mmod$ is a coherent sheaf on $\Delta^n$ whose restriction to $\dstn{n}$ is $F^{-p}
\shHO$. For the purposes of our construction, the important point is that $F_p \Mmod$
has more sections than $F^{-p} \shHOe$; the following lemma exhibits the ones that we
will use.

\begin{lemma} \label{lem:sections-Fp}
For any subset $I \subseteq \{1, \dotsc, n\}$ of cardinality $m$, and for any
vector $v \in F^{-p}$, the formula
\[
	\sigma_{I,v}(z) = e^{X(z)} \prod_{j \in I} \frac{N_j}{s_j} \cdot v
\]
defines a holomorphic section of the coherent sheaf $F_{p+m} \Mmod$ on $\Delta^n$.
\end{lemma}

\begin{proof}
We work by induction on the cardinality $m$ of the set $I$. Then case $m = 0$ is
clear from the definition of $F_p \Mmod$. We may therefore suppose that the assertion
has been proved for all subsets of cardinality at most $m$, and consider $I \subseteq
\{1, \dotsc, n\}$ with $\abs{I} = m+1$. Let $k = \max I$ and $J = I - \{k\}$. Then
\[
	\sigma_{J,w}(z) = e^{X(z)} \prod_{j \in J} \frac{N_j}{s_j} \cdot w
\]
is a section of $F_{p+m} \Mmod$ for every $w \in F^{-p}$. 

Using that $k \not\in J$, the first identity in Lemma~\ref{lem:Gamma} shows that 
\begin{align*}
	\vfeld{s_k} \sigma_{J,v}(z) 
		&= e^{X(z)} \left( \frac{N_k}{2 \pi i s_k} + \pder{\Gamma_{-1}}{s_k} \right) \cdot
		\prod_{j \in J} \frac{N_j}{s_j} v \\
		&= \frac{\sigma_{I,v}(z)}{2 \pi i} + e^{X(z)}
			\pder{\Gamma_{-1}(s)}{s_k} \prod_{j \in J} \frac{N_j}{s_j} \cdot v.
\end{align*}
The second half of Lemma~\ref{lem:Gamma}, applied recursively, implies that
\[
	C(s) = \bigglie{\Gamma_{-1}(s)}{\prod_{j \in J} \frac{N_j}{s_j}} \qquad \text{and}
		\qquad \pder{C(s)}{s_k} = \bigglie{\pder{\Gamma_{-1}(s)}{s_k}}%
			{\prod_{j \in J} \frac{N_j}{s_j}}
\]
are holomorphic on $\Delta^n$; as operators, they map $F^{-p}$ to $F^{-p-m-1}$. In
\[
	\vfeld{s_k} \sigma_{J,v}(z) = \frac{\sigma_{I,v}(z)}{2 \pi i} + e^{X(z)} 
		\biggl( \pder{C(s)}{s_k} v + \prod_{j \in J} \frac{N_j}{s_j} 
			\pder{\Gamma_{-1}(s)}{s_k} v \biggr),
\]
the left-hand side defines a holomorphic section of $F_{p+m+1} \Mmod$; by induction, the same
is true for the second term on the right-hand side. We conclude that
$\sigma_{I,v}(z)$ is itself a section of $F_{p+m+1} \Mmod$, thus completing
the induction.
\end{proof}

\subsection{The main estimate}
\label{subsec:main-estimate}

Now fix a norm on the vector space $\HC$. Let $\sigma_1, \dotsc, \sigma_m$ be a
collection of sections that generate the coherent sheaf $F_0 \Mmod$ over $\Delta^n$.
To prove the closedness of $\eps(\TZ)$ inside of $T(F_0 \Mmod)$, our strategy is to
show that the norm of any vector $h \in \HR$
is bounded uniformly by the values of $Q \bigl( h, \sigma_j(z) \bigr)$, provided that
the imaginary parts of $z_1, \dotsc, z_n$ are sufficiently large.

As a matter of fact, we will prove a slightly stronger statement, involving only the
special sections $\sigma_{I,v}$ from Lemma~\ref{lem:sections-Fp}.  Given
a real vector $h \in \HR$, and a point $z \in \HH^n$, we thus introduce the quantity
\[
	B(z,h) = \sup \Menge{\bigabs{Q(h, \sigma_{I,v}(z))}}%
		{\text{$I \subseteq \{1, \dotsc, n\}$ and $v \in F^{\abs{I}}$ with $\norm{v} \leq 1$}},
\]
which gives a norm on $\HR$ for every $z \in \HH^n$. Since we are trying to control
the size of $h$ in terms of $B(z,h)$, we also let $N = y_1 N_1 + \dotsb + y_n N_n$,
and define
\[
	Z(y,h) = \max_{k \geq 0} \norm{N^k h},
\]
noting that $N$ is nilpotent. After a few preliminary
results on decompositions in $\RR$-split mixed Hodge structures in
\subsecref{subsec:decompositions}, the following estimate will be proved in
\subsecref{subsec:proof}.

\begin{theorem} \label{thm:main-estimate}
Let $\Phit(z) = e^{i \delta} e^{\sum z_j N_j} e^{\Gamma(s)} F$ be the normal form of a
variation of polarized Hodge structure of weight $-1$ on $\dstn{n}$. Fix a norm on
the underlying vector space $\HC$. Then there are constants $C > 0$ and $\alpha > 0$,
such that we have
\[
	Z(y,h) \leq C \cdot B(z,h)
\]
for every $h \in \HR$ and every $z \in \HH^n$ with $y_j = \Im z_j \geq \alpha$ and
$0 \leq \Re z_j \leq 1$.
\end{theorem}

\subsection{The closure of the set of integral points}

Granting Theorem~\ref{thm:main-estimate} for the time being, we shall now show that
Condition~\ref{eq:main-condition} is true: the map $\eps \colon \TZ \to T(F_0 \Mmod)$
is injective and has closed image. The first result is that any sequence of points in
$\TZ$ over $\dstn{n}$ that converges in $T(F_0 \Mmod)$ has to be eventually constant
and invariant under monodromy. This is the only point in the proof where we use the
fact that we are dealing with integral classes.

\begin{theorem} \label{thm:closure}
Let $z(m) \in \HH^n$ be a sequence of points with $\Im z_j(m) \to \infty$ and $\Re
z_j(m) \in \lbrack 0, 1 \rbrack$ for $j = 1, \dotsc, n$. Let $h(m) \in \HZ$ be a
corresponding sequence of integral classes, such that
\[
	Q \bigl( h(m), \sigma_{I,v}(z(m)) \bigr)
\]
is convergent for every $I \subseteq \{1, \dotsc, n\}$ and every $v \in F^{\abs{I}}$.
Then the sequence $h(m)$ is eventually constant, and its constant value satisfies
$N_k h(m) = 0$ for $k = 1, \dotsc, n$.
\end{theorem}

\begin{proof}
The first step is to show that $N_k h(m) = 0$ for all $k = 1, \dotsc, n$, and all
sufficiently large $m$. We begin by finding a subsequence of $h(m)$ along which this
is true. By assumption, the quantity $B \bigl( z(m), h(m) \bigr)$ is
bounded, and so the inequality in Theorem~\ref{thm:main-estimate} implies that
$\norm{h(m)}$ is bounded. Since $h(m) \in \HZ$, the sequence can take on only finitely many
distinct values; let $h \in \HZ$ be one of them.  The inequality also implies that
$\sum z_j(m) N_j h(m)$ is bounded; according to Lemma~\ref{lem:subsequence} below, 
$h \in \Wn{n}_0 = W_{-1}$, and we can find a subsequence along which
\[
	\sum_{j=1}^n z_j(m) N_j h(m) \to \sum_{j=1}^n w_j N_j h,
\]
for some $w \in \HH^n$ with $\Im w$ large. We then have $e^{-\sum z_j(m) N_j} h(m)
\to e^{-\sum w_j N_j} h$; by taking $I = \{k\}$ and $v \in F^1$ arbitrary, it follows
that 
\begin{equation} \label{eq:Nk=0}
\begin{split}
	0 = \lim_{m \to \infty} s_k(m) \cdot Q \bigl( h(m), \sigma_{\{k\},v}(z(m)) \bigr) 
		&= Q \bigl( h, e^{i \delta} e^{\sum w_j N_j} N_k v \bigr) \\
		&= -Q \bigl( N_k h, e^{i \delta} e^{\sum w_j N_j} v \bigr).
\end{split}
\end{equation}
Now $\bigl( W, e^{i \delta} e^{\sum w_j N_j} F \bigr)$ is a mixed Hodge structure; because
the vector $N_k h$ is rational and belongs to $W_{-3}$, we easily conclude from
\eqref{eq:Nk=0} that $N_k h = 0$. 

The argument above actually proves that $N_k h(m) = 0$ for all sufficiently large $m$
(otherwise, we could find a subsequence along which $N_k h(m) \neq 0$, leading to a
contradiction).  Consequently, $e^{-\sum z_j(m) N_j} h(m) = h(m)$, and so we find that
\[
	\lim_{m \to \infty} Q \bigl( h(m), \sigma_{\emptyset,v}(z(m)) \bigr) =
		\lim_{m \to \infty} Q \bigl( h(m), e^{i \delta} v \bigr)
\]
for every $v \in F^0$. By looking at the mixed Hodge structure $\bigl( W, e^{i \delta} F
\bigr)$, we deduce from the convergence of all those expressions that the sequence of
integral vectors $h(m) \in W_{-1}$ is itself convergent, and hence eventually constant.
\end{proof}

\begin{lemma} \label{lem:subsequence}
Let $h \in \HR$, and suppose that $z_1(m) N_1 h + \dotsb + z_n(m) N_n h$
remains bounded for $m \to \infty$. Then $h \in \Wn{n}_0$. Moreover, for any $\alpha
> 0$, there is a point $w \in \CC^n$ with $\max_{1 \leq j \leq n} \Im w_j \geq
\alpha$, such that
\[
	\sum_{j=1}^n {z_j(m) N_j h} \to \sum_{j=1}^n w_j N_j h
\]
is true along a subsequence.
\end{lemma}

\begin{proof}
We borrow a technique introduced by E.~Cattani, P.~Deligne, and A.~Kaplan \cite{CDK}*{p.~494}.
Let $x_j(m) = \Re z_j(m)$, and $y_j(m) = \Im z_j(m)$. After passing to a subsequence,
we can find constant vectors $\theta^1, \dotsc, \theta^r \in \RR^n$, whose components
satisfy the inequalities $0 \leq \theta_j^1 \leq \theta_j^2 \leq \dotsb \leq
\theta_j^r$, such that
\[
	y_j(m) = t_1(m) \theta^1 + \dotsb + t_r(m) \theta^r + \eta(m),
\]
where the ratios $t_1(m) / t_2(m)$, \dots, $t_{r-1}(m) / t_r(m)$, and $t_r(m)$
are tending to infinity, and the remainder term $\eta(m)$ is convergent. We can take
every $\Im \eta_j(m) \geq \alpha$; moreover, we may clearly assume that the bounded
sequence $x(m)$ is also convergent. Let $w(m) = x(m) + i \eta(m)$.  Along the
subsequence in question, we then have
\[
	\sum_{j=1}^n z_j(m) N_j h = \sum_{j=1}^n w_j(m) N_j h + 
		i \sum_{k=1}^r t_k(m) \sum_{j=1}^n \theta_j^k N_j h.
\]
This expression can only be bounded if $\sum \theta_j^k N_j h = 0$ for every
$k$; it follows that $h \in \Wn{n}_0$, because $\sum \theta_j^r N_j \in \Cn{n}$.
We now obtain the second assertion with $w = \lim_{m \to \infty} w(m)$.
\end{proof}

\begin{corollary} \label{cor:main-condition}
The map $\eps \colon \TZ \to T(F_0 \Mmod)$ is injective, and $\eps(\TZ)$ is a closed
analytic subset; therefore Condition~\ref{eq:main-condition} is true for polarized
variations of Hodge structure on $\dstn{n}$ with unipotent monodromy.
\end{corollary}

\begin{proof}
The map $\eps$ is injective because the induced map $\TZ \to T(F^0 \shHOe)$ is injective.
Its image is a closed analytic subset because of Theorem~\ref{thm:closure}.
\end{proof}

\subsection{Decompositions in $\RR$-split mixed Hodge structures}
\label{subsec:decompositions}

In this section, we collect several auxiliary results about the primitive decomposition in an
$\RR$-split polarized mixed Hodge structure, and its relationship with the Hodge
decomposition in the associated nilpotent orbit. This is a preparation for the proof
of Theorem~\ref{thm:main-estimate} in \subsecref{subsec:proof} below.

Throughout, we let $(W,F)$ be an $\RR$-split mixed Hodge structure, polarized by a
nondegenerate bilinear form $Q$ and a nilpotent operator $N$, such that $W_{\bullet}
= W(N)_{\bullet-m}$. Then $e^{iN} F$ is a point in the corresponding period domain $D$, and
therefore a polarized Hodge structure of weight $m$. Let 
\[
	\HC = \bigoplus_{p,q} I^{p,q}
\]
be Deligne's decomposition; since the mixed Hodge structure is split over
$\RR$, we have $I^{p,q} = W_{p+q} \cap F^p \cap \overline{F^q}$. The operator $Y$,
which acts as multiplication by $p+q-m$ on $I^{p,q}$, is then a real splitting of the
filtration $W(N)$; let $\Npl$ be the real operator making $(N, Y, \Npl)$ into an
$\sltwo$-triple.

There are two natural decompositions of the vector space $\HC$, and one purpose of
this section is to relate the two. The first one is Deligne's decomposition by the
$I^{p,q}$, the second one the primitive decomposition determined by the
nilpotent operator $\Npl$. The reason for using $\Npl$ instead of $N$ will become
apparent below. We define the primitive subspaces for the operator $\Npl$ as
\[
	I_0^{p,q} = I^{p,q} \cap \ker N.
\]
Given a vector $h \in \HC$, we denote by $h^{p,q}$ its component in the space
$I^{p,q}$, and then $h = \sum_{p,q} h^{p,q}$. We can also write $h$ uniquely in the
form
\[
	h = \sum_{p,q} \sum_{b=0}^{m-p-q} (\Npl)^b h^{p,q}(b)
\]
where each vector $h^{p,q}(b) \in I_0^{p,q}$ is primitive for $\Npl$, meaning that $N
h^{p,q}(b) = 0$.

\begin{lemma} \label{lem:primitive}
There are constants $C(p,q,b,j) \in \QQ$, depending only on the Hodge numbers of the
$\RR$-split mixed Hodge structure $(W,F)$, such that
\[
	h^{p,q}(b) = \sum_{j \geq 0} C(p,q,b,j) (\Npl)^j N^{b+j} h^{p+b,q+b}.
\]
\end{lemma}

\begin{proof}
Since $\Npl$ is a morphism of type $(1,1)$, a short computation shows that
\[
	N^a h^{p+a,q+a} = \sum_{j \geq 0} R(a,a+j,m-p-q+2j) (\Npl)^j h^{p-j,q-j}(a+j),
\]
where the constants are as in Lemma~\ref{lem:constants-R} below. Since $R(a,a,m-p-q)
\neq 0$ for $0 \leq a \leq m-p-q$, we can solve those equations for the $h^{p,q}(b)$ by
descending induction on $b$ to arrive at the stated formulas.
\end{proof}

\begin{lemma} \label{lem:constants-R}
Let $v \neq 0$ be a vector satisfying $N v = 0$ and $Y v = - \ell v$ (and therefore
$\ell \geq 0$).  Then $N^a (\Npl)^b v = R(a,b,\ell) (\Npl)^{b-a}v$, with
\[
	R(a,b,\ell) = \frac{b! (\ell+a-b)!}{(\ell-b)! (b-a)!}
\]
for $0 \leq a \leq b \leq \ell$, and $R(a,b,\ell) = 0$ in all other cases.
\end{lemma}

\begin{proof}
This is well-known; but since the proof is short, we include it here. We have
\[
	N^{a+1} (\Npl)^b v = N \cdot R(a,b,\ell) (\Npl)^{b-a} v = R(a,b,\ell) R(1,b-a,\ell)
		(\Npl)^{b-a-1} v,
\]
from which the identity $R(a+1,b,\ell) = R(a,b,\ell) R(1,b-a,\ell)$ follows. We also
have
\[
	N (\Npl)^{b+1} v = (\Npl N - Y) \cdot (\Npl)^b v = 
		\bigl( R(1,b,\ell) - (2b-\ell) \bigr) (\Npl)^b v,
\]
from which one sees that $R(1,b+1,\ell) = R(1,b,\ell) + (\ell-2b)$. Together with the evident
condition that $R(1,0,\ell) = 0$, the two equations suffice to prove the formula for
$R(a,b,\ell)$ by induction.
\end{proof}

The formula in Lemma~\ref{lem:primitive} shows how the size of the primitive
components depends on the two operators $\Npl$ and $N$. Since we will need this fact
in \subsecref{subsec:proof}, we state it as a corollary.
\begin{corollary} \label{cor:primitive-bound}
Fix a norm on the vector space $\HC$, and define $Z(N,h) = \max_{k \geq 0} \norm{N^k
h}$. Then there is a constant $C > 0$ and an integer $d \in \NN$, both
depending only on the Hodge numbers of $(W,F)$, such that
\[
	\max_{p,q,b} \norm{h^{p,q}(b)} \leq C \norm{\Npl}^d \cdot Z(N,h)
\]
for every $h \in \HC$.
\end{corollary}

We now specialize to the case $m = -1$. Then $e^{iN} F$ is a polarized Hodge structure of
weight $-1$ by \cite{CKS}*{Lemma~3.12}, and so we have the decomposition
\begin{equation} \label{eq:decomposition}
	\HC = e^{iN} F^0 \oplus e^{-iN} \overline{F^0}.
\end{equation}
Any vector $h \in \HC$ can therefore be written uniquely as $h = e^{iN} v + e^{-iN}
w$, with
\[
	v \in F^0 = \bigoplus_{p \geq 0} I^{p,q} \quad \text{and} \quad
	w \in \overline{F^0} = \bigoplus_{q \geq 0} I^{p,q}.
\]
The \emph{uniqueness} of the decomposition has a useful consequence that we shall now
explain. Let $w = \sum (\Npl)^b w^{p,q}(b)$ be the primitive decomposition of the
vector $w \in \overline{F^0}$; note that 
\[
	w^{p,q} = \sum_{b \geq 0} (\Npl)^b w^{p-b,q-b}(b),
\]
which implies that $w^{p,q}(b) = 0$ unless $q+b \geq 0$. Set $g = e^{-iN} h$, and
similarly write $g = \sum (\Npl)^b g^{p,q}(b)$.
The decomposition in \eqref{eq:decomposition} becomes $g = v + e^{-2iN} w$, and since
$v \in F^0$, the vector $w$ is uniquely defined by the condition that
\[
	g^{p,q} = \bigl( e^{-2iN} w \bigr)^{p,q}
\]
for every $p \leq -1$ and every $q$. The right-hand side can be expanded as
\begin{align*}
	\bigl( e^{-2iN} w \bigr)^{p,q} &= \sum_{k,a \geq 0} \frac{(-2i)^k}{k!}
		N^k (\Npl)^a w^{p-a+k,q-a+k}(a) \\
		&= \sum_{k,b \geq 0} \frac{(-2i)^k}{k!} R(k,k+b,2b-1-p-q) 
			(\Npl)^b w^{p-b,q-b}(k+b).
\end{align*}
By equating primitive components, we obtain the set of equations
\begin{equation} \label{eq:equations-w}
	g^{p,q}(b) = \sum_{k \geq 0} \frac{(-2i)^k}{k!} R(k,k+b,-1-p-q)
			w^{p,q}(k+b)
\end{equation}
for $p+b \leq -1$. We point out one more time that $w^{p,q}(b) = 0$ unless $q+b \geq 0$.

We now consider \eqref{eq:equations-w} as a system of linear equations for the
vectors $w^{p,q}(b) \in I_0^{p,q}$ with $q+b \geq 0$. Since the decomposition $g = v + e^{-2iN} w$
is unique, the system must have a unique solution, which means that its coefficient
matrix has to be invertible. It follows that there are constants $\Gamma(p,q,b,a) \in
\QQ(i)$ such that
\[
	w^{p,q}(b) = \sum_{a=0}^{-1-p} \Gamma(p,q,b,a) g^{p,q}(a);
\]
the upper limit for the summation stems from the condition $p+a \leq -1$.
Since the proof of Theorem~\ref{thm:main-estimate} in
\subsecref{subsec:proof} is entirely based on the solution to the system of
equations in \eqref{eq:equations-w}, we summarize the result in the following proposition.

\begin{proposition} \label{prop:system}
Consider the system of equations (for $p+b \leq -1$)
\[
	y^{p,q}(b) = \sum_{k \geq 0} \frac{(-2i)^k}{k!} R(k,k+b,-1-p-q) \cdot x^{p,q}(k+b)
\]
in the unknowns $\{ x^{p,q}(b) \}_{q+b \geq 0}$.
Given any collection of vectors $\{ y^{p,q}(b) \}_{p+b \leq -1}$, the unique solution
to the system is given by the formula
\[
	x^{p,q}(b) = \sum_{a=0}^{-1-p} \Gamma(p,q,b,a) \cdot y^{p,q}(a),
\]
where $\Gamma(p,q,a,b) \in \QQ(i)$ are certain constants that depend only on the Hodge
numbers $\dim I^{p,q}$ of the $\RR$-split mixed Hodge structure $(W,F)$.
\end{proposition}

\subsection{Proof of the main estimate}
\label{subsec:proof}

After the preliminary work in the previous section, we now come to the proof of the
estimate from Theorem~\ref{thm:main-estimate}. Given a point $y \in \HH^n$, we set $N
= y_1 N_1 + \dotsb + y_n N_n$; note that the weight filtration $W(N)$ is independent
of $y$. Together with the bilinear form $Q$, the nilpotent operator $N$ polarizes the
$\RR$-split mixed Hodge structure $(W, F)$, where $W = W(N) \decal{-1}$.
Let $Y$ be the real splitting of $W(N)$ determined by Deligne's decomposition $\HC =
\bigoplus I^{p,q}$, and let $(N, Y, \Npl)$ be the corresponding $\sltwo$-triple. An
important observation is that the operator $\Npl$ is of order $1/y_n$; this is a
simple consequence of the $\SL_2$-Orbit Theorem of \cite{CKS}.

\begin{lemma} \label{lem:Npl}
There are constants $C > 0$ and $\alpha > 0$ such that $\norm{\Npl} \leq C/y_n$ holds for all $y_1
\geq \dotsb \geq y_n \geq \alpha$.
\end{lemma}

\begin{proof}
Since $y_n \Npl = (N/y_n)^{+}$, it follows from \cite{CKS}*{Theorem~4.20} that the
operator $y_n \Npl$ has a power series expansion in nonpositive powers of $y_1/y_2$,
\dots, $y_{n-1}/y_n$, convergent in a region of the form $y_2/y_1 > \beta$, \dots,
$y_n/y_{n-1} > \beta$ for some $\beta > 0$. The assertion follows from this via
dependence on parameters. More precisely, we argue as follows.

Suppose to the contrary that $y_n \Npl$ was not bounded. Since $y_n \Npl =
(N/y_n)^{+}$ depends only on the ratios $y_1/y_2$, \dots, $y_{n-1}/y_n$, we can
then find a sequence of points $y(m)$ with $y_1(m) \geq \dotsb \geq y_n(m)$ and
$y_n(m) \to \infty$, along which $\norm{y_n \Npl}$ diverges.
After passage to a subsequence, we can arrange that
\[
	y_1(m) N_1 + \dotsb + y_n(m) N_n = t_1(m) M_1(m) + \dotsb + t_r(m) M_r(m)
\]
where $t_1(m)/t_2(m)$, \dots, $t_{r-1}(m)/t_r(m)$, and $t_r(m) = y_n(m)$ are going to
infinity, and each $M_j(m)$ is a linear combination of $N_1, \dotsc, N_n$ with
coefficients that lie in a bounded interval $\lbrack 1, K \rbrack$. 
By \cite{CDK}*{Remark~4.65}, the data in the $\SL_2$-Orbit Theorem depend real
analytically on these coefficients; we can therefore use the convergence of the
series as above to conclude that 
\[
	y_n \Npl	= \biggl( \frac{t_1(m)}{t_r(m)} M_1(m) + \dotsb +
		\frac{t_{r-1}(m)}{t_r(m)} M_{r-1}(m) + M_r(m) \biggr)^{+}
\]
remains bounded as $m \to \infty$. But this clearly contradicts our original
assumption, and so the lemma is proved.
\end{proof}

We now use the boundedness of $y_n \Npl$, together with the analysis in
\subsecref{subsec:decompositions}, to prove the following important estimate. The
decomposition is based on the fact that $e^{i \delta} e^{iN} e^{\Gamma(s)} F =
e^{-\sum x_j N_j} \Phit(z)$ defines a Hodge structure of weight $-1$ on $\HC$.

\begin{proposition} \label{prop:main-estimate}
Let $h \in \HC$ be any vector, and define $u \in F^0$ and $v \in \overline{F^0}$
through the decomposition $h = e^{i \delta} e^{iN} e^{\Gamma(s)} u + e^{-i \delta}
e^{-iN} e^{\overline{\Gamma(s)}} v$.  Then there are constants $\alpha \geq 1$ and $C
> 0$, such that
\begin{equation} \label{eq:main-estimate}
	Z(y,v) \leq C \cdot Z(y,h),
\end{equation}
provided that $y_1 \geq \dotsb \geq y_n \geq \alpha$.
\end{proposition}

\begin{proof}
We let $g = e^{-\Gamma(s)} e^{-iN} e^{-i \delta} h$, and observe that $Z(y,g)$ is
bounded by a constant multiple of $Z(y,h)$. Let
\[
	g = \sum (\Npl)^b g^{p,q}(b)
\]
be the primitive decomposition of $g$ determined by $\Npl$, with $g^{p,q}(b) \in
I_0^{p,q}$ in the notation of \subsecref{subsec:decompositions}. 
According to Corollary~\ref{cor:primitive-bound}, the quantity $\max
\norm{g^{p,q}(b)}$ is still bounded by a fixed multiple of $Z(y,h)$.

Similarly write the primitive decomposition of the vector $v$ as
\[
	v = \sum (\Npl)^b v^{p,q}(b),
\]
keeping in mind that $v \in \overline{F^0}$ means that $v^{p,q}(b) = 0$ unless $q+b
\geq 0$. We will prove the estimate in \eqref{eq:main-estimate} by showing that 
$\max_{p,q,b} \norm{v^{p,q}(b)}$ is bounded by a constant multiple of $\max_{p,q,b}
\norm{g^{p,q}(b)}$, and hence by $Z(y,h)$; this clearly suffices because
$\norm{\Npl}$ is bounded due to Lemma~\ref{lem:Npl}.  

The vector $v$ in the decomposition is uniquely determined by the condition that
\[
	g - e^{-\Gamma(s)} e^{-2iN} e^{-2i \delta} e^{\overline{\Gamma(s)}} v \in F^0.
\]
If we set $w = e^{2iN} e^{-\Gamma(s)} e^{-2iN} \cdot e^{-2i \delta}
e^{\overline{\Gamma(s)}} v$, then we can use Deligne's decomposition $\HC = \bigoplus
I^{p,q}$ to recast that condition into the form 
\[
	g^{p,q} = \bigl( e^{-2iN} w \bigr)^{p,q} \quad 
		\text{for any $p \leq -1$ and any $q$.}
\]
We will show that this system of equations is a perturbation (of order $1/y_n$) of a
triangular system. The following convention greatly simplifies the book-keeping:

\begin{notation}
For two vectors $h_1, h_2 \in \HC$, we shall write $h_1 \equiv h_2$ to mean that
\[
	h_1 - h_2 = \sum_{p,q,b} P(p,q,b) v^{p,q}(b)
\]
for linear operators $P(p,q,b)$ that are allowed to depend on $z$ (but not on $v$),
and have to satisfy $\max \norm{P(p,q,b)} \leq B/y_n$ for a constant $B$ that
is independent of $z$. It is easy to see that if $X$ is a linear operator such that
$\norm{X}$ is bounded independently of $z$, then $h_1 \equiv h_2$ implies $X h_1
\equiv X h_2$.
\end{notation}

We begin our analysis by observing that the operator $\delta$ is nilpotent,
since it belongs to $L_{\RR}^{-1,-1}(W,F)$. Let $\Delta = e^{-2i \delta}$; then we
have
\[
	e^{-2i \delta} = \id + \sum_{p,q \geq 1} \Delta_{-p,-q},
\]
where $\Delta_{-p,-q}$ maps $I^{a,b}$ into $I^{a-p,b-q}$.

Next, we look more carefully at the relationship between $w$ and $v$. To begin with,
the boundedness of $y_n \Npl$, proved in Lemma~\ref{lem:Npl}, implies that
\begin{align*}
	N^b v 
		&= \sum_{p,q,a} N^b (\Npl)^a v^{p,q}(a) \\
		&= \sum_{p,q,a} R(b,a,-1-p-q) (\Npl)^{a-b} v^{p,q}(a) \\
	&\equiv \sum_{p,q} R(b,b,-1-p-q) v^{p,q}(b).
\end{align*}
According to the formula in Lemma~\ref{lem:primitive},
\[
	w^{p,q}(b) = \sum_{j \geq 0} C(p,q,b,j) (\Npl)^j N^{b+j} w^{p+b,q+b};
\]
to connect this with the primitive decomposition for the vector $v$, we compute
\[
	N^{b+j} w = N^{b+j} e^{2iN} e^{-\Gamma(s)} e^{-2iN} e^{-2i \delta}
				e^{\overline{\Gamma(s)}} v \\
		\equiv e^{-2i \delta} N^{b+j} v,
\]
using Lemma~\ref{lem:adN} to neglect the terms that arise when commuting $N^{b+j}$
past the two operators $e^{-\Gamma(s)}$ and $e^{\overline{\Gamma(s)}}$. Consequently,
\[
	N^{b+j} w \equiv e^{-2i \delta} \sum_{p,q} R(b+j,b+j,-1-p-q) v^{p,q}(b+j).
\]	
Again using the boundedness of $y_n \Npl$, this shows that we are allowed to write
\[
	w^{p,q}(b) \equiv C(p,q,b,0) N^b w^{p+b,q+b} 
		= C(p,q,b,0) \bigl( N^b w \bigr)^{p,q}.
\]
Combining the various pieces of information, and remembering that $C(p,q,b,0) \cdot
R(b,b,-1-p-q) = 1$, we find that there are constants $D(p,q,b,j,k) \in \QQ$ with the
property that
\begin{equation} \label{eq:relation-w-v}
	w^{p,q}(b) \equiv
		v^{p,q}(b) + \sum_{j,k \geq 1} D(p,q,b,j,k) \cdot \Delta_{-j,-k} v^{p+j,q+k}(b).
\end{equation}

Since we have $g^{p,q} = (e^{-2iN} w)^{p,q}$ for $p \leq -1$, the primitive
components of $g$ and $w$ are
related by the equations in \eqref{eq:equations-w}. Using the constants
$\Gamma(p,q,b,a) \in \QQ(i)$ introduced in Proposition~\ref{prop:system}, we define
\[
	G^{p,q}(b) = \sum_{a=0}^{-1-p} \Gamma(p,q,b,a) g^{p,q}(a).
\]
It follows that we can express each $w^{p,q}(b)$ with $q+b \geq 0$ as a linear
combination of $G^{p,q}(b)$ and the vectors $\{w^{p,q}(a)\}_{q+a<0}$. For $q+b \geq
0$, we therefore have 
\[
	w^{p,q}(b) = G^{p,q}(b) + \sum_{a < -q} E(p,q,b,a) w^{p,q}(a)
\]
with certain constants $E(p,q,b,a) \in \QQ(i)$ that again depend on nothing but the
Hodge numbers of $(W,F)$. Now we observe that for $q+a < 0$, the relation in
\eqref{eq:relation-w-v} simplifies to
\[
	w^{p,q}(a) \equiv
		\sum_{j,k \geq 1} D(p,q,a,j,k) \cdot \Delta_{-j,-k} v^{p+j,q+k}(a),
\]
due to the fact that $v^{p,q}(a) = 0$. When we combine the two formulas for
$w^{p,q}(b)$ from above, we obtain for $q+b \geq 0$ an equation of the form
\[
	v^{p,q}(b) \equiv G^{p,q}(b)  
		+ \sum_{j,k \geq 1} \sum_{a \geq -q-k} D(p,q,a,j,k) \Delta_{-j,-k}
			v^{p+j,q+k}(a).
\]
Recalling the definition of the symbol $\equiv$, this means that
there are linear operators $P_{j,k}(b,c)$, mapping $I^{p,q}$ to $I^{p+j,q+k}$, and of
size $\norm{P_{j,k}(b,c)} \leq B/y_n$ for a suitable constant $B > 0$, such that
\begin{equation} \label{eq:perturbed}
\begin{split}
	G^{p,q}(b) = v^{p,q}(b) &-
		\sum_{j,k \geq 1} \sum_{a \geq -q-k} D(p,q,a,j,k) \Delta_{-j,-k} v^{p+j,q+k}(a) \\
		&+ \sum_{j,k} \sum_{c \geq -q+k} P_{j,k}(b,c) v^{p-j,q-k}(c).
\end{split}
\end{equation}
Once again, we view this as a system of linear equations relating the primitive
components $\{v^{p,q}(b)\}_{q+b \geq 0}$ to the vectors $\{G^{p,q}(b)\}_{q+b \geq 0}$.

Here comes the crucial point: Consider the system of equations (for $q+b \geq 0$)
\[
	G^{p,q}(b) = v^{p,q}(b) - \sum_{j,k \geq 1} \sum_{a \geq -q-k} 
		D(p,q,a,j,k) \Delta_{-j,-k} v^{p+j,q+k}(a)
\]
in the vectors $\{v^{p,q}(b)\}_{q+b \geq 0}$. It is evidently triangular; written in
matrix form, the matrix of coefficients has determinant equal to $1$. Since
$\norm{P_{j,k}(b,c)} \leq B/y_n$, we can now choose $\alpha \geq 1$ sufficiently
large to guarantee that the coefficient matrix of the system in \eqref{eq:perturbed} has
determinant close to $1$ for $y_n \geq \alpha$. The system can then be solved for the
$v^{p,q}(b)$, in such a way that $\max_{p,q,b} \norm{v^{p,q}(b)}$ is bounded by a
constant multiple of $\max_{p,q,b} \norm{G^{p,q}(b)}$.  It follows that there is a
large constant $K > 0$ (depending on the Hodge numbers of $(W,F)$ and on $B$) such that
\[
	\sum_{p,q,b} \norm{v^{p,q}(b)} \leq K \cdot Z(y,h).
\]

The decomposition $v = \sum (\Npl)^b v^{p,q}(b)$ implies that each $N^k v$ can again
be written as a combination of vectors of the form $(\Npl)^b v^{p,q}(b+k)$. Since
Lemma~\ref{lem:Npl} bounds the size of $\Npl$, it is then easy to see that we have
$Z(y,v) \leq C \cdot Z(y,h)$ for a suitable constant $C > 0$, as long as $y_n \geq
\alpha$.
\end{proof}

\begin{note}
If we look more carefully at the calculation above, we find that each $P_{j,k}(b,c)$
is one of the Hodge components of an operator that is built up from $\delta$, $N$,
$\Npl$, $\Gamma(s)$, and $\overline{\Gamma(s)}$. What the proof actually shows 
is that $v$ can be expressed by a very complicated formula in the Hodge components of
those operators and the $h^{p,q}$. Similar reasoning can be used to prove that
the entire Hodge decomposition of $h$ in the Hodge structure $e^{i \delta} e^{\sum
z_j N_j} e^{\Gamma(s)} F$ is given by formulas of this type.
\end{note}

Having completed the main technical step, we can now prove Theorem~\ref{thm:main-estimate}.

\begin{proof}
Fix a real vector $h \in \HR$, and let $z \in \HH^n$ be any point with
$x_j = \Re z_j \in \lbrack 0, 1 \rbrack$. Without loss of generality, we may assume
that $y_1 \geq \dotsb \geq y_n \geq \alpha$, where $y_j = \Im z_j$. We will specify
shortly how large $\alpha$ needs to be to obtain the asserted inequality
between $Z(y,h)$ and $B(z,h)$. By definition, the various pairings
\[
	Q \biggl( h, e^{i \delta} e^{\sum z_j N_j} e^{\Gamma(s)} 
		\prod_{j \in I} \frac{N_j}{s_j} v \biggr)
\]
are bounded by $B(z,h)$ for $v \in F^{\abs{I}}$ with $\norm{v} \leq 1$. Since $0 \leq
x_j \leq 1$ for each $j$, we may replace $h$ by $e^{-\sum x_j N_j} h$ without affecting
the statement we are trying to prove. For the same choices of $I$ and $v$ as above,
we then have
\[
	\biggl\lvert Q \biggl( h, e^{i \delta} e^{iN} e^{\Gamma(s)}
		\prod_{j \in I} \frac{N_j}{s_j} v \biggr) \biggr\rvert \leq B(z,h).
\]
Let us introduce the auxiliary vector $w = e^{-\Gamma(s)} e^{-iN} e^{-i \delta} h$.
Since $N = \sum y_j N_j$ and $\abs{s_j} = e^{-2\pi y_j}$, it is easy to
deduce that
\[
	Q(N^k w, v) = (-1)^k Q \bigl( h, e^{i \delta} e^{iN} e^{\Gamma(s)} N^k v \bigr)
\]
is bounded by a constant times $B(z,h)$, for any $v \in F^k$ with $\norm{v} \leq 1$.
The fact that the pairing is nondegenerate and compatible with the decomposition $\HC
= \bigoplus I^{p,q}$ now implies that the norm of each vector $N^k w^{p,q}$ with $p
\leq -1$ is bounded by a constant multiple of $B(z,h)$. To exploit
this information, we define
\[
	h' = e^{i \delta} e^{i N} e^{\Gamma(s)} \sum_{p \leq -1} w^{p,q} 
		= e^{i \delta} \cdot e^{iN} e^{\Gamma(s)} e^{-iN} \cdot  
		\sum_{p \leq -1} e^{iN} w^{p,q},
\]
and observe that, as a consequence of Lemma~\ref{lem:adN}, $Z(y,h') \leq C_1 \cdot
B(z,h)$ for some constant $C_1 > 0$.

By construction, $h = h' + r$, where $r$ belongs to $e^{i \delta} e^{iN}
e^{\Gamma(s)} F^0$. Because $h$ is real, it follows that
$\overline{h'} - h' = r - \overline{r}$. This is a partial Hodge decomposition for
the vector $\overline{h'} - h' \in \HC$, relative to the Hodge structure of weight $-1$
defined by the point $e^{i \delta} e^{iN} e^{\Gamma(s)} F = e^{-\sum x_j N_j} \Phit(z) \in
D$. Proposition~\ref{prop:main-estimate}, applied to $\overline{h'} - h'$, shows that we have
$Z(y,r) \leq C_2 \cdot Z(y,h')$ for another constant $C_2 > 0$. The asserted bound on
$Z(y,h)$ is now a consequence of the identity $h = h' + r$ and the inequality
$Z(y,h') \leq C_1 \cdot B(z,h)$.
\end{proof}

\subsection{Graphs of admissible normal functions}
\label{subsec:graphs-NC}

With very little additional effort, the method in \subsecref{subsec:proof} extends to
the study of normal functions with possibly nontrivial singularities.  Continuing
with the notation from \subsecref{subsec:normal-form}, let $\shH$ be a polarized
variation of Hodge structure of weight $-1$ on $\dstn{n}$, and let $\nu$ be a normal
function, admissible relative to $\Delta^n$. We represent $\nu$ by an
admissible variation of mixed Hodge structure $\shV$, in the form of an extension
\begin{diagram}
0 &\rTo& \shH &\rTo& \shV &\rTo& \ZZ(0) &\rTo& 0.
\end{diagram}
Since $\shHZ$ has unipotent monodromy, the same is clearly true for $\shVZ$. Let
$\shVd$ denote the dual variation of Hodge structure. As in
\subsecref{subsec:Jacobians}, we have an isomorphism $F_{-1} \shVOd \simeq F_0
\shHO$, using the fact that $\shH$ is polarized. The extension therefore gives rise
to a map of sheaves $\shVZ \to (F_0 \shHO)^{\vee}$ on $\dstn{n}$. Let $\Tnu$ be the
subset of the \'etal\'e space of $\shVZ$, consisting of those points that map to $1
\in \ZZ$. We then have a holomorphic embedding 
\[
	\varphi \colon \Tnu \into T(F_0 \shHO)
\]
over $\dstn{n}$, and the goal of this section is to prove that the closure of
$\varphi(\Tnu)$ inside the bigger space $T(F_0 \Mmod)$ is an analytic subset.

Just as in the pure case, we let $\VC$ denote the fiber of the pullback of $\shV$ to
$\HH^n$.  Let $W$ be the resulting weight filtration on $\VZ$, with $W_{-1} = \HZ$
and $\Gr_0^W \simeq \ZZ$. Let $N_1', \dotsc, N_n' \in \End{\VQ}$ be the logarithms of
the monodromy operators; note that $\im N_j' \subseteq \HQ$, and that the restriction
of $N_j'$ to $\HQ$ equals $N_j$.

\begin{notation}
It will be convenient to let $\VCone \subseteq \VC$ denote the subset of elements
that map to $1 \in \CC \simeq \VC / \HC$. We similarly define $\VRone$ and $\VZone$.
\end{notation}

The lifting of the period map will be denoted by $\Phit' \colon \HH^n \to D'$; since
the original variation is admissible, we have $e^{-\sum z_j N_j'} \Phit'(z) =
\Psi'(s)$ with $\Psi'$ holomorphic on $\Delta^n$. In addition, the relative monodromy
weight filtration $M = M(N_1', \dotsc, N_n'; W)$ exists and is constant on the open
cone $C(N_1', \dotsc, N_n')$, and the pair $\bigl( M, \Psi'(0) \bigr)$ is a
mixed Hodge structure \cite{Kashiwara-study}*{Proposition~5.2.1}. Let $\delta' \in
L_{\RR}^{-1,-1} \bigl( M, \Psi'(0) \bigr)$ be the unique element for which $(M, F)$
is $\RR$-split, where $F = e^{-i \delta} \Psi'(0)$. As in
\subsecref{subsec:normal-form}, we can now put the period map for the variation of
mixed Hodge structure into the standard form \cite{Pearlstein}*{Proof of Theorem~6.13}
\[
	\Phit'(z) = e^{i \delta'} e^{\sum z_j N_j'} e^{\Gamma'(s)} F,
\]
where $\Gamma'$ is holomorphic and satisfies $\Gamma'(0) = 0$. Since the period map
is again horizontal, Lemma~\ref{lem:Gamma} extends to this setting. Evidently, the
restriction of $\Gamma'$ to $\HC$ equals $\Gamma$, that of $\delta'$ equals $\delta$,
and so on. 

In the remainder of this section, we prove the following generalization of
Theorem~\ref{thm:closure}; note the similarity with the main result of E.~Cattani,
P.~Deligne, and A.~Kaplan \cite{CDK}*{Theorem~2.16}. As for notation, we
let $Q'$ denote the pairing between $\VC$ and sections (of the pullback to $\HH^n$)
of $F_0 \shHO$, induced by the map $\shVC \to (F_0 \shHO)^{\vee}$ described above.

\begin{theorem} \label{thm:closure-VRone}
Let $z(m) \in \HH^n$ be a sequence of points with $\Im z_j(m) \to \infty$ and $\Re
z_j(m) \in \lbrack 0, 1 \rbrack$ for $j=1, \dotsc, n$. Let $v(m) \in \VZone$ be a
corresponding sequence of integral classes, such that $Q' \bigl( v(m), \sigma_{I,u}(z(m))
\bigr)$ converges for every $I \subseteq \{1, \dotsc, n\}$ and every $u \in
F^{\abs{I}} \cap \HC$ (see \subsecref{subsec:minimal-extension}). 
Then the following three things are true:
\begin{enumerate}[label=(\roman{*}), ref=(\roman{*})]
\item The sequence $v(m)$ is bounded, hence takes only finitely many values. 
\label{en:VRone-i}
\item Let $v \in \VZone$ be a point of accumulation. Then there are positive integers
$a_1, \dotsc, a_n$ with the property that $a_1 N_1' v + \dotsb + a_n N_n' v = 0$.
\label{en:VRone-ii}
\item There is a vector $w \in \CC^n$ such that
\[
	e^{-\Gamma'(s(m))} e^{-\sum z_j(m) N_j'} e^{-i \delta'} v(m) \to
		e^{-\sum w_j N_j'} e^{-i \delta'} v
\]
along a subsequence of the original sequence.
\label{en:VRone-iii}
\item For each $k=1, \dotsc, n$, we have
\[
	e^{-\sum w_j N_j} e^{-i \delta} N_k' v = e^{-\sum \Re w_j N_j} N_k' v^{0,0} = N_k' v^{0,0},
\]
which implies that the vector $N_k' v$ is a rational Hodge class of type $(-1,-1)$ in
the mixed Hodge structure $\bigl( M \cap H, e^{\sum w_j N_j} \Psi(0) \bigr)$.
\label{en:VRone-iv}
\end{enumerate}
\end{theorem}

The proof proceeds through a sequence of lemmas. In analogy with the notation used in
\subsecref{subsec:main-estimate}, we define $N' = y_1 N_1' + \dotsb +
y_n N_n'$ and $N = y_1 N_1 + \dotsb + y_n N_n$, and observe that $M = M(N',W)$ is the
relative weight filtration for $N'$.  Consequently, we have $M_{-1} \subseteq W_{-1}$
and $M_0 + W_{-1} = W_0$, and $M \cap H = W(N) \decal{-1}$ is the shifted monodromy
weight filtration for $N$ on $H$. 

\begin{lemma}
There is a unique element $v_0 \in M_0 \cap F^0 \cap \VRone$ with $N' v_0 = 0$.
\end{lemma}

\begin{proof}
Given that $F^0 \cap \ker N' \cap \HR = \{0\}$, the uniqueness of such an element is clear; it
remains to show its existence. Since $M_0 + \HQ = \VQ$, we can certainly find an element
$v \in M_0 \cap \VQ$ that lifts $1 \in \QQ$. Since $(M,F)$ is $\RR$-split, we can
replace $v$ by its component in the space $I^{0,0}(M,F)$ and assume that $v$ is real
and lies in $I^{0,0}(M,F)$. Then $N'v$ belongs to $I^{-1,-1}(M \cap H, F \cap H)$ and
hence to $W(N)_{-1}$, and so there is an element $h \in \HR$ with $N'v = Nh$. Again
replacing $h$ by one of its components, we may assume that $h \in I^{0,0}(M \cap H, F
\cap H)$. But now $v_0 = v - h$ satisfies all the required conditions. 
\end{proof}

Fix a norm on the vector space $\VC$. As in \subsecref{subsec:proof}, the analysis in
this section depends mostly on a single difficult statement, namely that $\norm{v_0}$
remains bounded as $y_1, \dotsc, y_n \to \infty$. This is a special case of a more
general theorem due to P.~Brosnan and G.~Pearlstein \cite{BP3}, and as in their work,
relies on the $\SL_2$-Orbit Theorem of K.~Kato, C.~Nakayama,
and S.~Usui \cite{KNU-SL}. Observe that the pair $(W, e^{iN'} F)$ defines an
$\RR$-split mixed Hodge structure, due to the fact that $(M,F)$ splits over $\RR$.
Since $N' v_0 = 0$, it is obvious that $v_0$ is the unique real element in $I^{0,0}
\bigl( W, e^{iN'} F \bigr)$ that maps to $1 \in \Gr_0^W$; said differently, $v_0$ is
the image of $1$ under the canonical splitting of $(W, e^{iN'} F)$
\cite{KNU-SL}*{Section~1.2}.

\begin{lemma} \label{lem:v0}
There are constants $C > 0$ and $\alpha > 0$, such that $\norm{v_0} \leq C$ for all
$y_1, \dotsc, y_n \geq \alpha$.
\end{lemma}

\begin{proof}
Without loss of generality, we may suppose that $y_1 \geq \dotsb \geq y_n \geq
\alpha$. \cite{KNU-SL}*{Theorem~0.5} implies that the canonical splitting of $(W,
e^{iN'} F)$ has a power series expansion in nonpositive powers of $y_1/y_2$, \dots,
$y_{n-1}/y_n$, and $y_n$; the series converges provided that $y_1/y_2 > \beta$,
\dots, $y_{n-1}/y_n > \beta$, and $y_n > \beta$. Arguing as in the proof of
Lemma~\ref{lem:Npl}, we conclude that the canonical splitting is uniformly bounded
for all $y_1, \dotsc, y_n \geq \alpha$, once we take $\alpha$ sufficiently large. The
same is therefore true for the image of $1 \in \Gr_0^W$ under the canonical
splitting; but this image is precisely $v_0$.
\end{proof}

For $v \in \VR$, define $Z(y,v) = \max_{k \geq 0} \norm{(N')^k v}$. As before, we
have to show that the norm $\norm{v}$ of a real vector $v \in \VRone$ is controlled
by the size of the pairings $Q' \bigl( v, \sigma_{I,u}(z) \bigr)$, once $y_1, \dotsc,
y_n$ are sufficiently large.

\begin{lemma} \label{lem:estimate-VRone}
Let $B(z,v)$ denote the supremum of $\bigl\lvert Q' \bigl( v, \sigma_{I,u}(z) \bigr)
\bigr\rvert$, taken over $I \subseteq \{1, \dotsc, n\}$ and $u \in F^{\abs{I}} \cap
\HC$ with $\norm{u} \leq 1$.  Then there are constants $C > 0$ and $\alpha > 0$, such
that
\[
	Z(y,v) \leq C \cdot B(z,v)
\]
for every $v \in \VRone$ and every $z \in \HH^n$ with $y_j = \Im z_j \geq \alpha$ and
$0 \leq \Re z_j \leq 1$.
\end{lemma}

\begin{proof}
Given a vector $v \in \VC$, we let $v^{p,q} \in I^{p,q}(M,F)$ denote its components
relative to Deligne's decomposition.  As in \subsecref{subsec:proof}, we may replace
$v$ by $e^{-\sum x_j N_j'} v$ without affecting the statement we are trying to prove.
Setting $w = e^{-\Gamma'(s)} e^{-iN'} e^{-i \delta'} v$, we easily see that the norm
of each vector $(N')^k w^{p,q}$ with $p \leq -1$ is bounded by a constant times
$B(z,v)$. We again define
\[
	v' = e^{i \delta'} e^{iN'} e^{\Gamma'(s)} e^{-iN'} \cdot 
		\sum_{p \leq -1} e^{iN'} w^{p,q},
\]
and observe that $Z(y,v')$ is bounded by a fixed multiple of $B(z,h)$ by a version of
Lemma~\ref{lem:adN}. A useful observation is that $v' \in \HC$; this is because
$\Gr_0^W$ is of type $(0,0)$ at every point $z \in \HH^n$. By construction, we have
$v - v' \in e^{i \delta'} e^{i N'} e^{\Gamma'(s)} F^0$; since $v \in \VRone$, it is
therefore possible to write
\[
	v = v' + e^{i \delta'} e^{iN'} e^{\Gamma'(s)} (v_0 + h)
\]
for a unique choice of $h \in F^0 \cap \HC$. 

To continue, we let $g = v' + e^{i \delta'} e^{iN'} e^{\Gamma'(s)} v_0 - v_0$; note
that this vector belongs to $\HC$. Since $N' v_0 = 0$, and since $\norm{v_0}$ is
uniformly bounded due to Lemma~\ref{lem:v0}, we still have $Z(y,g)$ bounded by
a constant multiple of $B(z,v)$. We can now rewrite the equation from above as
$v - v_0	= g + e^{i \delta} e^{iN} e^{\Gamma(s)} h$. Remembering that $v - v_0$ is
a real vector, we obtain the relation
\[
	\overline{g} - g = e^{i \delta} e^{iN} e^{\Gamma(s)} h - 
		e^{-i \delta} e^{-iN} e^{\overline{\Gamma(s)}} \overline{h}.
\]			
From Proposition~\ref{prop:main-estimate}, we deduce that $Z(y,h)$ is bounded by a
constant times $Z(y,g)$, and hence by a constant multiple of $B(z,v)$, provided that
$y_1, \dotsc, y_n \geq \alpha$. But now the formula $v = v_0 + g + e^{i \delta}
e^{iN} e^{\Gamma(s)} h$, together with Lemma~\ref{lem:v0}, shows that the same is
true for $Z(y,v)$.
\end{proof}

Once again, this single inequality is all that one needs to prove Theorem~\ref{thm:closure-VRone}
\begin{proof}
The inequality in Lemma~\ref{lem:estimate-VRone} shows that $\norm{v(m)}$ remains
bounded as $m \to \infty$. Since $v(m) \in \VZone$, the sequence can take only
finitely many values, proving \ref{en:VRone-i}. We can then pass to a subsequence,
and assume for the remainder of the argument that $v(m) = v$ for some $v \in \VZone$.
Arguing as in the proof of Lemma~\ref{lem:subsequence}, we conclude from the
boundedness of $\sum y_j(m) N_j' v$ that $v$ satisfies \ref{en:VRone-ii}. We also
see that there is a further subsequence along which $\sum z_j(m) N_j' v = \sum w_j(m)
N_j' v$, where the sequence of $w(m) \in \CC^n$ converges to a vector $w \in \CC^n$.
This implies \ref{en:VRone-iii}.

Finally, we need to establish \ref{en:VRone-iv}. From the convergence of $Q' \bigl( v,
\sigma_{\{k\},u}(z(m)) \bigr)$ for $u \in F^1 \cap \HC$, we deduce as in the proof of
Theorem~\ref{thm:closure} that
\[
	N_k' v \in e^{\sum w_j N_j} e^{i \delta} (F^{-1} \cap \HC).
\]
This means that the vector $e^{-\sum \Re w_j N_j} N_k' v$ is a real Hodge class
of type $(-1,-1)$ in the mixed Hodge structure $\bigl( M \cap H, e^{i \delta + i
\sum \Im w_j N_j} (F \cap H) \bigr)$. Lemma~\ref{lem:MHS-Hodge} implies that $N_k' v$
lies in the kernel of the operator $\delta + \sum \Im w_j N_j$, and that
\[
	e^{-\sum w_j N_j} e^{-i \delta} N_k' v = e^{-\sum \Re w_j N_j} N_k' v = N_k' v^{0,0}.
\]
This gives \ref{en:VRone-iv} and concludes the proof.
\end{proof}

\begin{lemma} \label{lem:MHS-Hodge}
Let $(W, F)$ be an $\RR$-mixed Hodge structure, and let $v \in W_{2p} \cap F^p$ be a
real Hodge class of type $(p,p)$. Let $\delta \in L_{\RR}^{-1,-1}(W, F)$ be the
unique element for which $\bigl( W, e^{-i \delta} F \bigr)$ splits over $\RR$. Then
$\delta v = 0$, and consequently $v \in I^{p,p} \bigl( W, e^{-i \delta} F \bigr)$.
\end{lemma}

\begin{proof}
Since $v$ defines a morphism of $\RR$-mixed Hodge structures $\RR(-p) \to
(W, F)$, the functoriality of $\delta$ (see \cite{KNU-SL}*{Lemma~1.6} for a
proof) implies that $\delta v = 0$. It follows that $v$ is also a real Hodge class of
type $(p,p)$ in $\bigl( W, e^{-i \delta} F \bigr)$.
\end{proof}

Recall that we defined $\Tnu$ as the subset of the \'etal\'e space of $\shVZ$,
consisting of those points that map to $1 \in \ZZ$. Theorem~\ref{thm:closure-VRone}
is strong enough to conclude that $\Tnu$ has an analytic closure inside of $T(F_0
\Mmod)$. Note that Corollary~\ref{cor:main-condition}, to the effect that $\TZ \subseteq
T(F_0 \Mmod)$ is closed analytic, can be viewed as the special case $\nu = 0$.

\begin{theorem} \label{thm:closure-Tnu}
The topological closure of $\Tnu$ inside $T(F_0 \Mmod)$ is an analytic subset.
\end{theorem}

\begin{proof}
We shall use both the space $T(F_0 \Mmod)$, as well as the space $T(F_0 \shHOe)$
coming from the canonical extension; since $F_0 \shHOe \subseteq F_0 \Mmod$, they are
related by a holomorphic mapping $g
\colon T(F_0 \Mmod) \to T(F_0 \shHOe)$. We also have holomorphic mappings $\varphi
\colon \Tnu \to T(F_0 \Mmod)$ and $\psi \colon \Tnu \to T(F_0 \shHOe)$ with $\psi = g
\circ \varphi$. Let $\Tnu(v)$ denote the connected component of $\Tnu$ containing a
given vector $v \in \VZone$. It suffices
to show that the image of each $\Tnu(v)$ under the holomorphic mapping $\varphi$ has
an analytic closure; this is because \ref{en:VRone-i} in
Theorem~\ref{thm:closure-VRone} assures us that only finitely many of these image
closures can meet at any given point of $T(F_0 \Mmod)$.

Fix a vector $v \in \VZone$. We may clearly assume for the remainder of the argument
that $v$ satisfies \ref{en:VRone-ii}--\ref{en:VRone-iv} in
Theorem~\ref{thm:closure-VRone}, for otherwise, the image of $\Tnu(v)$ is already
closed in a neighborhood of $0 \in \Delta^n$ and there is nothing to prove.
In particular, we have $a_1 N_1' v + \dotsb + a_n N_n' v = 0$, and $e^{-\sum
w_j N_j} e^{-i \delta} N_k' v = N_k' v^{0,0}$ for some $w \in \CC^n$. Replacing $v$
by $e^{-\sum \Re w_j N_j'} v$, we arrange that $w = 0$, at the cost of having $v \in
\VRone$.

We can use this information to show that the image of $\Tnu(v)$ under $\psi$ has an
analytic closure inside $T(F_0 \shHOe)$.  Let $D = \Delta^n - \dstn{n}$, let $p
\colon T(F_0 \shHOe) \to \Delta^n$ denote the projection map, and set $E =
p^{-1}(D)$. Clearly, $\psi(\Tnu(v))$ is a closed analytic subset of
$T(F_0 \shHOe) - E$, of pure dimension $n$. To prove that its closure remains
analytic, it suffices to show that the intersection of the closure with $E$ is
contained in a countable union of images of complex manifolds of dimension at most
$n-1$. Indeed, this implies that the intersection has $2n$-dimensional Hausdorff
measure equal to zero, and we conclude by a result of E.~Bishop's
\cite{Bishop}*{Lemma~9}.

By construction of the canonical extension, the mapping $\psi \colon \Tnu(v) \to
T(F_0 \shHOe)$ is given in coordinates by the formula
\[
	\psi \colon \HH^n \to \Delta^n \times \Hom(F^0 \HC, \CC), \qquad
		(z_1, \dotsc, z_n) \mapsto \bigl( e^{2 \pi i z_1}, \dotsc, e^{2 \pi i z_n}, 
		f_z \bigr)
\]
where $f_z \colon F^0 \HC \to \CC$ is the linear functional
\[
	u \mapsto f_z(u) = Q' \bigl( v, e^{i \delta} e^{\sum z_j N_j} e^{\Gamma(s)} u \bigr).
\]
Now we compute that
\[
	e^{-\Gamma'(s)} e^{-\sum z_j N_j'} e^{-i \delta'} v =
		e^{-\Gamma'(s)} e^{-i \delta'} v + \sum_{k=1}^{\infty} (-1)^k e^{-\Gamma(s)} 
			(z_1 N_1' + \dotsb + z_n N_n')^k v^{0,0},
\]
and hence we have
\begin{align*}
	f_z(u) &= Q' \bigl( e^{-\Gamma'(s)} e^{-\sum z_j N_j'} e^{-i \delta'} v, u \bigr) \\
	&= Q' \bigl( e^{-\Gamma'(s)} e^{-i \delta'} v, u \bigr)
		+ \sum_{k=1}^{\infty} (-1)^k Q' \bigl( e^{-\Gamma(s)} (z_1 N_1' +
		\dotsb + z_n N_n')^k v^{0,0}, u \bigr)
\end{align*}
for every $u \in F^0 \HC$. In particular, remembering that $\Gamma(s) \in \qlie$, we
find that when $u$ belongs to the subspace $I^{1,1}(M \cap H, F \cap H)$, then
\begin{equation} \label{eq:fzu}
	f_z(u) = Q' \bigl( e^{-\Gamma'(s)} e^{-i \delta'} v, u \bigr)
		- Q' \bigl( z_1 N_1' v^{0,0} + \dotsb + z_n N_n' v^{0,0}, u \bigr).
\end{equation}

It is now easy to determine the points in the closure. Fix a subset $I \subseteq \{1,
\dotsc, n\}$ of size $k$, and consider the stratum $D_I \subseteq D$ of points with
$s_j \neq 0$ for $j \in I$ and $s_j = 0$ for $j \not\in I$; note that $\dim D_I = k$.
Suppose that $z(m) \in \HH^n$ is a sequence of points for which $\psi(z(m))$
converges to a point over $s_0 \in D_I$. We have $N_j' v = 0$ for $j \in I$, and since
$a_1 N_1' v + \dotsb + a_n N_n' v = 0$, the span of $N_1' v, \dotsc, N_n' v$
has dimension at most $n-k-1$. Using the convergence of $f_{z(m)}(u)$ in
\eqref{eq:fzu} and arguing as in the proof of Lemma~\ref{lem:subsequence}, we see
that $\sum z_j(m) N_j' v^{0,0}$ converges to $\sum w_j N_j' v^{0,0}$ for some vector
$w \in \CC^n$. Thus every limit point over $s_0 \in D_I$ is of the form $(s_0, f)$,
where $f \colon F^0 \HC \to \CC$ is given by the formula
\[
	f(u) = Q' \bigl( e^{-\Gamma'(s_0)} e^{-i \delta'} v, u \bigr)
		+ \sum_{k=1}^{\infty} (-1)^k Q' \bigl( e^{-\Gamma(s_0)} (w_1 N_1' +
		\dotsb + w_n N_n')^k v^{0,0}, u \bigr)
\]
for some choice of $w \in \CC^n$. Evidently, such points are parametrized by a linear
space of dimension at most $n-k+1$. It follows that the intersection of $p^{-1}(D_I)$
with the closure of $\Tnu(v)$ is contained in a complex-analytic subset of dimension
at most $\dim D_I + (n-1+k) = n-1$; as explained above, this suffices to conclude
that $\psi(\Tnu(v))$ has an analytic closure inside $T(F_0 \shHOe)$.

To finish the proof, we observe that the preimage of $\overline{\psi(\Tnu(v))}$ under
$g$ is an analytic subset of $T(F_0 \Mmod)$ whose intersection with $T(F_0 \shHO)$
equals $\varphi(\Tnu(v))$. By a well-known result in complex analysis,
this implies that the closure of $\varphi(\Tnu(v))$ inside $T(F_0
\Mmod)$ is itself analytic, and thereby concludes the proof.
\end{proof}

Similarly, the graph of the normal function $\nu \colon X \to J(\shH)$ has an
analytic closure inside of $\Jb(\shH)$.

\begin{corollary} \label{cor:closure-Gammanu}
Let $\Gammanu \subseteq J(\shH)$ denote the graph of an admissible normal function
$\nu \colon \dstn{n} \to J(\shH)$. Then the closure of $\Gammanu$
inside of $\Jb(\shH)$ is analytic.
\end{corollary}

\begin{proof}
Since the quotient map $T(F_0 \Mmod) \to \Jb(\shH)$ is open by Lemma~\ref{lem:quotients-3}, 
this follows immediately from Theorem~\ref{thm:closure-Tnu}.
\end{proof}

When the singularity of $\nu$ is a nonzero torsion class, then the graph of $\nu$ is
already closed (this observation is due to P.~Brosnan).

\begin{lemma} \label{lem:torsion}
Suppose that the singularity of $\nu$ at $0 \in \Delta^n$ is a nonzero torsion class.
Then the closure of $\Gammanu$ in $\Jb(\shH)$ does not meet $p^{-1}(0)$.
\end{lemma}

\begin{proof}
It suffices to show that if the singularity of $\nu$ is torsion, but the closure of
$\Tnu$ in $T(F_0 \Mmod)$ contains a point over $0 \in \Delta^n$, then $\nu$ is
actually nonsingular. If there is such a point in the closure, then by
Theorem~\ref{thm:closure-VRone}, there exists a class $v \in \VZone$ that
satisfies the conditions listed there; in particular, $a_1 N_1' v + \dotsb + a_n N_n'
v = 0$ for positive integers $a_1, \dotsc, a_n$. Consider now the complex
\begin{diagram}
\HQ &\rTo& \bigoplus_j N_j(\HQ) &\rTo& \bigoplus_{j < k} N_j N_k(\HQ) &\rTo& \dotsm
\end{diagram}
that computes the intersection cohomology of the local system $\shHQ$ on $\Delta^n$
\cite{CKS-L2}*{p.~219}. With rational coefficients, the singularity of $\nu$ is
represented by the class of $(N_1' v, \dotsc, N_n' v)$ in the first cohomology group
of the complex. Since the singularity is torsion, there is a vector $h \in \HQ$ such
that $N_j' v = N_j h$ for all $j$. Now $a_1 N_1 h + \dotsb + a_n N_n h = 0$, and so
we have $h \in M_{-1} \cap \HQ$, which implies that $N_j' v \in M_{-3} \cap \HQ$.
Assertion \ref{en:VRone-iv} in Theorem~\ref{thm:closure-VRone} then forces 
$N_j' v = 0$ for all $j=1, \dotsc, n$, and so $\nu$ has no singularities on $\Delta^n$.
\end{proof}

%%%%%%%%%%%%%%%%%%%%%%%%%%%
\section{Admissible normal functions}

\subsection{A lemma about holonomic modules}
\label{subsec:holonomic}

In order to extend nonsingular admissible normal functions on $X$ to holomorphic sections of
$\Jb(\shH) \to \Xb$, we need a more general version of Lemma~\ref{lem:shHC-action}
that works for arbitrary holonomic $\Dmod$-modules. This section describes such a
generalization. Let $\Nmod$ be a holonomic left $\Dmod$-module
on a complex manifold $X$ of dimension $n$. Then the de Rham complex (in degrees
$-n, \dotsc, 0$)
\begin{diagram}
	\DR(\Nmod) = \Big\lbrack \Nmod &\rTo& \OmX{1} \tensor \Nmod &\rTo& \OmX{2} \tensor \Nmod
		&\rToDots& \OmX{n} \tensor \Nmod \Big\rbrack \decal{n}
\end{diagram}
is a perverse sheaf on $X$ by M.~Kashiwara's theorem, and we let $\shHC^k = \shH^k
\bigl( \DR(\Nmod) \bigr)$. Let $\DD = \shExt_{\Dmod_X}^n \bigl( \argbl, \Dmod_X
\tensor \omega_X^{-1} \bigr)$ be the duality functor for holonomic left
$\Dmod$-modules, and define $\Nmod' = \DD(\Nmod)$. The following result is well-known.

\begin{lemma} \label{lem:sol-DR}
For each $k \in \ZZ$, we have $\shExt_{\Dmod_X}^k(\Nmod', \OX) \simeq \shH^{k-n}
\bigl( \DR(\Nmod) \bigr)$.
\end{lemma}

A simple consequence is the isomorphism
\[
	\shHom_{\Dmod_X}(\Nmod', \OX) \simeq \ker \bigl( \nabla \colon \Nmod \to \OmX{1}
		\tensor \Nmod \bigr) \simeq \shHC^{-n}.
\]
Let $j \colon U \into X$ be an open subset over which $\Nmod$ is a flat vector
bundle; equivalently, such that $j^{-1} \shHC^{-n}$ is a local system. Since sections
of $j^{-1} \shHC^{-n}$ are flat, we then have an injection
\begin{equation} \label{eq:gen-action}
	\jl \bigl( j^{-1} \shHC^{-n} \bigr) \into \shHom_{\Dmod_X}(\Nmod', \OX).
\end{equation}

To relate this with Lemma~\ref{lem:shHC-action}, suppose that we have a holonomic
$\Dmod$-module $\Mmod$ that comes from a polarized Hodge module $M$ of weight $w$.
Being polarized, $M$ satisfies $\DD(M) \simeq M(w)$, which implies that $\Mmod'
\simeq \Mmod$. Let $j \colon U \into X$ be an open subset on which $M$ is the Hodge
module associated to a polarized variation of Hodge structure $\shH$ (of weight
$w-n$). The injective map from above now becomes
\[
	\jl \shHC \into \ker \bigl( \nabla \colon \Mmod \to \OmX{1} \tensor \Mmod \bigr)
		\simeq \shHom_{\Dmod_X}(\Mmod, \OX).
\]
Let $Q$ be the polarization of the variation; then the map is given by sending a flat
section $h$ of $\jl \shHC$ to the map of $\Dmod$-modules $Q(h, \argbl) \colon \Mmod
\to \OX$. This follows from M.~Saito's construction of the isomorphism $\Mmod \simeq
\Mmod'$ in terms of the polarization \cite{Saito-HM}.

\subsection{Extending admissible normal functions without singularities}
\label{subsec:nub}

We now look at the problem of extending admissible normal functions with no
singularities to holomorphic sections of the space $\Jb(\shH) \to \Xb$.

Let $\nu$ be a normal function on $X$ for the variation $\shH$, admissible relative
to $\Xb$. By M.~Saito's theory \cite{Saito-ANF}*{p.~243}, it corresponds to a mixed
Hodge module $\Nnu$ on $\Xb$, with $W_{n-1} \Nnu = M$, and $\Gr_n^W \Nnu$ the trivial
Hodge module of weight $n$. On $X$, we have an extension of integral local systems
\begin{equation} \label{eq:ext-shVZ} \begin{diagram} 0 &\rTo& \shHZ &\rTo& \shVZ
&\rTo& \ZZX &\rTo& 0, \end{diagram} \end{equation} and therefore a cohomology class
$\class{\nu} \in H^1(X, \shHZ)$. Using the Leray spectral sequence for the inclusion
$j \colon X \into \Xb$, we obtain an exact sequence \begin{diagram} 0 &\rTo& H^1
\bigl( \Xb, \jl \shHZ \bigr) &\rTo& H^1(X, \shHZ) &\rTo& H^0 \bigl( \Xb, R^1 \jl
\shHZ \bigr).  \end{diagram} The following concept has been introduced by M.~Green
and P.~Griffiths \cite{GG1}.  \begin{definition} The image of $\class{\nu}$ in the
space $H^0 \bigl( \Xb, R^1 \jl \shHZ \bigr)$ is called the \define{singularity} of
the normal function $\nu$. When the image is zero, we shall say that $\class{\nu}$ is
locally trivial, or that $\nu$ has no singularities.  \end{definition}

When $\nu$ has no singularities, we evidently have $\class{\nu} \in H^1 \bigl( \Xb,
\jl \shHZ \bigr)$.  The relationship of these definitions with \eqref{eq:ext-shVZ} is
the following: Taking direct images, we have a long exact sequence
\begin{diagram}[midshaft] 0 &\rTo& \jl \shHZ &\rTo& \jl \shVZ &\rTo& \ZZXb
&\rTo^{\delta}& R^1 \jl \shHZ &\rTo& \dotsb, \end{diagram} and local triviality of
$\class{\nu}$ is equivalent to the vanishing of the connecting homomorphism $\delta$.
Thus if the normal function has no singularities, we obtain from it an extension of
sheaves of abelian groups on $\Xb$, namely \begin{equation} \label{eq:ext-NF-ZZ}
\begin{diagram} 0 &\rTo& \jl \shHZ &\rTo& \jl \shVZ &\rTo& \ZZXb &\rTo& 0.
\end{diagram} \end{equation}

On the other hand, the mixed Hodge module $\Nnu$ is part of an extension
\begin{diagram} 0 &\rTo& M &\rTo& \Nnu &\rTo& \QHXb \decal{n} &\rTo& 0, \end{diagram}
with $\QHXb \decal{n}$ the trivial Hodge module of weight $n$ on $\Xb$. 
% By passing to the underlying filtered $\Dmod$-modules, we obtain an extension
% \begin{equation} \label{eq:ext-NF-Dmod-1} \begin{diagram} 0 &\rTo& (\Mmod, F)
% &\rTo& (\Nmodnu, F) &\rTo& \bigl( \OXb, F \bigr) &\rTo& 0; \end{diagram}
% \end{equation} as usual, $\OXb$ has the trivial filtration, with $\Gr_k^F \OXb = 0$
% unless $k = 0$.  By strictness, the sequence remains exact upon passing to graded
% quotients for the Hodge filtrations.
Let $\DD(\argbl)$ denote the Verdier dual on the category of mixed Hodge modules.
Since $M$ is polarized, and of weight $n-1$, we have $\DD(M) \simeq M(n-1)$; also,
$\DD \bigl( \QHXb \decal{n} \bigr) \simeq \QHXb \decal{n}(n)$. Dualizing the
extension above, and applying a Tate twist, we thus get an exact sequence
\begin{diagram} 0 &\rTo& \QHXb \decal{n} &\rTo& \Nnu' &\rTo& M(-1) &\rTo& 0,
\end{diagram} with $\Nnu' = \DD(\Nnu)(-n)$. Let $(\Nmodnu', F)$ be the underlying
filtered $\Dmod$-module; then we have an extension of filtered $\Dmod$-modules
\begin{equation} \label{eq:ext-NF-Dmod-2} \begin{diagram} 0 &\rTo& \bigl( \OXb, F
\bigr) &\rTo& \bigl( \Nmodnu', F \bigr) &\rTo^{\alpha}& \bigl( \Mmod, F \decal{-1}
\bigr) &\rTo& 0.  \end{diagram} \end{equation} Morphisms of mixed Hodge modules are
strictly compatible with the Hodge filtrations on the underlying $\Dmod$-modules;
because $F_{-1} \OXb = 0$, it follows that $\alpha$ induces an isomorphism $F_{-1}
\Nmodnu' \simeq F_0 \Mmod$.

Just as in \subsecref{subsec:Jacobians}, we can now compare the two extensions in
\eqref{eq:ext-NF-ZZ} and \eqref{eq:ext-NF-Dmod-2} to obtain a section of $T(F_0
\Mmod) / \TZ$. But note that this has to be done carefully, since $F_0 \Mmod$ is in
general not locally free near points of $\Xb - X$.

\begin{proposition} \label{prop:nub} Any admissible normal function $\nu \colon
J(\shH) \to X$ without singularities can be canonically extended to a holomorphic
section of $\Jb(\shH) \to \Xb$.  \end{proposition}

\begin{proof} Since $\nu$ has no singularities, it gives rise to an extension of
sheaves of abelian groups as in \eqref{eq:ext-NF-ZZ}. Now cover the space $\Xb$ by
open subsets $U_i$, such that \eqref{eq:ext-NF-ZZ} is locally split. This means that
we have $v_i \in H^0 \bigl( U_i, \jl \shVZ \bigr)$, mapping to the constant section
$1 \in H^0(U_i, \ZZ)$; it follows that $v_j - v_i = h_{ij}$ for certain $h_{ij} \in
H^0 \bigl( U_i \cap U_j, \jl \shHZ \bigr)$. Note that $v_i$ is well-defined up to a
section of $\jl \shHZ$ on $U_i$.

Let $j_i \colon U_i \into \Xb$ be the inclusion maps.  Because of
\eqref{eq:gen-action}, each $v_i$ defines a map of $\Dmod$-modules $\phi_i \colon
j_i^{-1} \Nmodnu' \to \shO_{U_i}$.  Moreover, since $v_j - v_i = h_{ij}$, we have
$\phi_j - \phi_i = Q(h_{ij}, \argbl) \circ \alpha$. Restricting to $F_{-1} \Nmodnu'
\simeq F_0 \Mmod$, and noting that the isomorphism is given by $\alpha$, we thus have
holomorphic sections \[ \psi_i \in H^0 \bigl( U_i, (F_0 \Mmod)^{\vee} \bigr) \quad
\text{satisfying} \quad \psi_j - \psi_i = Q(h_{ij}, \argbl).  \] By definition of the
analytic structure on $T(F_0 \Mmod) / \TZ$, this means exactly that we have a
holomorphic section of $\Jb(\shH) \to \Xb$. It is clear from the construction that
the section is independent of the choices made. That we recover the original normal
function on $X$ is a straightforward consequence of Lemma~\ref{lem:J1-J2}.
\end{proof}

From now on, let us write $\nub$ for the section of $\Jb(\shH) \to \Xb$ constructed
in the lemma; we refer to it as the \define{extension} of the original normal
function $\nu$. It would be interesting to know the set of points in $\Jb(\shH)$ that
can lie on the graph of an extended normal function.

\subsection{The horizontality condition} \label{subsec:horizontality}

It is clear that the extension $\nub$ constructed in Proposition~\ref{prop:nub} is
far from being an arbitrary section of the quotient. In fact, the proof shows that
there are local liftings $\psi \colon F_0 \Mmod \restr{U} \to \shO_U$ that are
compatible with differentiation: for any $k \geq 0$, any differential operator $D \in
H^0(U, F_k \Dmod_U)$, and any section $s \in H^0(U, F_{-k} \Mmod)$, the lifting
satisfies $\psi(Ds) = D \psi(s)$.  This appears to be the correct notion of
horizontality for sections of $\Jb(\shH) \to \Xb$.

\begin{definition} A holomorphic section of $\Jb(\shH) \to \Xb$ will be called
\define{horizontal} if it admits local holomorphic liftings $\psi \colon F_0 \Mmod
\restr{U} \to \shO_U$ with the property that $\psi(\xi s) = d_{\xi} \bigl( \psi(s)
\bigr)$ for any holomorphic tangent vector field $\xi \in H^0(U, \Theta_U)$ and any
section $s \in H^0(U, F_{-1} \Mmod)$ \end{definition}

It follows that $\psi(Ds) = D \psi(s)$ for $D \in H^0(U, F_k \Dmod_U)$ and $s \in
H^0(U, F_{-k} \Mmod)$ as above. Over $X$, the definition clearly recovers the usual
definition of horizontality.  We now prove the converse to
Proposition~\ref{prop:nub}.

\begin{proposition} \label{prop:sections} Let $\mu \colon \Xb \to \Jb(\shH)$ be a
holomorphic section that is horizontal. Then $\mu$ is the extension of an admissible
normal function on $X$ with locally trivial cohomology class.  \end{proposition}

\begin{proof} The restriction of $\mu$ to $X$ is a horizontal and holomorphic section
of $J(\shH)$, and therefore a normal function $\nu$. We have to prove that it is
admissible, and that its cohomology class is locally trivial. To begin with the
latter, consider the exact sequence of sheaves \begin{diagram} 0 &\rTo& \jl \shHZ
&\rTo& (F_0 \Mmod)^{\vee} &\rTo& (F_0 \Mmod)^{\vee} / \jl \shHZ &\rTo& 0,
\end{diagram} and note that the quotient is the sheaf of sections of $\Jb(\shH)$. Via
the connecting homomorphism, the section $\mu$ determines an element in $H^1 \bigl(
\Xb, \jl \shHZ \bigr)$, whose image in $H^1(X, \shHZ)$ is the class of the normal
function. By construction, $\class{\nu}$ goes to zero in $H^0 \bigl( \Xb, R^1 \jl
\shHZ \bigr)$; this means that $\nu$ has no singularities.

Since admissibility is defined by a curve test \cite{Kashiwara-study}, we let $f
\colon \Delta \to \Xb$ be an arbitrary holomorphic curve with $f(\dst) \subseteq X$,
such that $\shH' = \fu \shH$ has unipotent monodromy. By
Proposition~\ref{prop:Jb-functoriality}, we have a holomorphic map \[ \Delta
\times_{\Xb} \Jb(\shH) \to \Jb(\shH') \] over $\Delta$, and so $\mu$ induces a
holomorphic section of $\Jb(\shH')$ whose restriction to $\dst$ is the pullback of
the normal function. Since it suffices to prove the admissibility of the latter, we
have reduced the problem to the case of a disk, where we can apply the following
lemma.  \end{proof}

\begin{lemma} Let $\shH$ be a polarized variation of Hodge structure of weight $-1$
on $\dst$, whose monodromy is unipotent. Then any holomorphic and horizontal section
of $\Jb(\shH) \to \Delta$ is the extension of an admissible normal function.
\end{lemma}

\begin{proof} Shrinking the radius of the disk, if necessary, we may assume that the
section can be lifted to a map $\psi \colon F_0 \Mmod \to \shO$ that satisfies the
condition in the definition of horizontality. As before, it defines a normal function
with trivial cohomology class on $\dst$. Let \begin{diagram} 0 &\rTo& \shH &\rTo&
\shV &\rTo& \ZZ(0) &\rTo& 0 \end{diagram} be the corresponding extension of
variations of mixed Hodge structure on $\dst$; we need to show that $\shV$ is
admissible.

At this point, we do not know that $\shV$ can be extended to a mixed Hodge module on
$\Delta$---in fact, this is equivalent to admissibility by \cite{Saito-ANF}*{p.~243}.
Nevertheless, we can use $\psi$ to reconstruct at least the filtered $\Dmod$-module
$(\Nmodnu, F)$. Since the cohomology class of $\nu$ should be trivial, we define the
$\Dmod$-module $\Nmodnu' = \ODelta \oplus \Mmod$, and introduce a filtration on it by
setting \[ F_p \Nmodnu' = \begin{cases} \im \bigl( (\psi, \id) \colon F_{p+1} \Mmod
\to \ODelta \oplus \Mmod \bigr) & \text{for $p \leq -1$,} \\ \ODelta \oplus F_{p+1}
\Mmod & \text{for $p \geq 0$.} \end{cases} \] The condition on $\psi$ ensures that
the filtration is good, and therefore that $(\Nmodnu', F)$ is a filtered
$\Dmod$-module. Because we should have $\Nnu = \DD(\Nnu')(-1)$, we now define the
filtered $\Dmod$-module $(\Nmodnu, F)$ by dualizing; evidently, $\Nmodnu = \Mmod
\oplus \ODelta$, and after a short calculation, one finds that $F_p \Nmodnu = F_p
\Mmod$ for $p < 0$, while \[ F_p \Nmodnu = \menge{(h,f) \in \Mmod \oplus \ODelta}%
{\text{$Q(h,s) = f \psi(s)$ for every $s \in F_{-p} \Mmod$}} \] for $p \geq 0$.

Now we verify admissibility. Note first that each sheaf $F_p \Mmod$ is locally free
(since $\Mmod$ is always torsion-free, and $\Delta$ has dimension one). Let $\shHOe$
be Deligne's canonical extension the flat vector bundle; we then have $\shHOe \into
\Mmod$. Evidently, the canonical extension of $\shVO$ is given by $\shHOe \oplus
\ODelta \into \Nmodnu$. The description of $F_p \Nmodnu$ above shows that the Hodge
bundles $F_p \shVO$ extend to holomorphic subbundles of the canonical extension,
which is one of the two conditions for admissibility. The second one, existence of
the relative weight filtration, is trivially satisfied because the underlying local
system $\shVZ = \shHZ \oplus \ZZ_{\dst}$ is a direct sum.  \end{proof}

\subsection{Graphs of admissible normal functions} \label{subsec:graphs}

In this section, we consider admissible normal functions on $X$ with possibly
nontrivial singularities. By Proposition~\ref{prop:sections}, such a normal function
cannot be extended to a section of $\Jb(\shH) \to \Xb$. Nevertheless, the following
surprising result is true.

\begin{theorem} \label{thm:graphs} Let $\nu \colon X \to J(\shH)$ be a normal
function, admissible relative to $\Xb$.  Then the topological closure of the graph of
$\nu$ is an analytic subset of $\Jb(\shH)$.  \end{theorem}

\begin{proof} This follows from the corresponding statement in the normal crossing
case, contained in Corollary~\ref{cor:closure-Gammanu}, by the same argument as in
\subsecref{subsec:reduction}.  \end{proof}

One consequence is an alternative proof for Conjecture~\ref{conj:GG}. It is quite
different from the existing proof by P.~Brosnan and G.~Pearlstein \cite{BP3}, but
similar in spirit to the treatment of the one-variable case in M.~Saito's paper
\cite{Saito-GGK}.

\begin{corollary} If a normal function $\nu \colon X \to J(\shH)$ is admissible
relative to $\Xb$, then the closure of its zero locus $Z(\nu)$ is an analytic subset
of $\Xb$. In particular, when $X$ is an algebraic variety, the zero locus $Z(\nu)$ is
an algebraic subvariety.  \end{corollary}

\begin{proof} The closure of $Z(\nu)$ is contained in the intersection of the closure
of the graph of $\nu$ with the zero section of $\Jb(\shH)$, and is therefore analytic
as well.  When $X$ is an algebraic variety, we take $\Xb$ to be
projective---admissibility is independent of the choice of compactification in that
case---and the algebraicity of $Z(\nu)$ follows by Chow's Theorem.  \end{proof}

We also note the following property of normal functions with torsion singularities,
suggested by P.~Brosnan. In the statement, $p \colon \Jb(\shH) \to \Xb$ is the
projection map, and $\Gammanu \subseteq J(\shH)$ is the graph of $\nu$.

\begin{proposition} \label{prop:torsion} Suppose that $\nu \colon X \to J(\shH)$ is
an admissible normal function, whose singularity at a point $x \in \Xb - X$ is
torsion. If the closure of $\Gammanu$ meets $p^{-1}(x)$, then $\nu$ extends to a
holomorphic section of $\Jb(\shH)$ in a neighborhood of $x$.  \end{proposition}

\begin{proof} This follows from Lemma~\ref{lem:torsion}, by a similar argument as in
\subsecref{subsec:curves}. Namely, let \begin{diagram} 0 &\rTo& M &\rTo& \Nnu &\rTo&
\QHXb \decal{n} &\rTo& 0 \end{diagram} be the extension of mixed Hodge modules on
$\Xb$ corresponding to the normal function $\nu$, let $i \colon \{x\} \to \Xb$ be the
inclusion, and set $V = H^{-n} \iu \Nnu$ and $H = H^{-n} \iu M$. Then both $H$ and
$V$ are mixed Hodge structures, defined over $\ZZ$, and since the singularity of
$\nu$ at $x$ is torsion, we get $\VQ / \HQ \simeq \QQ$. Also note that $\HC / F_0 \HC
\to \VC / F_0 \VC$ is an isomorphism by strictness.

Now suppose that the closure of $\Gammanu$ meets $p^{-1}(x)$. Let $\Tnu \subseteq
T(F_0 \Mmod)$ be the preimage of the graph. Since the closure of $\Tnu$ is analytic,
we can find a holomorphic curve $g \colon \Delta \to T(F_0 \Mmod)$ with $g(\dst)
\subseteq \Tnu$, such that $g(0)$ lies over the point $x$. Let $\nu'$ denote the
pullback of $\nu$ to $\dst$; then $\nu'$ has trivial singularity on $\Delta$, and so
the corresponding extension of mixed Hodge modules \begin{diagram} 0 &\rTo& M' &\rTo&
N_{\nu'} &\rTo& \QHDelta \decal{1} &\rTo& 0 \end{diagram} splits over $\ZZ$. If we
let $i_0 \colon \{0\} \to \Delta$, $H' = H^{-1} \iku{0} M'$, and $V' = H^{-1} \iku{0}
N_{\nu'}$, we have $\VZ' \simeq \HZ' \oplus \ZZ$. As in \subsecref{subsec:curves}, we
obtain a commutative diagram with exact rows: \begin{diagram} 0 &\rTo& \HZ &\rTo& \VZ
&\rTo& \ZZ &\rTo& \dotsm \\ & & \dTo && \dTo && \dEqual \\ 0 &\rTo& \HZ' &\rTo& \VZ'
&\rTo& \ZZ &\rTo& 0 \end{diagram}

Let $v' \in \VZ'$ be a lifting of $1 \in \ZZ$. It determines a point in $T(F_0
\Mmod')_0$, and hence in the quotient $\HC' / F_0 \HC' \simeq \VC' / F_0 \VC'$. As
before, the fact that we have $g(0) \in T(F_0 \Mmod)_x$, together with the
compatibility of the maps, implies that \[ v' \in F_0 \VC' + \im(\VC \to \VC').  \]
Since the singularity of $\nu$ is a torsion class, the quotient $V'/V$ is a mixed
Hodge structure of weight $\leq -1$. We again conclude that $v' \in \VQ$, and hence
$v' \in \VZ$. But then $\VZ \simeq \HZ \oplus \ZZ$, and so $\nu$ has no singularity
at $x$.  \end{proof}

\subsection{A N\'eron model for torsion singularities} \label{subsec:torsion}

The analytic space $\Jb(\shH)$ has all the properties that are expected for the
identity component of the N\'eron model. In this section, we extend the construction
to produce an analytic space that graphs admissible normal functions with torsion
singularities. This generalizes work by M.~Saito \cite{Saito-GGK} in the case where
$\dim \Xb = 1$.

\begin{theorem} \label{thm:Neron-torsion} There is an analytic space $\Jbtor(\shH)
\to \Xb$, whose horizontal and holomorphic sections are precisely the admissible
normal functions with torsion singularities.  \end{theorem}

We obtain the space $\Jbtor(\shH)$ by a gluing construction as in
\cite{BPS}*{Section~2.3}; local models are given by locally defined admissible normal
functions with torsion singularities.  To introduce some notation, suppose that we
have an open subset $U \subseteq \Xb$, and an admissible normal function $\nu$ on $U
\cap X$ with only torsion singularities. Then $\nu$ defines a section of $p \colon
\Jb(\shH) \to \Xb$ over $U \cap X$. By Proposition~\ref{prop:torsion}, there is a
maximal open subset $U(\nu) \subseteq U$ to which this section can be extended; the
important fact is that the graph of $\nu \colon U(\nu) \to p^{-1}(U)$ is a closed
analytic subset.  For every such pair, we let $Y(U, \nu)$ be a copy of $p^{-1}(U)
\subseteq \Jb(\shH)$, and write $p \colon Y(U,\nu) \to U$ for the projection map.

Let $Y$ be the disjoint union of all the spaces $Y(U, \nu)$, and define an
equivalence relation on $Y$ by setting \[ y \sim y' \quad \text{if and only if} \quad
\left\{ \begin{aligned} &\text{$x = p(y) = p(y')$ lies in $U(\nu) \cap U'(\nu')$,} \\
&\text{and $y + \nu(x) = y' + \nu'(x)$} \end{aligned}\right.  \] for $y \in Y(U,
\nu)$ and $y' \in Y(U', \nu')$. 

\begin{lemma} The quotient map $q \colon Y \to Y/\sim$ is an open map, and the
topology on $Y/\sim$ is Hausdorff.  \end{lemma}

\begin{proof} To prove that $q$ is an open map, it suffices to show that the image of
each $Y(U,\nu)$ is an open subset of the quotient. One easily sees that the preimage
of $q \bigl( Y(U,\nu) \bigr)$ intersects $Y(U',\nu')$ in the open subset $p^{-1}
\bigl( U(\nu) \cap U'(\nu') \bigr)$; this implies the first assertion

For the second, it is again enough to prove that $\sim$ defines a closed subset of $Y
\times Y$. So suppose that we have two sequences of points $y_n \in Y(U,\nu)$ and
$y_n' \in Y(U',\nu')$ with $y_n \sim y_n'$ for all $n \in \NN$, such that $(y_n,y_n')
\to (y,y')$. Letting $x_n = p(y_n) = p(y_n')$, we obtain $x_n \to x$, where $x = p(y)
= p(y')$. Since the graphs of $\nu$ and $\nu'$ are closed by
Proposition~\ref{prop:torsion}, we have $x \in U(\nu) \cap U'(\nu')$. The continuity
of $\nu$ and $\nu'$ now implies that $y + \nu(x) = y' + \nu'(x)$, proving that $y
\sim y'$.  \end{proof}

Here is the proof of Theorem~\ref{thm:Neron-torsion}.

\begin{proof} In the notation from above, let $\Jbtor(\shH) = Y / \sim$, with the
obvious projection map to $\Xb$. Evidently, no two distinct points of $Y(U, \nu)$ are
identified by the equivalence relation, and so $Y(U, \nu)$ is isomorphic to its image
in the quotient. Since the quotient is Hausdorff, it follows that it is an analytic
space, with local analytic charts given by the $Y(U, \nu)$. It is clear from the
construction that any admissible normal function $\nu$ on $X$ with torsion
singularities extends to a holomorphic section of $\Jbtor(\shH)$: the extension is
given by the zero section of $Y(\Xb, \nu) \to \Xb$, followed by the inclusion into
the quotient.  \end{proof}

Let $\shG = \ker \bigl( R^1 \jl \shHZ \to R^1 \jl \shHQ \bigr)$ be the sheaf of
torsion sections in $R^1 \jl \shHZ$; if the singularity of an admissible normal
function on $X$ is torsion, then it is an element of $H^0 \bigl( \Xb, \shG \bigr)$.
Note that $\shG$ is a constructible sheaf of finite abelian groups, with support
contained in $\Xb - X$.

\begin{lemma} For a point $x \in \Xb$, let $G_x$ denote the stalk of the sheaf $\shG$
at $x$. Then every element of $G_x$ is the singularity of an admissible normal
function that is defined in a neighborhood of $x$.  \end{lemma}

\begin{proof} Fix an element $g \in G_x$. After replacing $\Xb$ by a small open
neighborhood of $x$, if necessary, we may assume that $g$ belongs to $H^1 \bigl( X,
\shHZ \bigr)$ and therefore corresponds to an extension of local systems
\begin{diagram} 0 &\rTo& \shHZ &\rTo& \shVZ &\rTo& \ZZ_X &\rTo& 0 \end{diagram} on
$X$. The extension splits over $\QQ$ because $g$ is torsion.  Since $\shVQ \simeq
\shHQ \oplus \QQ$, it follows that $\shVZ$ underlies the variation of mixed Hodge
structure $\shV = \shH \oplus \QQ(0)$. Now $\shV$ is clearly admissible, and
therefore corresponds to an admissible normal function, whose singularity equals the
original element $g \in H^1 \bigl( X, \shHZ \bigr)$.  \end{proof}

\subsection{Impossibility of a general analytic N\'eron model}
\label{subsec:impossible}

We now describe the implications of Theorem~\ref{thm:graphs} for the construction of
the full N\'eron model. As mentioned in the introduction, it should have the property
that its sections are the admissible normal functions.

\begin{lemma} \label{lem:impossible} Let $X \subseteq \Xb$ be a Zariski-open subset,
and let $\shH$ be a polarized variation of Hodge structure of weight $-1$ on $X$.
Suppose that there is a topological space $Y$ with the following three properties:
\begin{enumerate}[label=(\roman{*}), ref=(\roman{*})] \item The topology on $Y$ is
Hausdorff, and there is a continuous map $Y \to \Xb$.  \item There is a continuous
injective map $\Jb(\shH) \into Y$ over $\Xb$ that is a homeomorphism over $X$.  \item
Admissible normal functions on $X$ extend to a continuous sections of $Y$.
\end{enumerate} Then the the closure of the graph of an admissible normal function
inside $\Jb(\shH)$ can meet every fiber of $p \colon \Jb(\shH) \to \Xb$ in at most
one point.  \end{lemma}

\begin{proof} Let $\nu \colon X \to J(\shH)$ be an admissible normal function. By
assumption, it extends to a continuous section $\mu \colon \Xb \to Y$, and since $Y$
is Hausdorff, its graph $\mu(\Xb)$ has to be closed. It follows that the preimage of
$\mu(\Xb)$ in $\Jb(\shH)$ is also closed, and therefore contains the closure of
$\Gammanu$. But this implies that $\overline{\Gammanu}$ intersects each fiber
$p^{-1}(x)$ in at most one point.  \end{proof}

Now the problem is that, for a general admissible normal function with nontorsion
singularities, the closure of the graph typically has fibers of positive dimension
over $\Xb$. This can happen even in the simplest of examples: \subsecref{subsec:ex3}
exhibits a family of elliptic curves over $\dstn{2}$, pulled back from $\dst$ via the
map $(s_1, s_2) \mapsto s_1 s_2$, where the central fiber of $\Jb(\shH) \to \Delta^2$
is a copy of $\CC^{\ast}$. One can then easily write down an admissible normal
function on $\dstn{2}$ that extends holomorphically to $\Delta^2 - \{(0,0)\}$, but
such that the closure of its graph has a one-dimensional fiber over the origin. 

In my eyes, this example makes the existence of a N\'eron model that is Hausdorff as
a topological space very unlikely, for the following reason: For a family of elliptic
curves on $\dstn{2}$ with unipotent monodromy, any reasonable candidate for the
N\'eron model should have $\Jb(\shH)$ as its identity component, since the latter
agrees with the classical construction \cite{Namikawa}. By
Lemma~\ref{lem:impossible}, this means that the normal function in the example cannot
be a continuous section of a N\'eron model that is also Hausdorff. Thus it appears
that one cannot do any better than Theorem~\ref{thm:Neron-torsion}, at least if one
is interested in producing analytic spaces or Hausdorff spaces.

\subsection{Comparison with Brosnan-Pearlstein-Saito} \label{subsec:BPS}

We now make the comparison of our construction with the N\'eron model defined by
P.~Brosnan, G.~Pearlstein, and M.~Saito \cite{BPS}. We denote the identity component
of their model by $\JBPS(\shH)$.

We begin by constructing a map on fibers. Let $i \colon \{x\} \to \Xb$ be the
inclusion of an arbitrary point; following \subsecref{subsec:points}, define two
mixed Hodge structures $H = H^{-n} \iu M$ and $P = H^n \ius M(n) \simeq \Hd(1)$.
Lemma~\ref{lem:restriction-points} provides us with a surjection from $\Jb(\shH)_x$
to the generalized intermediate Jacobian $J(H)$. As explained in
\subsecref{subsec:Jacobians}, $J(H) \simeq \Ext_{\MHS}^1 \bigl( \ZZ(0), H \bigr)$ is
exactly the fiber of $\JBPS(\shH)$ over the point $x$. In this way, we obtain for
every point $x \in \Xb$ a surjective map of complex Lie groups \[ \Jb(\shH)_x \onto
\JBPS(\shH)_x.  \]

\begin{lemma} \label{lem:blowdown} The resulting map of sets $\pi \colon \Jb(\shH)
\to \JBPS(\shH)$ is continuous.  \end{lemma}

\begin{proof} Because of how the topology on $\JBPS(\shH)$ is defined in \cite{BPS},
and because of the functoriality of our construction, it suffices to prove the
statement in the case when $\Xb - X$ is a divisor with normal crossings and the local
monodromy of $\shHZ$ is unipotent. Let $\shHOe$ be Deligne's canonical extension of
$\shHO$; then $\shHOe \into \Mmod$. The Hodge bundles extend to locally free
subsheaves $F_p \shHOe = \shHOe \cap F_p \Mmod$. Let $E \to \Xb$ be the holomorphic
vector bundle defined by $\shHOe$, and $F_0 E \subseteq E$ the subbundle given by
$F_0 \shHOe$. We then have a holomorphic map $T(F_0 \Mmod) \to T(F_0 \shHOe) \simeq
(F_0 E)^{\ast}$. Using the polarization, we see that $(F_0 \shHOe)^{\vee} \simeq
\shHOe / F_0 \shHOe$; this means that we get a holomorphic map \[ T(F_0 \Mmod) \to E
/ F_0 E \] from the analytic space on the left to the vector bundle on the right.
Since the topology on $\JBPS(\shH)$ is induced from that on $E / F_0 E$, and topology
on $\Jb(\shH)$ from that on $T(F_0 \Mmod)$, the continuity of $\Jb(\shH) \to
\JBPS(\shH)$ is immediate.  \end{proof}

In the case $\dim \Xb = 1$, a very precise description of the map $\pi$ as a
composition of blowups has been given in \cite{SS}.

\begin{note} The map $\Jb(\shH)_x \to \JBPS(\shH)_x$ constructed in
\subsecref{subsec:points} has a splitting (and this is what guarantees its
surjectivity): in fact, we have maps \[ \HC / F_0 \HC \to T(F_0 \Mmod)_x \to (F_0
\PC)^{\vee}, \] and the first and last space are naturally isomorphic. This turns out
to be something of a red herring though, because the resulting map $\JBPS(\shH) \to
\Jb(\shH)$ is neither continuous, nor compatible with normal functions (as pointed
out by M.~Saito).  \end{note}

Now let $\nu$ be an admissible normal function on $X$ with locally trivial cohomology
class. We can also show that its extension $\nub$ to a holomorphic section of
$\Jb(\shH) \to \Xb$ is mapped to the extension constructed in \cite{BPS}.

\begin{lemma} Let $\nub \colon \Xb \to \Jb(\shH)$ be the extension of an admissible
normal function $\nu$ without singularities. Then the induced section $\pi \circ
\nub$ of $\JBPS(\shH)$ agrees with the extension of $\nu$ defined by P.~Brosnan,
G.~Pearlstein, and M.~Saito.  \end{lemma}

\begin{proof} Associated to the normal function, we have an extension of variations
of mixed Hodge structure \begin{diagram} 0 &\rTo& \shH &\rTo& \shV &\rTo& \ZZ_X(0)
&\rTo& 0.  \end{diagram} Because of admissibility, $\shV$ can be extended to a mixed
Hodge module $\Nmodnu$ on $\Xb$ with $W_{n-1} \Nmodnu \simeq M$ and $\Gr_n^W \Nmodnu
\simeq \QQ_{\Xb}^H \decal{n}$. 

Fix a point $i \colon \{x\} \to \Xb$, and let $H = H^{-n} \iu M$ and $V = H^{-n} \iu
\Nmodnu$. Also define $P = H^n \ius M(n) \simeq \Hd(1)$. Since the cohomology class
of $\nu$ is trivial near $x$, it is easy to see that we obtain an extension of mixed
Hodge structures \begin{equation} \label{eq:ext-MHS-2} \begin{diagram} 0 &\rTo& H
&\rTo& V &\rTo& \ZZ(0) &\rTo& 0, \end{diagram} \end{equation} and therefore a point
in $\Ext_{\MHS}^1 \bigl( \ZZ(0), H \bigr) \simeq J(H)$; it is the value of the
extended normal function in $\JBPS(\shH)_x$. According to
\subsecref{subsec:Jacobians}, this point is obtained by choosing a lifting $\vZ \in
\VZ$ for $1 \in \ZZ$, and restricting it to a linear operator on $F_0 \PC$. If we
take $v$ equal to the value at $x$ of a locally defined flat section of $\shVZ$
splitting the extension of local systems, then it follows that this prescription is
compatible with the definition of the extended normal function $\nub$ in
Proposition~\ref{prop:nub}. This means that $\pi \bigl( \nub(x) \bigr)$ gives the
same point in $J(H)$, as claimed.  \end{proof}

\begin{note} A shorter proof is the following: Both the extension of $\nu$
constructed in \cite{BPS} and $\pi \circ \nub$ are continuous sections of
$\JBPS(\shH)$. Since they agree over $X$, and since $X$ is dense in $\Xb$, it follows
that they agree everywhere.  \end{note}

\section{Examples} \label{sec:examples}

\subsection{Non-unipotent monodromy}

In this section, we describe a simple one-parameter family of elliptic curves in
which the local monodromy is not unipotent. This illustrates the difference between
$\Jb(\shH)$ and the identity component of the N\'eron model constructed in
\cite{BPS}. Let $E = \CC / (\ZZ + \tau \ZZ)$ be the elliptic curve with an
automorphism of order six; here $\tau = e^{i \pi/3}$, and the automorphism is given
by multiplication by $\tau$. Note that $\tau^2 = \tau - 1$.

We consider the trivial family $E \times \dst$, as well as its quotient by $\ZZ / 6
\ZZ$; a generator acts on $\Delta$ as multiplication by $\tau$, and on $E$ by the
automorphism. We denote the local system corresponding to the quotient by $\shH$; our
aim is to describe the structure of $\Jb(\shH) \to \Delta$.

We first work out the monodromy. Let $\alpha$ and $\beta$ be the standard basis for
$H_1(E, \ZZ)$; in the usual fundamental domain inside $\CC$, the cycle $\alpha$ is
the image of the line segment from $0$ to $1$, and the cycle $\beta$ that of the
segment from $0$ to $\tau$. Drawing a picture, it is clear that the automorphism acts
by \[ \alpha \mapsto \beta, \qquad \beta \mapsto \beta - \alpha.  \] Letting
$\alphast$ and $\betast$ denote the dual basis for $H^1(E, \ZZ)$, we also have \[
\alphast \mapsto -\betast, \qquad \betast \mapsto \alphast + \betast.  \] Thus the
monodromy operator $T$ is given by \[ T = \begin{pmatrix} 0 & 1 \\ -1 & 1
\end{pmatrix}, \] and one easily checks that it has eigenvalues $\tau$ and
$\bar{\tau} = - \tau^2$.  Also, $\det(T - \id) = 1$, and so the local system (over
$\ZZ$) has vanishing $H^0$ and $H^1$. It is clear from the construction that
$\alphast + \tau \betast$ is an eigenvector for $\tau$ (indeed, it restricts to a
holomorphic $1$-form on each fiber).

Let $s$ be the holomorphic coordinate on $\Delta$. For our construction, we need the
minimal extension of the flat vector bundle with monodromy $T$; according to
\cite{Saito-HM}, this is given by the Deligne lattice on which the residues of the
connection lie in $(-1, 0 \rbrack$. Thus the correct extension is given by $\shO e_1
\oplus \shO e_2$, with connection \[ \nabla e_1 = -e_1 \tensor \frac{ds}{6s}, \qquad
\nabla e_2 = -e_2 \tensor \frac{5ds}{6s}.  \] Let $\HH \to \dst$, with $s = e^{2 \pi
i z}$, be the universal covering space; on $\HH$, a flat section $\sigma(z)$ with
$\sigma(z+1) = \tau \sigma(z)$ is then found by solving $f'(z) - \pi i f(z) / 3 = 0$,
and so \[ \sigma(z) = e^{\pi i z / 3} e_1.  \] Neglecting constants, we have
$\sigma(z) = \alphast + \tau \betast$; thus $\omega = e_1$ is a section of $F^0$ of
the canonical extension (since it gives a holomorphic $1$-form on each fiber), and \[
\alphast + \tau \betast = e^{\pi i z / 3} \omega.  \]

Thus we see that \[ \int_{m \alpha + n \beta} \omega = (m + \tau n) \cdot e^{-\pi i z
/ 3} = (m + \tau n) \cdot e^{\pi y / 3} \cdot e^{-\pi i x / 3}, \] which goes to
infinity with $y$ unless $m = n = 0$. It follows that the closure of the family of
integral lattices inside the line bundle (dual to $\shO \omega$) only adds one point;
thus the fiber of the N\'eron model $\Jb(\shH)$ over $0 \in \Delta$ is a copy of
$\CC$. This is what it should be, given that we started from a family of elliptic
curves.

Next, we look at admissible normal functions and their extensions. By definition,
admissibility can be tested by pulling back along a branched cover ($s = t^6$ in our
case) to make the monodromy unipotent \cite{Saito-ANF}.  Thus we only need to
consider the family $E \times \dst$. Admissibility implies that the normal function
extends to a holomorphic map $\Delta \to E$. Lifting this to a map $g \colon \Delta
\to \CC$, we have \[ g(\tau t) - \tau g(t) \in \ZZ + \ZZ \tau, \] because the normal
function is pulled back from the original family. It is easy to see that we can
choose $g$ so that, in fact, $g(\tau t) = \tau g(t)$.  This choice of $g$ represents
the pullback of the extended normal function; its value over the origin is $g(0) =
0$, and so the pullback of any admissible normal function to the family $E \times
\Delta$ has to go through the origin in $\CC$. This is consistent with the N\'eron
model constructed by P.~Brosnan, G.~Pearlstein, and M.~Saito in \cite{BPS}: its fiber
over the origin is a single point, because the local system $\shH$ has no nontrivial
sections on $\dst$.

It should be noted, however, that there are no constraints on the graphs of normal
functions in our N\'eron model $\Jb(\shH) \to \Delta$. In fact, as shown in
Proposition~\ref{prop:sections}, any holomorphic section of $\Jb(\shH) \to \Delta$ is
an admissible normal function; the reason why the pullback of such a section to $E
\times \Delta$ has to pass through the origin is that the image of $\Jb(\shH)_0 \to
E$ is a point.

\subsection{A singular N\'eron model}

The example in this section was suggested by M.~Saito; it shows that the analytic
space $\Jb(\shH)$ can have singularities once $\dim X \geq 2$.

We let $\HZ = \ZZ^4$, with $\RR$-split mixed Hodge structure given by $I^{1,-1}
\oplus I^{-1,1} \oplus I^{0,-2} \oplus I^{-2,0}$. We let the alternating bilinear
form $Q$ be given by the matrix
\[
	S = \begin{pmatrix}
		0 & 0 & -1 & 0 \\
		0 & 0 & 0 & -1 \\
		1 & 0 & 0 & 0 \\
		0 & 1 & 0 & 0
	\end{pmatrix}
\]
and define nilpotent operators
\[
	N_1 = N_2 = \begin{pmatrix}
		0 & 0 & 1 & 0 \\
		0 & 0 & 0 & 1 \\
		0 & 0 & 0 & 0 \\
		0 & 0 & 0 & 0
	\end{pmatrix}.
\]
Let $\omega \in \CC$ have $\Im \omega \neq 0$. If we are satisfied with having the
mixed Hodge structure split over $\ZZ$, we may set
\[
	I^{1, -1} = \CC \begin{pmatrix}
		0 \\ 0 \\ 1 \\ \omega
	\end{pmatrix}, \quad
	I^{-1, 1} = \CC \begin{pmatrix}
		0 \\ 0 \\ 1 \\ \overline{\omega}
	\end{pmatrix}, \quad
	I^{0,-2} = \CC \begin{pmatrix}
		1 \\ \omega \\ 0 \\ 0
	\end{pmatrix}, \quad
	I^{0,-2} = \CC \begin{pmatrix}
		1 \\ \overline{\omega} \\ 0 \\ 0
	\end{pmatrix}.
\]
These data define an $\RR$-split nilpotent orbit on $\dstn{2}$, by the rule $(z_1,
z_2) \mapsto e^{z_1 N_1 + z_2 N_2} F$, where $F$ is given by the $I^{p,q}$.
Evidently, it is the pullback of a nilpotent orbit from $\dst$, by the map $(z_1,
z_2) \mapsto z_1 z_2$.

We now describe the sheaf $F_0 \Mmod$ and the analytic space $T = T(F_0 \Mmod)$ over
$\Delta^2$. Let the coordinates on $\Delta^2$ be $(s_1, s_2)$. The
Deligne extension is a trivial vector bundle of rank $4$, with Hodge filtration given
by the $I^{p,q}$. Thus $F_0 \Mmod$ is spanned by three sections,
\[
	e_0 = \begin{pmatrix}
		0 \\ 0 \\ 1 \\ \omega
	\end{pmatrix}, \quad
	e_1 = \frac{1}{s_1} \begin{pmatrix}
		1 \\ \omega \\ 0 \\ 0
	\end{pmatrix}, \quad
	e_2 = \frac{1}{s_2} \begin{pmatrix}
		1 \\ \omega \\ 0 \\ 0
	\end{pmatrix}.
\]
This gives a presentation for $F_0 \Mmod$ in the form
\begin{diagram}
	\shO &\rTo^{\left( \begin{smallmatrix} 0 \\ -s_1 \\ s_2 \end{smallmatrix} \right)}& 
		\shO^3 &\rTo& F_0 \Mmod &\rTo& 0,
\end{diagram}
and so $T(F_0 \Mmod)$ is the subset of $\Delta^2 \times \CC^3$ given by the equation
$s_1 v_1 = s_2 v_2$, using coordinates $(s_1, s_2, v_0, v_1, v_2)$. Thus $T$ is a
vector bundle of rank $2$ outside the origin, while the fiber over the origin is
$\CC^3$. Moreover, $T$ is clearly singular along the line $\CC (0,0,v_0,0,0)$.

Next, we look at the embedding of the set of integral points $\TZ$. Let $h \in \ZZ^4$
be any integral vector. We compute that
\begin{align*}
	Q \bigl( e_0, e^{-(z_1 N_1 + z_2 N_2)} h \bigr) 
		& = (z_1 + z_2) (h_3 + h_4 \omega) - (h_1 + h_2 \omega), \\
	Q \bigl( e_j, e^{-(z_1 N_1 + z_2 N_2)} h \bigr) 
		& = -\frac{h_3 + h_4 \omega}{s_j} \qquad \text{(for $j = 1, 2$)}.
\end{align*}
This means that $\TZ \subseteq T$ is the closure of the image of the map $\HH^2 \times \ZZ^4 \to
\Delta^2 \times \CC^3$, given by the formula
\[
	\left( e^{2 \pi i z_1}, e^{2 \pi i z_2}, 
		(z_1 + z_2) (h_3 + h_4 \omega) - (h_1 + h_2 \omega),
		- \frac{h_3 + h_4 \omega}{e^{2 \pi i z_1}},
		- \frac{h_3 + h_4 \omega}{e^{2 \pi i z_2}} 
	\right).
\]
Over $s_1 s_2 = 0$, the points in the closure are of the form $\bigl( s_1, s_2, -(h_1
+ h_2 \omega), 0, 0 \bigr)$. Let $J_0 = \CC / (\ZZ + \ZZ \omega)$ be the torus
corresponding to the monodromy-invariant part of the mixed Hodge structure. Then the
quotient $T / \TZ$ has the following structure: over $\dstn{2}$, the fibers are
the two-dimensional intermediate Jacobians; over $(0,0)$, the fiber is $J_0 \times
\CC^2$; over the remaining points with $s_1 s_2 = 0$, the fiber is $J_0 \times \CC$.
Moreover, $T / \TZ$ is singular along the torus $J_0 \times \{(0,0)\}$ over the origin.

\begin{note}
In this case, the Zucker extension is not Hausdorff. In fact, the integral
points are embedded into the ambient space $\Delta^2 \times \CC^2$ via the map $\HH^2
\times \ZZ^4 \to \Delta^2 \times \CC^2$, given by the formula
\[
	\bigl( e^{2 \pi i z_1}, e^{2 \pi i z_2}, 
		(z_1 + z_2) (h_3 + h_4 \omega) - (h_1 + h_2 \omega), 
		-(h_3 + h_4 \omega)
	\bigr).
\]
The closure of the image is much bigger than just the set of monodromy-invariant classes in
$\HZ$; to obtain the Zucker extension, therefore, one is taking a quotient by a
non-closed equivalence relation, which can never produce a Hausdorff space.
\end{note}

\subsection{A normal function with nontorsion singularity}
\label{subsec:ex3}

In this section, we shall look at a simple example of a normal function on $\dstn{2}$
with a nontorsion singularity at the origin in $\Delta^2$. The interesting point here
is that the closure of its graph has a one-dimensional fiber over the origin.

The example is a family of elliptic curves; the corresponding variation
of Hodge structure of weight $-1$ is a nilpotent orbit, which we describe by giving
its limit mixed Hodge structure. So let $\HZ = \ZZ^2$, with nilpotent operators
\[
	N_1 = N_2 = \begin{pmatrix} 0 & 1 \\ 0 & 0 \end{pmatrix},
\]
and define the limit mixed Hodge structure by letting $I^{0,0} = \CC (0,1)$ and
$I^{-1,-1} = \CC (1,0)$. The period mapping of the associated variation of Hodge
structure is then given by $\Phit(z) = e^{z_1 N_1 + z_2 N_2} F$, and so the vector
$(z_1+z_2,1)$ spans $\Phit(z)^0$.

We now introduce an admissible normal function through its variation of mixed Hodge
structure. Let $\VZ = \HZ \oplus \ZZ$, and define
\[
	N_1' = \begin{pmatrix} 0 & 1 & 1 \\ 0 & 0 & 0 \\ 0 & 0 & 0 \end{pmatrix}
	\quad \text{and} \quad
	N_2' = \begin{pmatrix} 0 &1 & -1 \\ 0 & 0 & 0 \\ 0 & 0 & 0 \end{pmatrix};
\]
Thus the vector $v = (0,0,1)$ belongs to $\VZone$, and satisfies $N_1' v + N_2' v =
0$. Let $W_{-1} = \HZ$ and $W_0 = \VZ$. The $\RR$-split mixed Hodge structure $(M,F)$
with $I^{0,0}(M,F) = \CC (0,1,0) \oplus \CC(\lambda,0,1)$ and $I^{-1,-1}(M,F) =
\CC(1,0,0)$ defines a mixed nilpotent orbit $\bigl( W, e^{z_1 N_1' + z_2 N_2'} F
\bigr)$, and one can easily check that it is admissible. Let $\nu$ denote the
corresponding admissible normal function on $\dstn{2}$.

We will now determine the closure of $\Tnu$ inside $T(F_0 \Mmod)$. In this situation,
$F_0 \Mmod$ is a trivial line bundle on $\Delta^2$, whose pullback to $\HH^2$ is
spanned by the section $(z_1+z_2,1)$. The map $\Tnu \into T(F_0 \Mmod)$ now takes the
form
\[
	\HH^2 \times \VZone \to \Delta^2 \times \CC,
\]
and is given by the formula
\[
	(z_1,z_2, a,b,1) \mapsto \bigl( e^{2\pi i z_1}, e^{2 \pi i z_2}, 
		a - b(z_1 + z_2) + (z_2 - z_1) \bigr).
\]
From this, it is easy to determine the closure of the graph. Over a point $(s_1,0)$
with $s_1 \neq 0$, we only get points in the
closure when $b = 1$, and so the fiber consists of all points $a -2 z_1$ with $e^{2 \pi
i z_1} = s_1$. Similarly, the fiber over $(0,s_2)$ with $s_2 \neq 0$ is the discrete
set of points $a + 2 z_2$ with $e^{2 \pi i z_2} = s_2$. More interesting is the fiber
over $(0,0) \in \Delta^2$. By taking $a = b = 0$ and $z_2 = z_1 + w$ with $w \in \CC$
arbitrary and $\Im z_1 \to \infty$, we see that the fiber consists of all of $\CC$.

The quotient $\Jb = T(F_0 \Mmod) / \TZ$ is a family of elliptic curves over
$\dstn{2}$, with fibers over $s_1 s_2 = 0$ copies of $\CC^{\ast}$. The discussion
above shows that $\nu$ extends to an admissible normal function over $\Delta^2 -
\{(0,0)\}$, but that the closure of the graph of $\nu$ inside $\Jb$ contains the
entire fiber $\CC^{\ast}$ over $(0,0)$. As mentioned in \subsecref{subsec:graphs}, this is
evidence that there can probably not exist a N\'eron model that is Hausdorff as a topological
space.

%---------------------------------------------

\section*{References}

%\begin{bibsection}

%\end{bibsection}

\end{document}